\documentclass[12pt, a4paper]{article}
\usepackage{amsmath}
\usepackage{amssymb}
\usepackage{amsthm}
\usepackage[totalwidth=17cm, totalheight=25cm]{geometry}
\usepackage{graphicx}
\usepackage{picins}
\usepackage{tabularx}
	\newcolumntype{C}[1]{>{\centering\arraybackslash}m{#1}} 
	\newcolumntype{R}[1]{>{\raggedleft\arraybackslash}m{#1}} 

\newtheoremstyle{boldplain}
{9pt}
{9pt}
{\itshape}
{}
{\bfseries}
{.}
{.5em}
{\thmname{#1}\thmnumber{ #2}\thmnote{ (#3)}}%

\newtheoremstyle{bolddefinition}
{9pt}
{9pt}
{}
{}
{\bfseries}
{.}
{.5em}
{\thmname{#1}\thmnumber{ #2}\thmnote{ (#3)}}%

\theoremstyle{boldplain}
\newtheorem{conj}[equation]{Conjecture}
\newtheorem{cor}[equation]{Corollary}
\newtheorem{lem}[equation]{Lemma}
\newtheorem{prop}[equation]{Proposition}

\newtheorem{sublem}[equation]{Sublemma}
\newtheorem{thm}[equation]{Theorem}

\theoremstyle{bolddefinition}
\newtheorem{rem}[equation]{Remark}

\newfont{\bigbf}{cmbx10 scaled\magstep1}
\normalbaselines
\numberwithin{equation}{section}

\def\no{\noindent}

\def\R{{\mathbb R}}

\def\Si{\Sigma}

\def\acts{\curvearrowright}

\def\half{\frac{1}{2}}
\def\third{\frac{1}{3}}

\def\ora{\overrightarrow}
\def\pihalf{\frac{\pi}{2}}
\def\pithird{\frac{\pi}{3}}
\def\2pithird{\frac{2\pi}{3}}
\def\piquart{\frac{\pi}{4}}
\def\3piquart{\frac{3\pi}{4}}
\def\quart{\frac{1}{4}}

\def\phm{\phantom{-}}
\def\aast{\ast\ast}
\def\ac{\arccos}
\def \CAT {\text{CAT}}
\def\wht{\widehat}

\title{The Center Conjecture for thick spherical buildings} 
\author{Carlos Ramos-Cuevas\footnote{
Mathematisches Institut, Ludwig-Maximilians-Universit\"at,
Theresienstr. 39,
D-80333 M\"unchen, Germany,
cramos@mathematik.uni-muenchen.de}}
\date{October 17, 2012}
\begin{document}
 \maketitle

\begin{abstract}
In this paper we show that a convex subcomplex of a 
spherical building
of type $E_6$, $E_7$ or $E_8$ is a subbuilding or the 
automorphisms of the subcomplex fix
a point on it. Together with previous results 
of M\"uhlherr-Tits, and Leeb and the author, this completes
the proof of Tits' Center Conjecture for 
spherical buildings without factors of type $H_4$, 
in particular,
for thick spherical buildings.
\end{abstract}

\section{Introduction}

The Center Conjecture was first proposed by J.\ Tits in the 50's and 
it is now formulated as follows
(compare \cite{MuehlherrTits} and \cite[Conjecture 2.8]{Serre}).

\begin{conj}[Center Conjecture]\label{centerconj}
Suppose that $B$ is a spherical building 
and that $K\subseteq B$ is a convex subcomplex. 
Then $K$ is a subbuilding 
or the action $Stab_{Aut(B)}(K)\acts K$ 
of the automorphisms of $B$ preserving $K$ has a fixed point. 
\end{conj}

A building automorphism is an isometry, which preserves
the polyhedral structure of the building.
A fixed point of the action $Stab_{Aut(B)}(K)\acts K$
is called a {\em center} of the subcomplex $K$.

Part of Tits' motivation for 
the Center Conjecture
came from algebraic group theory
(cf. \cite[Lemma 1.2]{Titscenterconj}).
A special case of the Center Conjecture 
is also considered in Geometric Invariant Theory
related to the discussion of instability
(see \cite{Mumford}). 
This special case 
was proven by Rousseau \cite{Rousseau} and Kempf \cite{Kempf}.
The conjecture regained interest through the work of Serre 
on {\em complete reducible} subgroups of algebraic groups (see \cite{Serre}).
A more general geometric form of the Center Conjecture 
concerning arbitrary convex {\em subsets} and their circumradii
arises when studying the asymptotic
geometry of convex sets in symmetric spaces of noncompact type 
(cf. \cite{KleinerLeebinvconv}). 
We refer to the introduction of \cite{LeebRamos}
for a detailed discussion of several aspects concerning the Center
Conjecture and other related questions.

Recently, after the publication of a preprint version of this paper,
applications of the Center Conjecture
have appeared in e.g.\ \cite{BateMartinRoehrle3}, \cite{Struyve}.

If $dim(K) \leq 1$ then it is easy to see that the Center Conjecture
holds, namely, a one-dimensional convex subcomplex 
is a building or a tree with circumradius $\leq\pihalf$. In the last case, it has a 
unique circumcenter which is fixed by $Isom(K)$. 
If $dim(K)\leq 2$ a stronger version of the conjecture has
been shown in \cite{BalserLytchak}.

If $K$ is a convex subcomplex of a reducible building
$B=B_1\circ\dots\circ B_k$, then $K$ also decomposes as a spherical join
$K=K_1\circ\dots\circ K_k$ where each $K_i\subset B_i$ is a convex subcomplex
for $i=1,\dots,k$.
Thus, the Center Conjecture easily reduces to the case of irreducible buildings.
For irreducible buildings of classical type
(i.e.\ $A_n$, $B_n$ and $D_n$) the Center Conjecture was shown by M\"uhlherr and Tits 
\cite{MuehlherrTits}. 
The $F_4$-case was announced 
by Parker and Tent in a talk in Oberwolfach \cite{ParkerTent_conv}.
In both cases the arguments are based on the incidence-geometric
realizations of the corresponding different types of buildings.

Our approach to this problem
is of differential-geometric nature, using methods
of the theory of metric spaces with curvature bounded above.
It develops further the techniques in \cite{LeebRamos} where
the cases of buildings of type $F_4$ and $E_6$
were settled.

The main result in this paper is:

\begin{thm}
 The Center Conjecture~\ref{centerconj} holds for spherical buildings
of type $E_6$, $E_7$ and $E_8$.
\end{thm}

The proof of the $E_7$-case in this paper uses 
the result for the $E_8$-case. A direct proof can
be found in \cite{Ramos}. 
The $E_6$-case is also deduced here from the $E_8$-case, 
however, this case is less complex and the direct
proof given in \cite{LeebRamos} is much shorter.
The case of buildings of type $H_3$ is a consequence
of the main result in \cite{BalserLytchak}. It can also 
be easily treated with our methods (see \cite{Ramos}).
Hence we obtain the following result.

\begin{cor}
The Center Conjecture~\ref{centerconj} holds for
spherical buildings without factors of type $H_4$.
\end{cor}

While any spherical Coxeter complex is a spherical building,
not all spherical Coxeter complexes occur
as Coxeter complexes for {\em thick} spherical buildings
(\cite{Tits_w}). Namely,
there are no thick spherical buildings of type
$H_3$ (\includegraphics[scale=0.4]{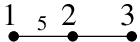}) and
$H_4$ (\includegraphics[scale=0.4]{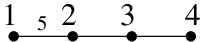}),
these being the only cases.
On the other hand,
any spherical building has a canonical 
thick structure (depending only on its isometry type)
which results from
restricting to a subgroup of the Weyl group
(\cite{Scharlau}, \cite[Sec.\ 3.7]{KleinerLeeb}).
The polyhedral structure thus obtained is (possibly) coarser.
The Center Conjecture is more natural when posed
for thick spherical buildings, because then
$K$ is a subcomplex of the natural 
polyhedral structure of $B$.
In this case we have:

\begin{cor}
 The Center Conjecture~\ref{centerconj} holds for
 all thick spherical buildings.
\end{cor}

A completely different approach 
to the special case of the Center Conjecture
for spherical buildings $B$ associated to
algebraic groups $G$ and 
subcomplexes $K$ which are fixed point sets of the action
of a subgroup $H\subset G$ can be found in \cite{BateMartinRoehrle}.
They show that such a subcomplex is a subbuilding or the action
$Stab_G(K)\acts K$ fixes a point. In \cite{BateMartinRoehrleTange}
this result is extended to the action $Stab_{Aut(G)}(K)\acts K$.
In \cite{BateMartinRoehrle2} a special kind of convex
subsets of such fixed point subcomplexes is also considered.

Our arguments actually yield a more general and intrinsic version
of the Center Conjecture (cf. \cite[Thm.\ 3.24]{LeebRamos}).

\begin{cor}
 If $B$ is a spherical building without factors of type $H_4$ 
and $K\subseteq B$ is a convex subcomplex, then $K$
is a subbuilding or the action $Aut_B(K)\acts K$ has a fixed point.
\end{cor}

We define the automorphisms in $Aut_B(K)$ to be isometries of $K$
preserving its polyhedral structure induced by $B$
and such that the changes of types of vertices of $K$ 
correspond to a symmetry of the Dynkin diagram of $B$.
They need not
be extendable to automorphisms of $B$
(see Remark \ref{rem:ccE8}).
Hence, we have $Stab_{Aut(B)}(K)\subset Aut_B(K)\subset Isom(K)$.

Recently an alternative proof of Corollary~\ref{lem:no8Apt} has been given
in \cite{MuehlherrWeiss} in the case of {\em thick} spherical buildings. 
They prove the analogous statement of Corollary~\ref{lem:no8Apt}
for thick irreducible buildings of any 
type, thus, giving a proof of the Center Conjecture 
in the case of top-dimensional
convex subcomplexes of thick spherical buildings.

\tableofcontents

\section{Preliminaries}

In this section we will fix the terminology
that we will use through this paper.
A more detailed treatment of the basic geometric notions
and facts concerning $\CAT(1)$ spaces, spherical Coxeter
complexes and spherical buildings
needed for our proof of the Center Conjecture
can be found in \cite[Sec. 2]{LeebRamos} (see also \cite{KleinerLeeb}).

\subsection{$\CAT(1)$ spaces}

Recall that a complete metric space is called 
a {\em $\CAT(1)$ space} 
if it is $\pi$-geodesic 
and geodesic triangles of perimeter less than $2\pi$
are not {\em thicker} than those in the
round sphere with curvature $\equiv 1$.

For points $x,y$ in a $\CAT(1)$ space $X$ at distance $<\pi$, 
we denote by $xy$ the unique segment connecting both points. 
Let $m(x,y)$ denote the midpoint of the segment $xy$.
Two points at distance $\geq\pi$ are called {\em antipodal}.

The {\em link} $\Si_x X$ at a point $x\in X$ is the space of directions at $x$
with the angle metric. It is again a $\CAT(1)$ space. If $y\neq x$ 
and $y$ is not antipodal to $x$, we denote
with $\ora{xy}\in\Si_x X$ the direction at $x$ of the segment $xy$. 

A subset $C$ of a $\CAT(1)$ space is called {\em convex}, 
if for any 
$x,y\in C$ at distance $<\pi$ the segment $xy$ 
is also contained in $C$.
A closed convex subset of a $\CAT(1)$ space is itself $\CAT(1)$.

Let $A$ be a subset of a $\CAT(1)$ space $X$ and let $x\in X$. 
The {\em radius of $A$ with respect to $x$} is defined as 
$rad(x,A):=\sup\{d(x,y) | y \in A\}$ and the {\em circumradius of $A$
in $X$} is 
$rad_X(A)=\inf\{rad(x,A) | x\in X\}$.
A point $x\in X$, such that $rad(x,A)=rad_X(A)$ is called a
{\em circumcenter}.

For more information on $\CAT(1)$ spaces we refer to \cite{BridsonHaefliger}.

\subsection{Coxeter complexes}

A spherical Coxeter complex $(S,W)$ is a pair consisting in
a round sphere $S$ with curvature $\equiv 1$
together with a finite group of isometries $W$, called
the {\em Weyl group}, generated by reflections on
great spheres of codimension 1. 

There is a natural structure of spherical 
polyhedral complex on $S$ induced by $W$.
The spheres of codimension 1, that are the fixed point sets of the reflections in $W$ 
are called the {\em walls}.
The {\em Weyl chambers} are the closures of the connected components of 
$S$ minus the union of all the walls.
A Weyl chamber is a convex spherical polyhedron.
The Weyl chambers are fundamental domains for the action
of the Weyl group on $S$ and therefore isometric to the
{\em model Weyl chamber} $\triangle_{mod}:=S/W$.
A {\em root} is a top-dimensional hemisphere bounded by a wall.
A {\em singular} sphere is an intersection of walls.
A segment contained in a singular 1-sphere is called a 
{\em singular segment}.

The geometry of a spherical Coxeter complex can be encoded in a graph, the so-called
{\em Dynkin diagram} (cf. \cite[Ch. 5.1]{GroveBenson}). 
A labelling by an index set $I$
of the vertices of the Dynkin diagram induces a labelling of
the vertices of $S$. We say that a vertex in $S$ is of
{\em type} $i$ or that it is an $i$-vertex for $i\in I$,
if it has label $i$.

An {\em automorphism} of $(S,W)$ is an isometry of $S$ preserving its polyhedral
structure. 
The group $Aut(S,W)/W$ can be canonically identified with the isometries of
the model Weyl chamber $\triangle_{mod}$
and therefore also with the symmetries of the Dynkin diagram.
Notice that the antipodal involution of $S$ is always an automorphism of
$(S,W)$. The {\em canonical involution} of $\triangle_{mod}$ is the image
of the antipodal involution under the identification mentioned above. 

We refer to \cite{GroveBenson} and \cite[Sec. 3.1, 3.3]{KleinerLeeb} for further
information on spherical Coxeter complexes.

\subsection{Spherical buildings}

We refer to \cite{AbramenkoBrown}, \cite{KleinerLeeb} and \cite{Tits_bn}
for more information on spherical buildings. We will consider spherical
buildings from the point of view of $\CAT(1)$ spaces as presented in
\cite{KleinerLeeb}.

A spherical building $B$ modelled on a spherical
Coxeter complex $(S,W)$ is a $\CAT(1)$ space
together with an atlas $\cal A$ of isometric embeddings $S \hookrightarrow B$
(the images of these embeddings are called {\em apartments})
with the following properties: any two points in $B$ are contained in a common
apartment, the atlas $\cal A$ is closed under precomposition with isometries in $W$
and the coordinate changes are restrictions of isometries in $W$. The empty set
is considered to be a building.

The polyhedral structure of $(S,W)$ induces a 
polyhedral structure on the building $B$. 
A labelling by an index set $I$
of the vertices of the Dynkin diagram of $(S,W)$ induces a labelling of
the vertices of $B$.
The objects (walls, roots,... ) defined for spherical Coxeter complexes
can be defined for the building $B$ as the corresponding images in $B$. 

For any point $x\in B$, the link $\Si_x B$
is again a spherical building (cf. \cite[Prop. 3.6.4]{KleinerLeeb}).
Its Dynkin diagram its obtained from the one of $B$ by
deleting those vertices corresponding to the vertices of the face of $B$
spanned by $x$.

A building is called {\em thick} if every wall is the boundary of at least 
three different roots.

A {\em subbuilding} is a convex subcomplex $K$ of a building, such that any
two points in $K$ are contained in a singular sphere 
$s\subset K$ of the same dimension as $K$.
A subbuilding carries a natural structure as a spherical building
induced by its ambient building 
(cf. \cite[Proposition 2.13]{LeebRamos}).

\section{Spherical Coxeter complexes}\label{sec:coxeter}

This section contains some geometric properties of spherical Coxeter complexes.

In our arguments later, we will need some information on singular spheres
of codimension $\leq 2$ in the different Coxeter complexes.

If the Coxeter complex  $(S,W)$ is irreducible and its Dynkin diagram
has no weights on its edges, i.e.\ if it is of type
$A_n$, $D_n$, $E_6$, $E_7$ or $E_8$, then it is easy to see, that
the Weyl group acts transitively on the set of roots 
(\cite[Proposition 5.4.2]{GroveBenson}).
In particular all walls (singular spheres of codimension 1) are equivalent modulo the
action of $W$.
If there is more than one orbit of roots, then we define the {\em type} of a wall 
as the type of the center of a root bounded by this wall. Note that this
definition is independent of which of both roots we take. 
A vertex is of {\em root type} if it is the center of a root.

A singular sphere of codimension 2 is the intersection of two different walls.
We define the {\em type} of a sphere of codimension 2 as the type of the circle 
spanned by the centers of the corresponding roots.

We gather in the next sections some of the geometric properties of
the different Coxeter complexes.
This information can be deduced from the data in the Appendix~\ref{app:coxeter}
(see Remarks~\ref{rem:1-2-spheres}-\ref{rem:tables} 
for descriptions of how to do the respective computations). 

\subsection{The Coxeter complex of type $D_n$}\label{sec:Dncoxeter}

For $n\geq 4$ let $(S,W_{D_n})$ be the spherical Coxeter complex of type $D_n$
with Dynkin diagram \hpic{\includegraphics[scale=0.4]{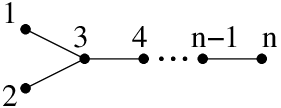}}.
It has dimension $n-1$.

The $(n-1)$-vertices are the vertices of {\em root type}.
All hemispheres bounded by walls are centered at a $(n-1)$-vertex
(see Section~\ref{app:Dn}).

For $n\geq 5$ the Dynkin diagram has one symmetry: it exchanges the vertices
$1 \leftrightarrow 2$ and fixes the others.
This symmetry is induced by the canonical involution 
of the Weyl chamber $\triangle^{D_n}_{mod}$
if $n$ is odd.
If $n$ is even, then the canonical involution is trivial.
For $n=4$ the Dynkin diagram has six symmetries, they permute the
vertices 1, 2, 4 and fix the vertex 3.

We describe now the possible lengths and types 
(modulo the action of the Weyl group) of segments
between vertices. 
In Remarks~\ref{rem:segmenttype}, \ref{rem:tables}
we explain how these lengths can be 
deduced from the data in the Appendix~\ref{app:Dn}.
The simplicial convex hulls of the segments can be determined
by identifying the smallest singular sphere containing the circle
spanned by the corresponding two vertices.
We list only the lengths and types that we will need later.

\begin{itemize}
\item Distances between two $n$-vertices $x$ and $x'$:
\begin{center}
\begin{tabular}{|C{4cm}|C{10cm}|}
\hline
Distance & Simplicial convex hull  of segments $xx'$ \\ \hline
$0,\pi$&\\ \hline
$\pihalf$ & singular segment of type $n(n-1)n$ \\ \hline
\end{tabular}
\end{center}

\item Distances between two $(n-1)$-vertices $x$ and $x'$:
\begin{center}
\begin{tabular}{|C{4cm}|C{10cm}|}
\hline
Distance & Simplicial convex hull  of segments $xx'$ \\ \hline
$0,\pi$&\\ \hline
$\pihalf$ & singular segment of type $(n-1)n(n-1)$ for $n\geq 4$/\hfill\linebreak
		singular segment of type $(n-1)(n-3)(n-1)$, if $n\geq6$;\hfill\linebreak
		\hpic{\includegraphics[scale=0.35]{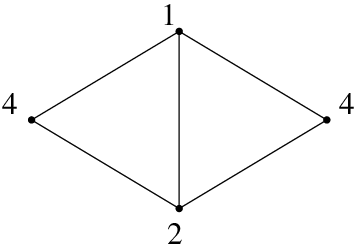}}, 
			if $n=5$. \hfill\linebreak
		singular segment of type 313 or 323, if $n=4$.\\ \hline
$\pithird$ ($\2pithird$) & 
		\hpic{\includegraphics[scale=0.55]{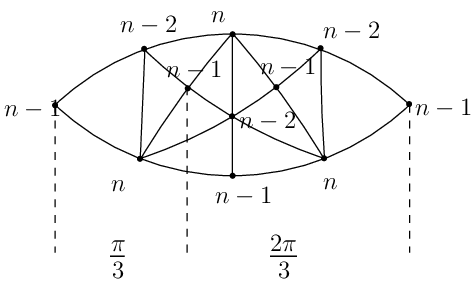}}, 
			if $n\geq5$;\hfill\linebreak
		if n=4, the simplicial convex hull of a segment $xx'$ 
			of length $\pithird$ is 3-dimensional:\hfill\linebreak
			\hpic{\includegraphics[scale=0.55]{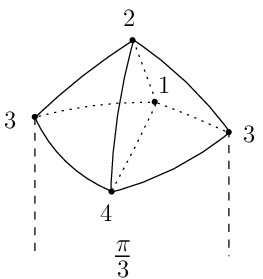}}\hfill\linebreak 
			A segment $xx'$ of length $\2pithird$ consists of two
			segments of length $\pithird$ as above.
			\\ \hline
\end{tabular}
\end{center}

\item Distances between two 1- (2)-vertices $x$ and $x'$:
\begin{center}
\begin{tabular}{|C{4cm}|C{10cm}|}
\hline
Distance & Simplicial convex hull  of segments $xx'$ \\ \hline
$\ac(\frac{n-4k}{n})$ for $k=0,1,\dots,[\frac{n}{2}]$
		& singular segment of type $1(2k+1)1$, ($2(2k+1)2$ resp.) \\ \hline
\end{tabular}
\end{center}

\item Distances between a 1- (2)-vertex $x$ and a $n$-vertex $y$:
\begin{center}
\begin{tabular}{|C{4cm}|C{10cm}|}
\hline
Distance & Simplicial convex hull  of segments $xy$ \\ \hline
$\ac(\frac{1}{\sqrt{n}})$ & singular segment of type $1n$, ($2n$ resp.) \\ \hline
$\ac(-\frac{1}{\sqrt{n}})$ & singular segment of type $12n$, ($21n$ resp.) \\ \hline
\end{tabular}
\end{center}
\end{itemize}

The following properties of singular spheres in $D_n$,
while not obvious,
can be easily seen in the vector space realization of
the Coxeter complex presented in Appendix~\ref{app:coxeter}.
We omit here their proofs.
\label{pag:spheresDn}

A wall in $S$ contains a singular sphere of codimension 1 spanned by $n-2$
pairwise orthogonal $n$-vertices, that is, the convex hull of these $n-2$
pairwise orthogonal $n$-vertices and their antipodes.

The convex hull of $n-1$ pairwise orthogonal $n$-vertices and their antipodes is
a $(n-2)$-sphere, but it is not a wall, in particular,
it is not a subcomplex. Its simplicial convex hull is $S$.

If $n\geq 5$ ($n=4$) there are three (four) types of singular spheres of
codimension 2. They correspond to the two (three) types of segments connecting
two $(n-1)$-vertices at distance $\pihalf$ and the unique type of segments
connecting two $(n-1)$-vertices at distance $\pithird$.
We say that a sphere of the last type is a $(n-3)$-sphere {\em of type $\pithird$}.

A singular sphere of codimension 2 always contains a singular 
$(n-5)$-sphere spanned by $n-4$ pairwise orthogonal $n$-vertices.

Let $h$ be a singular hemisphere of codimension 1 bounded 
by a singular $(n-3)$-sphere $s$.
It is the intersection of a wall and a root bounded by
a different wall.
If $n\geq 6$ and $s$ is of type
$(n-1)n(n-1)$ (or $(n-1)(n-3)(n-1)$),
 then $h$ is the join of $s$ with a point in the singular circle of type
$\dots(n-1)n(n-1)\dots$ (or $\dots(n-1)(n-3)(n-1)\dots$),
since $h$ is singular, this point is of type $n-1$. 
That is, $h$ is centered at a $(n-1)$-vertex $x$. The link $\Si_x h$ 
in the Coxeter complex $\Si_x S$ of type $D_{n-2}\circ A_1$ is a 
wall of type $n$ (or $(n-3)$).
If $n\geq 5$ and $s$ is of type $\pithird$, 
then $h$ is centered at a point contained
in a singular segment of type $n(n-2)$, it 
is the midpoint of a segment connecting
two $(n-1)$-vertices at distance $\pithird$.

\subsection{The Coxeter complex of type $E_6$}\label{sec:E6coxeter}

Let $(S,W_{E_6})$ be the spherical Coxeter complex of type $E_6$
with Dynkin diagram \hpic{\includegraphics[scale=0.4]{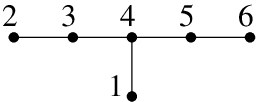}}.
It has dimension 5.

The 1-vertices are the vertices of {\em root type}. All hemispheres
bounded by walls are centered at a 1-vertex.

The Dynkin diagram has one symmetry, namely, the one that exchanges
 the vertices
$2 \leftrightarrow 6$, $3\leftrightarrow 5$ and fixes the 1-
and 4-vertices.
It corresponds to the canonical involution of the Weyl chamber
$\triangle^{E_6}_{mod}$. Therefore, the properties of $i$-
and $(8-i)$-vertices for $i=2,3,5,6$, are {\em dual} to each other. 

We describe now the possible lengths and types 
(modulo the action of the Weyl group) of segments
between vertices. 
In Remarks~\ref{rem:segmenttype}, \ref{rem:tables}
we explain how these lengths can be 
deduced from the data in the Appendix~\ref{app:E6}.
The simplicial convex hulls of the segments can be determined
by identifying the smallest singular sphere containing the circle
spanned by the corresponding two vertices.
Here this singular sphere is always the circle itself. 
See Remark~\ref{rem:1-2-spheres} for an explanation of how to determine
the singular circles.
We list only the lengths and types that we will need later.

\begin{itemize}
\item Distances between two 2- (6)-vertices $x$ and $x'$:
\begin{center}
\begin{tabular}{|C{4cm}|C{10cm}|}
\hline
Distance & Simplicial convex hull  of segments $xx'$ \\ \hline
$0$&\\ \hline
$\ac(\quart)$& singular segment of type 232 (656) \\ \hline
$\2pithird$ & singular segment of type 262 (626) \\ \hline
\end{tabular}
\end{center}
\item Distances between a 2-vertex $x$ and a 6-vertex $y$:
\begin{center}
\begin{tabular}{|C{4cm}|C{10cm}|} \hline
Distance & Simplicial convex hull  of segments $xy$ \\ \hline
$\pi$&\\ \hline
$\ac(-\quart)$& singular segment of type 216 \\ \hline
$\pithird$ & singular segment of type 26 \\ \hline
\end{tabular}
\end{center}

\end{itemize}

\subsection{The Coxeter complex of type $E_7$}\label{sec:E7coxeter}

Let $(S,W_{E_7})$ be the spherical Coxeter complex of type $E_7$
with Dynkin diagram \hpic{\includegraphics[scale=0.4]{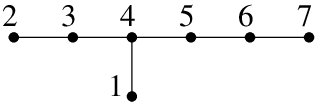}}.
It has dimension 6. 

The Dynkin diagram for $E_7$ has no symmetries, therefore all
automorphisms of $(S,W_{E_7})$ are {\em type preserving}.

These are the one dimensional singular spheres
(see Remark~\ref{rem:1-2-spheres}) in $(S,W_{E_7})$:
\begin{center}
\vspace{-0.3cm}
 \includegraphics[scale=0.65]{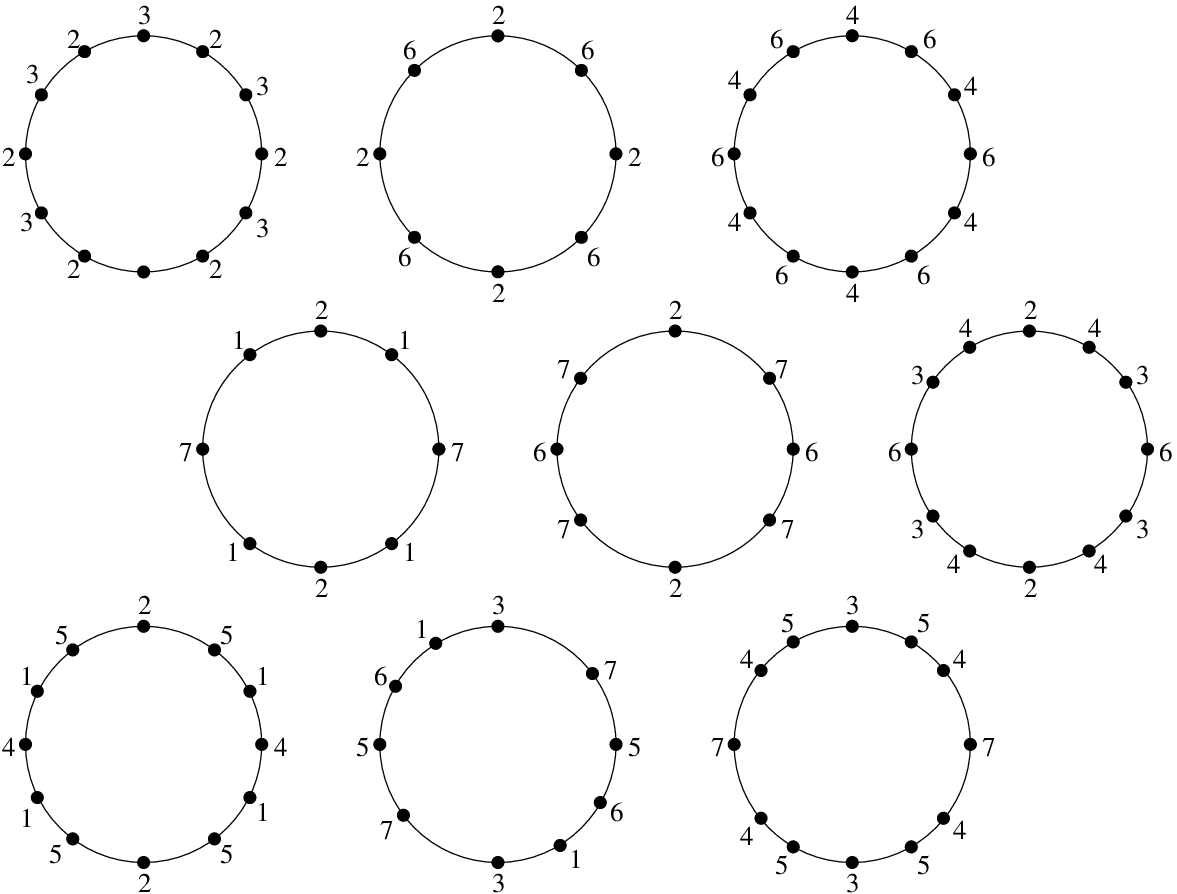}
\end{center}

The 2-vertices are the vertices of {\em root type}. All hemispheres
bounded by walls are centered at a 2-vertex.

We describe now the possible lengths and types 
(modulo the action of the Weyl group) of segments
between vertices. 
In Remarks~\ref{rem:segmenttype}, \ref{rem:tables}
we explain how these lengths can be 
deduced from the data in the Appendix~\ref{app:E7}.
The simplicial convex hulls of the segments can be determined
by identifying the smallest singular sphere containing the circle
spanned by the corresponding two vertices.
Here this singular sphere is always the circle itself. 
See Remark~\ref{rem:1-2-spheres} for an explanation of how to determine
the singular circles.
We list only the lengths and types that we will need later.

\begin{itemize}
\item 
Distances between two 2-vertices $x$ and $x'$:
\begin{center}
\begin{tabular}{|C{4cm}|C{10cm}|}
\hline
Distance & Simplicial convex hull  of segments $xx'$ \\ \hline
$0,\pi$&\\ \hline
$\pithird$& singular segment of type 232  \\ \hline
$\pihalf$& singular segment of type 262  \\ \hline
$\2pithird$ & singular segment of type 23232 \\
\hline
\end{tabular}
\end{center}
\item
Distances between two 7-vertices $x$ and $x'$:
\begin{center}
\begin{tabular}{|C{4cm}|C{10cm}|}
\hline
Distance & Simplicial convex hull  of segments $xx'$ \\ \hline
$0,\pi$&\\ \hline
$\arccos (\third )$& singular segment of type 767  \\ \hline
$\arccos (-\third )$& singular segment of type 727  \\ \hline
\end{tabular}
\end{center}

\item 
Distances between a 2-vertex $x$ and a 7-vertex $y$:
\begin{center}
\begin{tabular}{|C{4cm}|C{10cm}|}
\hline
Distance & Simplicial convex hull  of segments $xy$ \\ \hline
$\arccos (\frac{1}{\sqrt{3}})$& singular segment of type 27  \\ \hline
$\pihalf$& singular segment of type 217  \\ \hline
$\arccos (-\frac{1}{\sqrt{3}})$& singular segment of type 2767  \\ \hline
\end{tabular}
\end{center}
\end{itemize}

\subsection{The Coxeter complex of type $E_8$}\label{sec:E8coxeter}

Let $(S,W_{E_8})$ be the spherical Coxeter complex of type $E_8$
with Dynkin diagram \hpic{\includegraphics[scale=0.4]{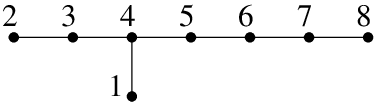}}.
It has dimension 7.

The Dynkin diagram for $E_8$ has no symmetries, therefore all
automorphisms of $(S,W_{E_8})$ are {\em type preserving}.

The 8-vertices are the vertices of {\em root type}. All hemispheres
bounded by walls are centered at an 8-vertex.

These are the one dimensional singular spheres 
(see Remark~\ref{rem:1-2-spheres}) in $(S,W_{E_8})$:
\begin{center}
\vspace{-0.3cm}
 \includegraphics[scale=0.8]{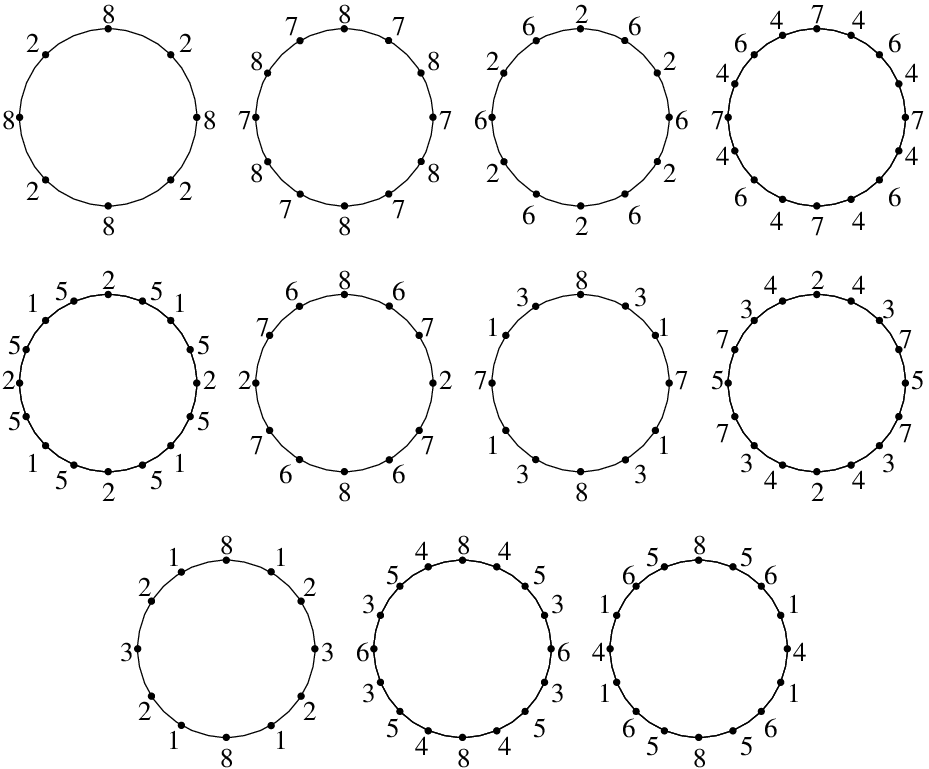}
\end{center}

We describe now the possible lengths and types 
(modulo the action of the Weyl group) of segments
between vertices. 
In Remarks~\ref{rem:segmenttype}, \ref{rem:segmenttypeII}
and \ref{rem:tables}
we explain how these lengths can be 
deduced from the data in the Appendix~\ref{app:E8coxeter}.
The simplicial convex hulls of the segments can be determined
by identifying the smallest singular sphere containing the circle
spanned by the corresponding two vertices.
Here these singular spheres are always 1- and 2-dimensional 
(see Remark~\ref{rem:1-2-spheres} for an explanation of how to determine
the singular 1- and 2-spheres) except for one case, where it is 3-dimensional
(this case is described in detail in the Appendix~\ref{app:E8coxeter}, 
page~\pageref{app:7ptsconvhull}).
We list only the lengths and types that we will need later.

\begin{itemize}
\item 
Distances between two 8-vertices $x$ and $x'$:
\begin{center}
\begin{tabular}{|C{4cm}|C{10cm}|}
\hline
Distance & Simplicial convex hull  of segments $xx'$ \\ \hline
$0,\pi$&\\ \hline
$\pithird$& singular segment of type 878  \\ \hline
$\pihalf$& singular segment of type 828  \\ \hline
$\2pithird$ & singular segment of type 87878 \\
\hline
\end{tabular}
\end{center}

\item
Distances between two 2-vertices $x$ and $x'$:
\begin{center}
\begin{tabular}{|C{4cm}|C{10cm}|}
\hline
Distance & Simplicial convex hull  of segments $xx'$ \\ \hline
$0,\pi$&\\ \hline
$\arccos(\frac{3}{4})$ & singular segment of type 232  \\ \hline
$\pithird$ & singular segment of type 262  \\ \hline
$\arccos(\frac{1}{4})$,
$\arccos(-\frac{1}{4})$ & \hpic{\includegraphics[scale=0.45]{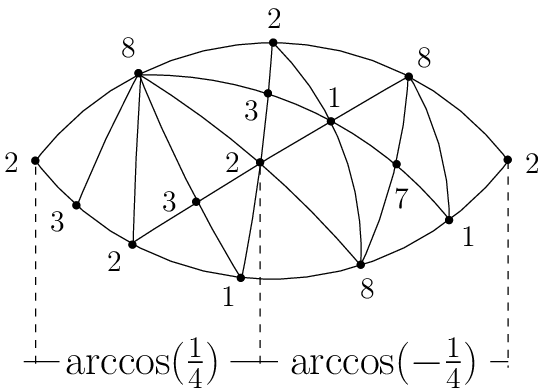}}  \\ \hline
$\pi/2$ & singular segment of type 282 / \hfill\linebreak
		singular segment of type 25152 \\ \hline
$\2pithird$ & singular segment of type 26262  \\ \hline
$\arccos(-\frac{3}{4})$ & singular segment of type 21812  \\ \hline
\end{tabular}
\end{center}

\item The possible distances between two 7-vertices $x$ and $x'$ 
are $\arccos(\frac{k}{6})$
for an integer $-6\leq k\leq 6$. 
Here we will just need to describe the following segments:
\begin{center}
\begin{tabular}{|C{2cm}|C{6.5cm}|m{5cm}|}
\hline
Distance & Simplicial convex hull  of segments $xx'$ & Comments\\ \hline
$\arccos(-\frac{1}{6})$ &
	\hpic{\includegraphics[scale=0.5]{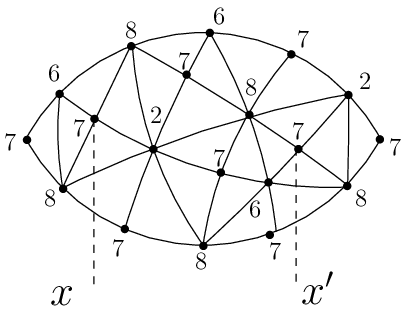}}
	\hpic{\includegraphics[scale=0.5]{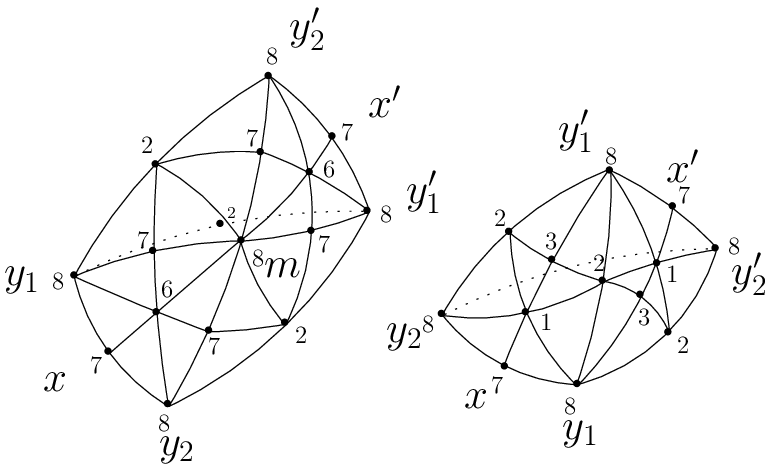}}	&
	There are two types of segments $xx'$.
	The simplicial convex hull $C$ of $xx'$
	is 2- or 3-dimensional, resp.

	For the case $dim(C)=3$,
	we present two perspectives
	from the {\em front} and from {\em behind}
	of a larger polyhedron $C'$.
	It is the simplicial convex hull 
	of $xx'\cup\{y_2,y_2'\}$. 
	(See Appendix~\ref{app:E8coxeter}, p.~\pageref{app:7ptsconvhull}.)
	
	We describe $\Si_m C'$ below. ($\dagger$)
\\ \hline
$\arccos(-\frac{1}{3})$ & singular segment of type 76867 / 
		\hpic{\includegraphics[scale=0.55]{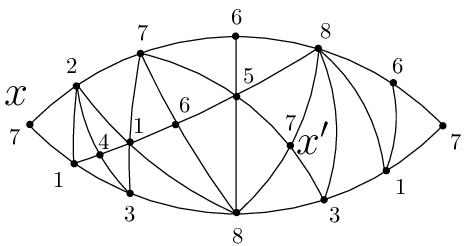}}&\\ \hline
$\arccos(-\frac{2}{3})$ & singular segment of type 7342437 & We present here only the segment
$xx'$, such that the simplex containing $\ora{xx'}$ does not contain 2- or 8-vertices.\\ \hline
\end{tabular}
\end{center}

\parpic[r]{
\hspace{0.5cm}$\Si_m C'$:\hspace{0.3cm}
\hpic{
\includegraphics[scale=0.4]{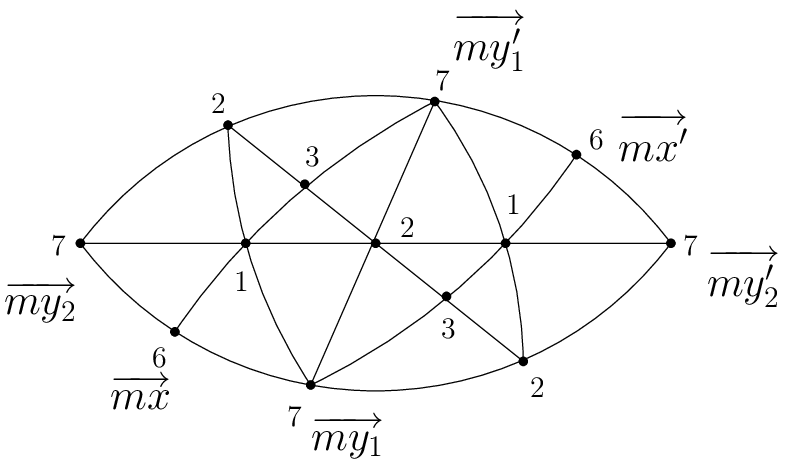}}
}
\bigskip($\dagger$)\\
For a detailed description of the 3-dimensional
spherical polyhedra $C$ and $C'$
we refer to Appendix~\ref{app:E8coxeter}, p.\pageref{app:7ptsconvhull}.
We notice here that $d(\ora{mx},\ora{mx'})=\arccos(-\frac{3}{4})$.

\picskip{0}
The possible lengths of segments $xx'$, 
such that $\pi>d(x,x')>\pihalf$ and the simplex
containing the direction $\ora{xx'}$ in its interior
does not contain a 2- or 8-vertex, are only
$\arccos(-\frac{1}{3})$ and $\arccos(-\frac{2}{3})$.
\end{itemize}

\begin{itemize}
\item Distances $>\pihalf$ and $<\pi$
between two 1-vertices $x$ and $x'$, such that the simplex
containing $\ora{xx'}$ in its interior has no 2-, 7- or 8-vertex
(see Remark~\ref{rem:segmenttypeII}):
\label{page:1vertices}
\begin{center}
\begin{tabular}{|C{4cm}|C{10cm}|}
\hline
Distance & Simplicial convex hull  of segments $xx'$ \\ \hline
$\ac(-\frac{3}{8})$ &  \hpic{\includegraphics[scale=0.65]{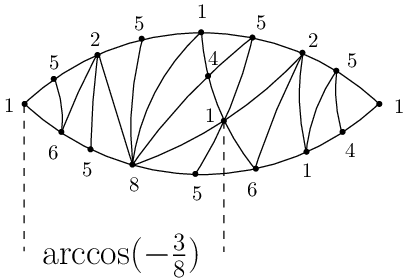}}  \\ \hline
$\2pithird$ & singular segment of type 13831  \\ \hline
$\ac(-\frac{5}{8})$ & \hpic{\includegraphics[scale=0.65]{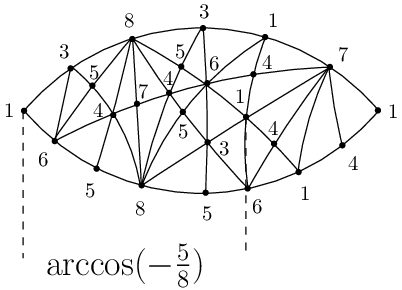}}\\ \hline
$\ac(-\frac{7}{8})$ & singular segment of type 1658561 \\ \hline
\end{tabular}
\end{center}

\item Distances $>\pihalf$ and $<\pi$
between two 6-vertices $x$ and $x'$, such that the simplex
containing $\ora{xx'}$ in its interior has no 1-, 2-, 7- or 8-vertices
(see Remark~\ref{rem:segmenttypeII}):
\label{page:6vertices}
\begin{center}
\begin{tabular}{|C{4cm}|C{10cm}|}
\hline
Distance & Simplicial convex hull  of segments $xx'$ \\ \hline
$\ac(-\quart)$ &  singular segment of type 65856  \\ \hline
$\2pithird$ & \hpic{\includegraphics[scale=0.65]{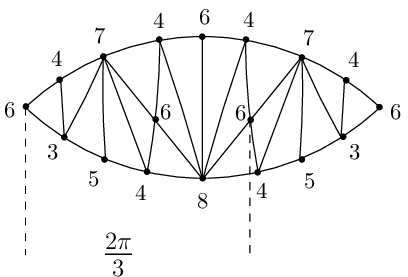}}  \\ \hline
$\ac(-\frac{3}{4})$ & \hpic{\includegraphics[scale=0.65]{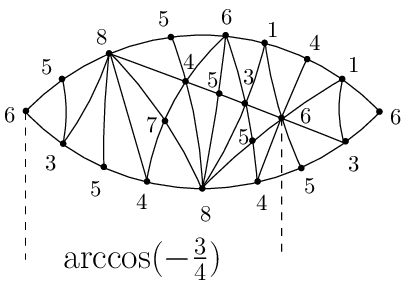}} \\ \hline
\end{tabular}
\end{center}

\item Distances between a 2-vertex $x$ and an 8-vertex $y$:
\begin{center}
\begin{tabular}{|C{4cm}|C{10cm}|}
\hline
Distance & Simplicial convex hull  of segments $xy$ \\ \hline
$\frac{\pi}{4}$ & singular segment of type 28  \\ \hline
$\arccos(\frac{1}{2\sqrt{2}})$ & singular segment of type 218  \\ \hline
$\pihalf$ & singular segment of type 2768 \\ \hline
$\arccos(-\frac{1}{2\sqrt{2}})$ & singular segment of type 23218 \\ \hline
$\frac{3\pi}{4}$ & singular segment of type 2828  \\ \hline
\end{tabular}
\end{center}

\item Distances $> \pihalf$ between a 7-vertex $x$ and an 8-vertex $y$:
\begin{center}
\begin{tabular}{|C{4cm}|C{10cm}|}
\hline
Distance & Simplicial convex hull  of segments $xy$ \\ \hline
$\frac{5\pi}{6}$ & singular segment of type 787878  \\ \hline
$\arccos(-\frac{1}{\sqrt{3}})$ & singular segment of type 72768  \\ \hline
$\arccos(-\frac{1}{2\sqrt{3}})$ & \hpic{\includegraphics[scale=0.55]{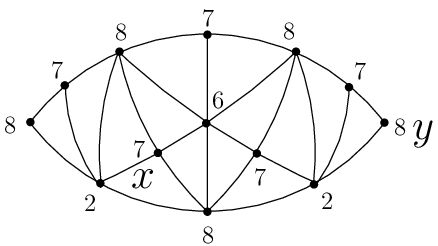}}  \\ \hline
\end{tabular}
\end{center}

\end{itemize}

\section{Convex subcomplexes}

In this section we will describe some general properties of convex subcomplexes
of buildings, as well as some results for buildings of specific types. 
These will be needed later in the proof of the Center Conjecture.

Let $K$ be a convex subcomplex of a spherical building $B$. 

Let $v\in\Si_x K$. We say that $v$ is {\em $d$-extendable}, if there is a segment
$xy \subset K$ of length $d$ and so that $v=\ora{xy}$. 
If $v\in\Si_x K$ is a vertex, we say that it is {\em extendable}, 
if there is a segment
$xy \subset K$ with a vertex adjacent to $x$ (of the same type as $v$)
in its interior and so that $v=\ora{xy}$.

Let $x,y,z\in B$, we define the {\em type of the angle $\angle_x(y,z)$}
as the type of the segment $\ora{xy}\ora{xz}$ in $\Si_x B$.

We say that a point $x \in K$ is {\em interior in $K$}, 
if the link $\Si_x K$ is a subbuilding
of $\Si_x B$.

\begin{prop}[Local conicality, {\cite[Lemma 3.6.1]{KleinerLeeb}}]
\label{prop:locconic}
Let $x,\hat x \in B$ be a pair of antipodal points. Then the union of the
geodesics of length $\pi$ from $x$ to $\hat x$ is isometric to the
spherical join $S^0*\Si_x B$ and contains a neighborhood of $\{x,\hat x\}$.
\end{prop}

The next lemma is an easy consequence of the local conicality property of 
spherical buildings.

\begin{lem}
 \label{lem:sphinlink}
Let $x\in K$ such that $\Si_x K$ contains a sphere 
(or hemisphere) of dimension $d$. Suppose that $x$ has an antipode in $K$.
Then $K$ contains a sphere (hemisphere) of dimension $d+1$
containing $x$ and its antipode. 
\end{lem}

The next two lemmas explain how two properties of points in $K$ 
(being interior or having an antipode) can propagate inside $K$.

\begin{lem}\label{lem:propagationI}
 Let $x\in K$ be an interior point and let $xy\subset K$ be a segment.
Then any point in the interior of the segment $xy$ is also interior in $K$
\end{lem}
\proof
Let $z$ be a point in the interior of the segment $xy$. 
Since $\Si_x K$ is a subbuilding, the link $\Si_{\ora{xz}} \Si_x K$
is also a building. 
In particular,we see that
$\Si_{\ora{zx}} \Si_z K\cong\Si_{\ora{xz}} \Si_x K$ contains a 
sphere of dimension $dim(K)-2$. 
The direction $\ora{zy}\in\Si_z K$ is antipodal to
$\ora{zx}$, then by Lemma~\ref{lem:sphinlink},
it follows that $\Si_z K$ contains a sphere
of dimension $dim(K)-1$. 
Proposition~\ref{lem:Ksubbuild} now implies
that $\Si_z K$ is a subbuilding of $\Si_z B$.
\qed

\begin{lem}\label{lem:antipodes}
 Let $x_1x_2\subset K$ be a segment. Suppose $z$ is a point in the interior of 
the simplicial convex hull of
$x_1x_2$, which has an antipode $\hat z\in K$. Then $x_i$ has also an antipode in $K$.
\end{lem}
\proof
Let $C$ be the simplicial convex hull of $x_1x_2$. Notice that $C$ is contained
in an apartment and $\Si_z C$ is a sphere.
Let $\gamma_i\subset K$ for 
$i=1,2$ be the geodesic connecting $z$ and $\hat z$, such that
the initial direction of $\gamma_i$ at $z$ is the antipode in $\Si_z C$ 
of $\ora{zx_{i}}$. Then
$x_iz \cup \gamma_i$ is a geodesic of length $>\pi$. It is clear that $\gamma_i$
contains an antipode of $x_i$.
\qed

The following results give us conditions, under which $K$ satisfies the conclusions
of the Center Conjecture~\ref{centerconj}.

The next Proposition puts together the results 
\cite[Prop.\ 2.14, Prop.\ 2.15]{LeebRamos}. 
Compare also \cite[Thm.\ 2.2]{Serre} and \cite[Prop.\ 3.10.3]{KleinerLeeb}.

\begin{prop}
\label{lem:Ksubbuild}
The following assertions are equivalent:
\begin{enumerate}
\item[(i)] $K$ is a subbuilding of $B$,
\item[(ii)] every {\em vertex} of $K$ has an antipode in $K$,
\item[(iii)] $K$ contains a sphere of dimension equal to the dimension of $K$.
\end{enumerate}
\end{prop}

The following result was stated in \cite[Cor.\ 2.20]{LeebRamos} for
convex subcomplexes, but the proof works also for closed convex subsets.
In \cite{BalserLytchak} a more general result is shown, namely, for an 
arbitrary $\CAT(1)$ space $C$ of finite dimension and the action
$Isom(C)\acts C$.

\begin{lem}
\label{lem:fixptoncompl}
Let $C\subset B$ be a closed convex subset. Suppose that $rad_C(C) \leq \pihalf$, 
then the action $Stab_{Aut(B)}(C)\acts C$ has a fixed point. 
\end{lem}

\begin{lem}[{\cite[Cor.\ 2.22]{LeebRamos}}]
\label{lem:codim1sph}
 If $K$ contains a singular sphere of dimension $dim(K)-1$, then $K$ is 
a subbuilding or $Stab_{Aut(B)}(K)\acts K$ has a fixed point.
\end{lem}

The following results about convex subcomplexes of type
$D_n$, $E_6$ and $E_7$ are rather technical and may seem
unmotivated at first.
They collect some special geometric information, which
will be repeatedly used in our main argument.
The reader may proceed directly to Section~\ref{sec:centerconj}
and return to them when needed,
however, their proofs are also a good preparation for the 
more complicated arguments in Section~\ref{sec:centerconj}.

\subsection{Convex subcomplexes of buildings of type $D_n$}

In this section let $L\subset B$ be a convex 
subcomplex of a building of type $D_n$ for $n\geq 4$. 
We use the following labelling of the Dynkin diagram
\hpic{\includegraphics[scale=0.45]{Dndynk}}.

\begin{lem} \label{lem:Dn_wall}
If $L$ contains a singular $(n-2)$-sphere $S$ (i.e.\ $S$ is a wall) 
and $x\in L$ is a $1$-, $2$- or $n$-vertex
without antipodes in $S$, then $\Si_x L$ contains an apartment. In particular,
$x$ is an interior vertex in $L$.
\end{lem}
\proof
Let first $x$ be an $n$-vertex. The sphere $S$ contains $n-2$ pairwise orthogonal 
$n$-vertices and their antipodes 
(see Section~\ref{sec:Dncoxeter}, p.~\pageref{pag:spheresDn}). 
They span a singular $(n-3)$-sphere $S'\subset S$.
Since $x$ has no antipodes in $S$ and the distances between two $n$-vertices
are $0,\pihalf,\pi$, 
$x$ must have
distance $\pihalf$ to all these $n$-vertices in $S'$.
Notice that $S'$ is the convex hull of its
$n$-vertices. This implies that $d(x,S')\equiv \pihalf$
 and $h:=CH(S',x)$ is a $(n-2)$-dimensional
hemisphere centered at $x$. 
Put $D_3:=A_3$.
The building $\Si_x B$ has type $D_{n-1}$. 
The link $\Si_x h$ is a 
$(n-3)$-sphere containing $n-2$ pairwise orthogonal $(n-1)$-vertices
(the directions of the segments between $x$ and the $n$-vertices in $S'$). 
This $(n-3)$-sphere
is not a subcomplex (cf. Section~\ref{sec:Dncoxeter}, p.~\pageref{pag:spheresDn}), 
its simplicial convex hull is an apartment contained in $\Si_x L$.

We may now assume w.l.o.g.\ that $x$ is a 1-vertex. 
We prove the assertion by induction on $n$.
Let $B$ be of type $D_3$
with Dynkin diagram \includegraphics[scale=0.4]{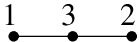}. 
In this case the 1-dimensional sphere $S$ contained in $L\subset B$
is a circle of type 1312321. 
Since the 1-vertex $x$ has no antipodes in
$S$, it must be adjacent to the 2-vertices in $S$ 
and therefore it is also adjacent to the
3-vertex $w$ between them. It follows 
by the local conicality property 
(Proposition~\ref{prop:locconic} applied to the 3-vertex $w$ and
its antipode in $S$)
that the convex 
hull $CH(S,x)$ is a 2-dimensional hemisphere
with $x$ in its interior. 
The link $\Si_x CH(S,x)$ is an apartment in $\Si_x L$.

Let now $B$ be of type $D_n$ for $n\geq 4$. 
Let $y_1,y_2\in S$ be two antipodal $n$-vertices. 
If $x$ lies on a geodesic of length $\pi$ connecting 
$y_1$ and $y_2$, then the geodesic 
$y_1xy_2$ is of type $n21n$. The link $\Si_{y_1}L$ is of type $D_{n-1}$
and contains the wall $\Si_{y_1}S$.
The vertex $\ora{y_1x}$ of type 1 or 2 cannot have antipodes in $\Si_{y_1}S$
because $x$ has no antipodes in $S$.
By induction
 it follows that $\Si_{\ora{xy_1}} \Si_{x}L\cong\Si_{\ora{y_1x}} \Si_{y_1}L$ 
contains an apartment, and 
therefore by Lemma~\ref{lem:sphinlink}, $\Si_x L$ contains also an apartment. 

\parpic{\includegraphics[scale=0.35]{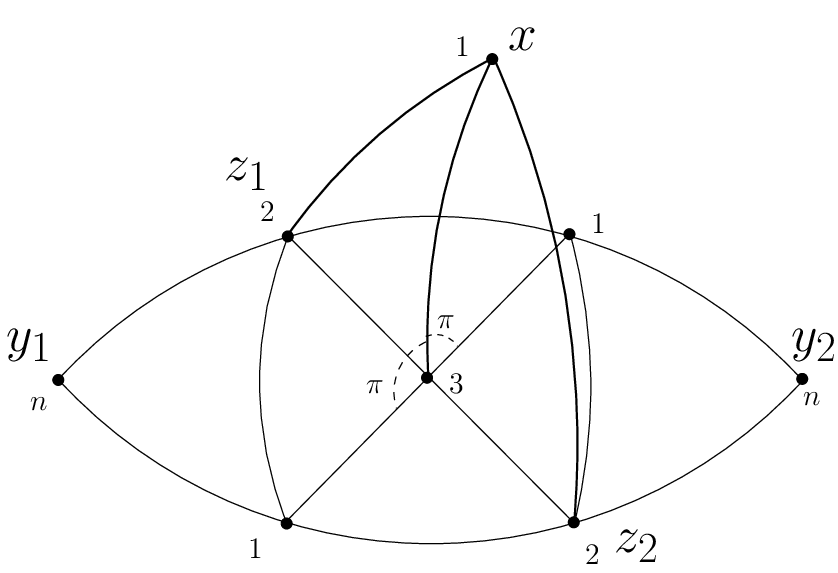}}
On the other hand, if 
$d(x,y_1)+d(x,y_2)>\pi$, then the segments $xy_i$ are of type $12n$.
Let $z_i$ be the 2-vertex
on the segment $xy_i$. Since $z_i$ is adjacent to $y_i$ we deduce that $z_1\neq z_2$.
 Since the link $\Si_x L$ has type $A_{n-1}$, it follows
that the segment $\ora{xz_1}\ora{xz_2}\subset \Si_x L$ is of type 232. 
Again by the induction
 hypothesis, $\Si_{\ora{y_iz_i}}\Si_{y_i}L$ contains an apartment, which in turn
implies that $\Si_{\ora{xz_i}}\Si_{x}L$ contains an apartment. In particular
the 2-vertices $\ora{xz_i}$ are interior vertices in $\Si_x L$. 
Thus, we can extend the segment
$\ora{xz_1}\ora{xz_2}$ to a geodesic in $\Si_x L$ of length $\pi$ and type $232n$. 
In particular $\ora{xz_1}$ has antipodes in $\Si_x L$.
The convex hull of
a small neighborhood in $\Si_x L$ of the interior vertex $\ora{xz_1}$ 
and an antipode in $\Si_x L$ contains the desired
apartment in $\Si_x L$.
\qed

\begin{lem}\label{lem:Dn_n-2sph}
Let $n\geq k \geq 3$.
Suppose that $L$ contains a singular $(n-k)$-sphere 
$S$ containing $n-k+1$ pairwise orthogonal 
$n$-vertices. Assume also that $L$ contains 
a $1$-vertex $x$ and an antipode $x'$ of $x$ (of
type $1$ or $2$ depending on the parity of $n$).
 Then $L$ contains a singular $(n-k+1)$-sphere containing a simplex of type 
$1k(k+1)\dots(n-1)n$.
Moreover, if $n=k$, we can choose the singular $1$-sphere in $L$ to contain 
the 1-vertex $x$ and its antipode $x'$
\end{lem}
\proof
We prove this again by induction on $n$. 
Let $B$ be of type $D_3:=A_3$ with Dynkin diagram
\includegraphics[scale=0.4]{A3_132dynk}, that is $n=k=3$.

The hypothesis in this case is that $L$ contains a pair of antipodal 3-vertices 
$a,a'$ and a pair of antipodal 1- and 2-vertices $x,x'$, respectively. 
If $x'$ lies on a geodesic connecting $a$ and $a'$, then 
this geodesic is of type $3123$.
The convex hull of $x$ and a small neighborhood of $x'$ in this geodesic
is a singular circle of type 2321312 as desired.
\parpic{\includegraphics[scale=0.3]{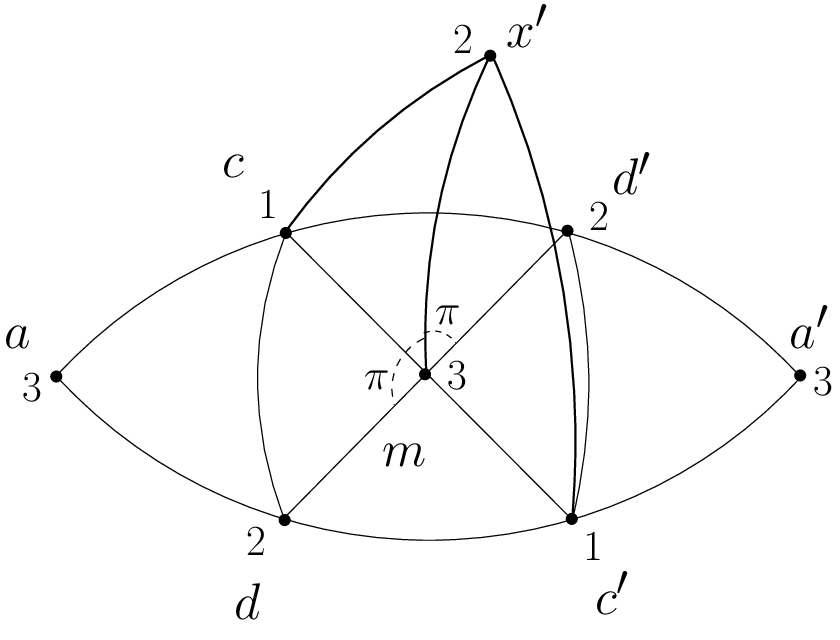}}
Let us suppose then, that $d(a,x')+d(x',a')>\pi$. 
The segments $x'a$ and $x'a'$ are of type 213. 
Let $c,c'$ be the 1-vertices on these segments. 
It is clear that $c\neq c'$ and therefore
the segment connecting them must be of type 131. Let $m:=m(c,c')$.
Since $c$ and $c'$ are adjacent to $x'$, it follows that $m$ is also
adjacent to $x'$.
Let $d,d'$ be the 2-vertices in the
segments $ac'$ and $a'c$ of type 321.
By considering the spherical triangles $CH(a,c,c')$ and
$CH(a',c,c')$, we see that $d$ and $d'$ are adjacent to $m$.
The segment $mx$ is of type 321.
It follows that $x$ must be antipodal to $d$ or $d'$ 
(either $xmd$ or $xmd'$ is a geodesic of length $\pi$)
because $\Si_m B$ is a building of type $A_1\circ A_1$. 
W.l.o.g.\ $x$ is antipodal to $d$, the the convex hull of $x$
and $ac'$ is a circle in $L$ of type 2321312.
The convex hull of $x'$ and a small neighborhood of $x$ in this
circle is again a circle in $L$ of type 2321312 containing both
$x$ and $x'$.

The argument for the induction step is very similar. Let $n\geq 4$.
Let $b,b'$ be a pair
of antipodal $n$-vertices in the $(n-k)$-sphere $S\subset L$.
The links $\Si_b B$ and $\Si_{b'}B$ are of type $D_{n-1}$.
If $b$ lies on a geodesic connecting $x$ and $x'$, then this geodesic is of type
$1n21$, $1n12$, $12n2$ or $12n1$ depending on the parity of $n$
and if $b$ is adjacent to $x$ or $x'$. 
It follows that $\Si_b L$ or $\Si_{b'} L$ 
contains a 1-vertex and an antipode of it.

\parpic{\includegraphics[scale=0.3]{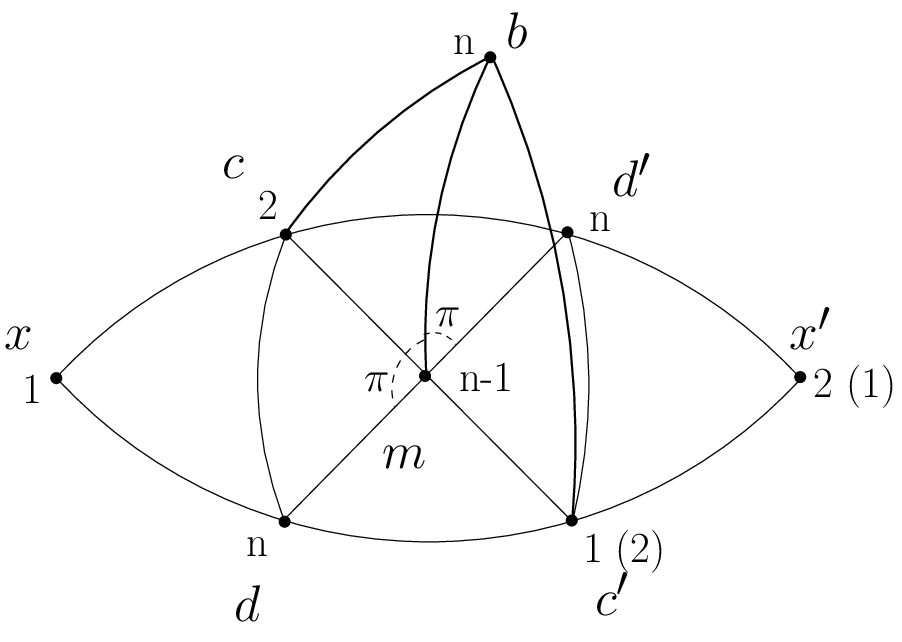}}
If $d(x,b)+d(b,x')>\pi$, then the segment $bx$ is of type $n21$ and the segment 
$bx'$ is of type $n12$ or $n21$. Let $c,c'$ be the vertices in the interior of
the segments $bx$, $bx'$ and let $d,d'$ be the $n$-vertices on the segments
$c'x$ and $cx'$. 
Since $c$ and $d$ are adjacent to $x$, then they are adjacent or $cxd$ is a segment.
In this last case, $c$ and $c'$ must be antipodal, but this cannot happen, because
they are adjacent to $b$. So $c$ and $d$ are adjacent.
This implies that the
segment $cc'$ is of type $2(n-1)1$ or $2(n-1)2$.
The $(n-1)$-vertex $m:=m(d,d')=m(c,c')$ is adjacent to $b$. It
follows that the segment $mb'$ is of type $(n-1)n(n-1)n$. Again we conclude
that $b'$ is antipodal to $d$ or $d'$. This implies that $b'$ lies in a circle 
in $L$ of type $n21n21n$ or $n21n12n$ containing also $x$ or $x'$. 
In particular $\Si_{b'} L$ or
$\Si_b L$ contains a 1-vertex and an antipode of it. 
Suppose w.l.o.g.\ that it holds for $\Si_b L$.
It follows, that $L$ contains a circle spanned by a simplex of type $1n$
containing $x$ and $x'$.
So, if $k=n$, we are done. Suppose then, that $k\leq n-1$.

We have seen that the link $\Si_b L$ of type $D_{n-1}$ contains
a 1-vertex and an antipode of it. It also contains the singular 
$(n-1-k)$-sphere $\Si_ b S$ spanned by $n-k$ pairwise orthogonal $(n-1)$-vertices.
By the induction assumption, $\Si_b L$ contains a singular 
$(n-k)$-sphere spanned by a simplex
of type $1k(k+1)\dots (n-1)$. 
Hence, $L$ contains a singular $(n-k+1)$-sphere spanned by a simplex of type 
$1k(k+1)\dots(n-1)n$.
\qed

\begin{rem}
 If $k=3$ in Lemma~\ref{lem:Dn_n-2sph}, 
then the conclusion is that $L$ contains a wall.
\end{rem}

\begin{cor}\label{lem:D4_circle}
 Let $n=4$, i.e.\ $B$ is of type $D_4$ 
and suppose that $L$ contains a pair of antipodal $i$-vertices
and a pair of antipodal $j$-vertices for $i\neq j$ and $i,j\in\{1,2,4\}$. 
Then $L$ contains a singular circle of type $1241241$.
Moreover, we can choose this singular circle in $L$ to contain 
the two antipodal $i$-vertices or the two antipodal $j$-vertices.
\end{cor}
\proof
Using the symmetries of the Dynkin diagram of type $D_4$
this is just the special version of Lemma~\ref{lem:Dn_n-2sph}
where $n=k=4$.
\qed

\subsection{Convex subcomplexes of buildings of type $E_6$}

In this section let $L\subset B$ be a convex subcomplex of a building of type $E_6$. 
We use the following labelling of the Dynkin diagram
\hpic{\includegraphics[scale=0.4]{E6dynk}}.

\begin{lem} \label{lem:E6_5sph}
If $L$ contains a singular $4$-sphere $S$ 
(i.e $S$ is a wall) and $x\in L$ is a $2$ or $6$-vertex
without antipodes in $S$, then $\Si_x L$ contains an apartment. In particular,
$x$ is an interior vertex in $L$. 
\end{lem}
\proof
By the symmetry of the Dynkin diagram for $E_6$ 
it suffices to show it for a 2-vertex $x\in L$.
The wall $S$ contains a pair of antipodal 2- and 6-vertices $a$ and $a'$, respectively.
The link $\Si_a B$ ($\Si_{a'}B$) is of type $D_5$ and Dynkin
diagram \hpic{\includegraphics[scale=0.4]{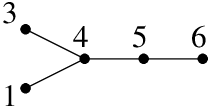}}
(\hpic{\includegraphics[scale=0.4]{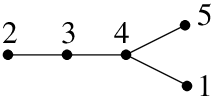}}).
$\Si_a L$ and $\Si_{a'} L$ contain a singular 3-sphere $\Si_a S$, 
respectively $\Si_{a'} S$.
Suppose first that $x$ lies on a geodesic $\gamma$ connecting $a$ and $a'$.
$\gamma$ is of type 23216 or 2626. 
Since $x$ has no antipodes in $S$, the vertex $\ora{ax}$ of type 3 or 6 
has no antipodes in $\Si_a S$.
It follows from Lemma~\ref{lem:Dn_wall}, that $\Si_{\ora{ax}}\Si_a L$
contains an apartment and this implies in turn, that $\Si_{\ora{xa}}\Si_x L$
contains also an apartment. Since $\ora{xa'}\in \Si_x L$ is antipodal to $\ora{xa}$,
this implies that $\Si_x L$ contains an apartment.
\parpic{\includegraphics[scale=0.6]{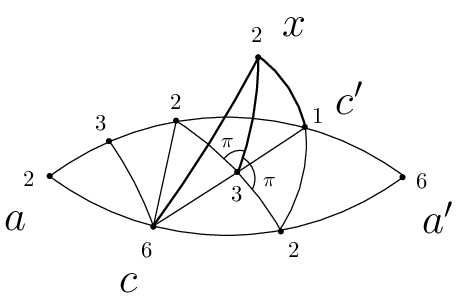}}
On the other hand, if $d(x,a)+d(x,a')>\pi$, then the segments
$xa$ and $xa'$ are of type 262 and 216. Let $c$ be the 6-vertex on $xa$
and let $c'$ be the 1-vertex on $xa'$. 
$c$ is adjacent to $a$ and $c'$ is adjacent to $a'$, therefore
$c$ and $c'$ cannot be adjacent and since both are adjacent to $x$,
it follows that the segment $\ora{xc}\ora{xc'}$ is of type 631.
It follows again from Lemma~\ref{lem:Dn_wall},
that $\Si_{\ora{ax}}\Si_a L$ and $\Si_{\ora{a'x}}\Si_{a'} L$ contain an apartment.
This implies that $\Si_{\ora{xc}}\Si_x L$ and $\Si_{\ora{xc'}}\Si_x L$ contain
an apartment, in particular, $\ora{xc}$ and $\ora{xc'}$ are interior vertices in $\Si_x L$.
The segment $\ora{xc}\ora{xc'}$ is of type 631 and
since $\ora{xc'}$ is interior, it can be extended in 
$\Si_x L$ to a segment of type 6316. This means that 
the interior vertex $\ora{xc}$ has an antipode in 
$\Si_x L$ implying that $\Si_x L$ contains an apartment as desired.
\qed

\subsection{Convex subcomplexes of buildings of type $E_7$}

In this section let $L\subset B$ be a convex subcomplex of a building of type $E_7$. 
We use the following labelling of the Dynkin diagram
\hpic{\includegraphics[scale=0.4]{E7dynk}}.

\begin{lem} \label{lem:E7_5sph}
If $L$ contains a singular $5$-sphere $S$ (i.e.\ $S$ is a wall) and $x\in L$ is a $7$-vertex
without antipodes in $S$, then $\Si_x L$ contains an apartment. In particular,
$x$ is an interior vertex in $L$.
\end{lem}
\proof
The wall $S$ contains a pair of antipodal 7-vertices $a_1,a_2$.
The link $\Si_{a_i} B$ is of type $E_6$ with Dynkin
diagram \hpic{\includegraphics[scale=0.4]{E6dynk}}.
$\Si_{a_i} L$ contains the wall $\Si_{a_i} S$.

Suppose w.l.o.g.\ that $d(x,a_1)=\ac(-\third)$. Then the segment
$xa_1$ is of type 727. Since $x$ has no antipodes in $S$ it follows
that the 2-vertex $\ora{a_1x}$ has no antipodes in $\Si_{a_1} S$. We apply now
Lemma~\ref{lem:E6_5sph} to deduce that $\Si_{\ora{a_1x}}\Si_{a_1} L$
contains an apartment.
This implies in turn, that $\Si_{\ora{xa_1}}\Si_x L$ contains an apartment.
Therefore, if we find an antipode in $\Si_x L$ of $\ora{xa_1}$, we are done.
This is trivial if $x$ lies on a geodesic connecting $a_1$ and $a_2$.

\parpic{\includegraphics[scale=0.6]{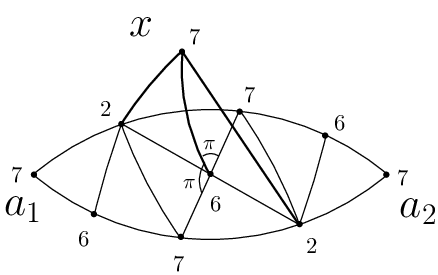}}
Otherwise also $d(x,a_2)=\ac(-\third)$. We may argue as above and conclude
that $\Si_{\ora{xa_2}}\Si_x L$ contains an apartment. 
In particular $\ora{xa_2}$ is an interior
vertex in $\Si_x L$.
Notice that the segment connecting $m(x,a_i)$ for $i=1,2$ cannot be of type 232,
otherwise we find a curve of length $<\pi$ connecting $a_1$ and $a_2$. Therefore,
the segment $\ora{xa_1}\ora{xa_2}$ is of type 262.
Since $\ora{xa_2}$ is interior, we can extend the segment $\ora{xa_1}\ora{xa_2}$
to a segment of type 2626 and length $\pi$ in $\Si_x L$. We have found an antipode
of $\ora{xa_1}$ in $\Si_x L$. 
\qed

\section{The Center Conjecture}\label{sec:centerconj}

Let $B$ be a spherical building and $K\subset B$ a convex subcomplex. We say that $K$
is a {\em counterexample} to the Center Conjecture, if $K$ is not a subbuilding and
the action of $G:=Stab_{Aut(B)}(K)$ has no fixed points in $K$.

From the Proposition~\ref{lem:Ksubbuild} and Lemmata~\ref{lem:fixptoncompl} and
\ref{lem:codim1sph} we can deduce some general properties of 
convex subcomplexes $K\subset B$, which are counterexamples to
the Center Conjecture:

\begin{enumerate}

\item If $x\in K$ and $y\in CH(G\cdot x)$, then there exists
$x'\in G\cdot x$, such that $d(y,x')>\pihalf$.
This is just Lemma~\ref{lem:fixptoncompl} applied to $CH(G\cdot x)$.
In particular, if $x\in K$,
then there exists $x'\in G\cdot x$, such that $d(x,x')>\pihalf$.

Another way to look at this is the following.
If $P$ is a property of vertices in $K$ invariant under the action of $G$,
then for every point $y$ in the convex hull of the $P$-vertices,
we can find a $P$-vertex $x$
with $d(x,y)>\pihalf$.

\item $K$ contains no sphere of dimension $dim(K)-1$
by Lemma~\ref{lem:codim1sph}.

\item If $K$ has dimension $\leq 1$ and is not a subbuilding, 
then by Proposition~\ref{lem:Ksubbuild}, it contains no circles. It follows that $K$
is a (bounded) tree and it has a unique circumcenter, which is fixed by $Isom(K)$.
Hence, a counterexample $K$ has dimension $\geq 2$.
By the main result in \cite{BalserLytchak} mentioned in the introduction,
a counterexample has actually dimension $\geq 3$, but we do not use this fact
in our proof.
\end{enumerate}

Let \textbf{$A$ be the property} of a point in $K$ of not having antipodes in $K$.
Let \textbf{$I$ be the property} of a point in 
$x\in K$ of being interior, i.e.\ $\Si_x K$
is a subbuilding of $\Si_x B$, or equivalently, $\Si_x K$ contains a singular
sphere of dimension $dim(K)-1$.

Notice that an interior point in $K$ cannot have antipodes in $K$, that is,
$I\Rightarrow A$. 
Otherwise $K$ would contain a singular sphere of dimension $dim(K)$
by Lemma~\ref{lem:sphinlink} and 
$K$ would be a subbuilding by Proposition~\ref{lem:Ksubbuild}.
Further, a point $x\in K$ with a sphere of dimension $dim(K)-2$ 
in its link $\Si_x K$ cannot have antipodes in $K$ by Lemma~\ref{lem:codim1sph}.

\subsection{The $E_8$-case}

Let $K$ be a convex subcomplex of a spherical building $B$ of type $E_8$, 
which is a counterexample to the center conjecture.

Our strategy is as follows. 
We focus our attention mainly on the vertices of type 2 and 8. 
The 8-vertices are the vertices of root type and there are few
possibilities for the types of segments between 8-vertices.
The 2-vertices have the second smallest orbit (after the 8-vertices)
under the action of the Weyl group
in the Coxeter complex of type $E_8$.
This implies that the types of the segments between 2-vertices are 
still manageable. 
Another reason to consider 2-vertices is that 
their links have a relatively simple geometry, they are
buildings of type $D_7$. In these buildings, there
is only one type of segments between two distinct non-antipodal 8-vertices, 
namely 878, and it has length $\pihalf$.
First we want to prove that $K$ cannot contain
2- or 8-vertices, whose links contain spheres of large dimension.
This is achieved in the Lemmata \ref{lem:int8pts}-\ref{lem:4sph2ptb}, see also
Lemma~\ref{lem:no8B3pt}.
Then under the assumption of existence of $8A$-vertices, we find
2- and 8-vertices in $K$, with links containing spheres of larger and larger
dimensions, see Corollary~\ref{cor:8Aptwithcircle} and Lemma~\ref{lem:exist8B3pt}. 
This allows us to conclude that all 8-vertices in $K$ have antipodes
in $K$
(Corollary~\ref{lem:no8Apt}). At this point the hard work is already done.
Finally we show that all other vertices in $K$ must also have antipodes in $K$.
This contradicts Proposition~\ref{lem:Ksubbuild} and the assumption that $K$ is not
a subbuilding.

\begin{lem}\label{lem:int8pts}
$K$ contains no $8I$-vertices.
\end{lem}
\proof
Suppose the contrary. There are $8I$-vertices $x_1, x_2\in K$ with distance 
$> \pihalf$. Clearly $I\Rightarrow A$, therefore, $d(x_1,x_2)=\2pithird$ and
the segment $x_1x_2$ is of type 87878. 
Since the $x_i$ are interior vertices, we
can find 7-vertices $y_i\in K$ adjacent to $x_i$ and such that 
$y_1x_1x_2y_2$ is a geodesic of length $\pi$ and type 7878787. The direction
$\ora{y_ix_i}$ is an interior 8-vertex in $\Si_{y_i} K$. Note that $\Si_{y_i} B$
is a building of type $E_6\circ A_1$ and with Dynkin diagram
\hpic{\includegraphics[scale=0.4]{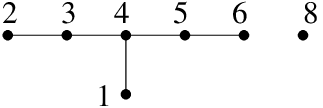}}. It follows that $\Si_{y_i} K$ contains
a top-dimensional (i.e.\ of dimension $dim(K)-1$)
 hemisphere centered at $\ora{y_ix_i}$.
This implies that $K$ contains a hemisphere of dimension $dim(K)$ 
by Lemma~\ref{lem:sphinlink} and in particular it contains a sphere
of dimension $dim(K)-1$.
A contradiction to the properties of a counterexample (cf. Lemma~\ref{lem:codim1sph}).
\qed

\begin{lem}\label{lem:6sph2pt}
$K$ contains no $2$-vertices $x$, such that $\Si_x K$ 
contains an apartment of $\Si_x B$.
\end{lem}
\proof
Let $x$ be such a $2$-vertex in $K$, in particular, $K$ is top-dimensional.
 Then there is another $2$-vertex $x'\in G\cdot x$ at
distance $>\pihalf$ to $x$. 
Notice that $x,x'$ are interior vertices in $K$ and 
therefore, they are $2A$-vertices, by $I\Rightarrow A$.

Case 1: $d(x,x')=\arccos(-\frac{3}{4})$. 
The segment $xx'$ is of type $21812$. 
Since $x$ is an interior vertex, it follows
that the 8-vertex $m(x,x')$ must be interior in $K$
by the propagation of property $I$ along segments (Lemma~\ref{lem:propagationI}).
This contradicts Lemma~\ref{lem:int8pts}.

Notice that again by Lemmata~\ref{lem:propagationI}
and \ref{lem:int8pts}, the 8-vertices
in $\Si_x K$ and $\Si_{x'}K$ are not extendable.
This fact will help us rule out the remaining cases below.

\parpic(5cm,2cm)[r]{\includegraphics[scale=0.4]{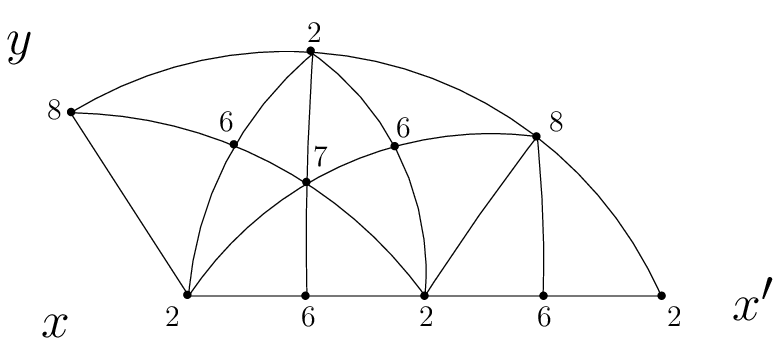}}
Case 2: $d(x,x')=\2pithird$. The segment $xx'$ is of type 26262.
Recall that $\Si_x B$ is of type
$D_7$ with Dynkin diagram \hpic{\includegraphics[scale=0.4]{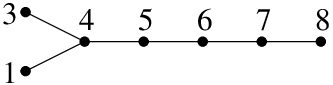}}.
Since $x$ is interior and $K$ is top-dimensional, the link $\Si_x K$ is a building
of type $D_7$ and we can find an 
8-vertex $y\in K$ adjacent to $x$ and such that
$\angle_x(y,x')>\pihalf$.
Then $\angle_x(y,x')=\ac(-\frac{1}{\sqrt{3}})$ 
and this angle must be of type 8676 (i.e.\ the segment $\ora{xy}\ora{xx'}$
is of type 8676).
Since the triangle $CH(y,x,x')$ is spherical (see figure),
it follows that $d(y,x')=\frac{3\pi}{4}$ and the segment
$yx'$ is of type 8282. 
But this cannot happen because the 8-vertices
in $\Si_{x'}K$ are not extendable, hence, a contradiction.

\parpic[r]{\includegraphics[scale=0.5]{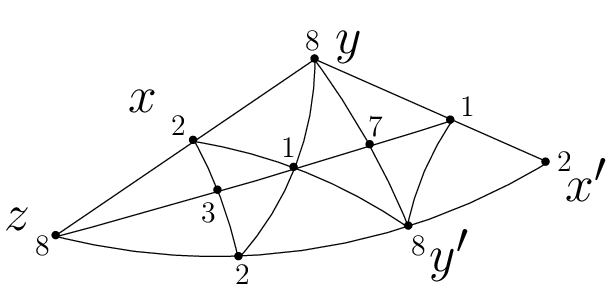}}
Case 3: $d(x,x')=\arccos(-\frac{1}{4})$. The simplicial convex hull of the segment $xx'$
is 2-dimensional and contains 8-vertices $y,y'\in K$ adjacent to $x,x'$. 
Let $z\in K$ be an 
8-vertex adjacent to $x$, such that $zxy$ is a segment of type 828,
which exists because $x$ is an interior vertex. 
One sees that $zy'x'$ is a segment of type 8282 and 
$d(z,x')=\frac{3\pi}{4}$. Again a contradiction 
because 8-vertices in $\Si_{x'}K$ are not extendable.
\qed

\begin{lem}\label{lem:6sph7pt}
$K$ contains no $7$-vertices $x$, such that $\Si_x K$ 
contains an apartment  of $\Si_x B$.
\end{lem}
\proof
Suppose there is such a
7-vertex $x\in K$ and let $y\in K$ be an 8-vertex. 
By $I\Rightarrow A$, $x$ is a $7A$-vertex.
If $d(x,y)=\frac{5\pi}{6}$, then the segment
$xy$ is of type 787878 and we would find interior 8-vertices in $K$
(by the propagation of the property $I$ along segments
Lemma~\ref{lem:propagationI}) which
contradicts Lemma~\ref{lem:int8pts}.
If $d(x,y)=\arccos(-\frac{1}{\sqrt{3}})$, then the segment
$xy$ is of type 72768 and we would find interior 2-vertices 
contradicting Lemma~\ref{lem:6sph2pt}.
So, $d(x,y)\leq\arccos(-\frac{1}{2\sqrt{3}})$.

Let $y_1,y_2\in K$ be 8-vertices adjacent to $x$, 
such that $y_1xy_2$ is a segment of
type 878. Let $x'\in G\cdot x$ with $d(x,x')>\pihalf$. Then 
$d(x',y_i)\leq\arccos(-\frac{1}{2\sqrt{3}})$ and triangle comparison 
with the triangle $(x',y_1,y_2)$ implies that 
$d(x,x')\leq\arccos(-\frac{1}{3})$.

Case 1: $d(x,x')=\arccos(-\frac{1}{3})$.
If the segment $xx'$ is singular of type 76867, then the 8-vertex $m(x,x')$ 
is interior, contradiction.
If the segment $xx'$ has 2-dimensional simplicial convex hull $C$, then there is an
8-vertex $y\in C$ adjacent to $x$ or $x'$. 
Since $x,x'$ are in the same $G$-orbit, we may suppose 
w.l.o.g.\ that $y$ is adjacent to $x$. Let $y'\in K$ be 
another 8-vertex adjacent to $x$ and such that $yxy'$ is a segment of type 878.
Then $d(x',y')=\arccos(-\frac{1}{\sqrt{3}})$ and this case cannot occur by the above.

Case 2: $d(x,x')=\arccos(-\frac{1}{6})$. 
Let $C$ be the simplicial convex hull of the segment
$xx'$. If $C$ is 2-dimensional, there are 8-vertices $y,y'\in C\subset K$ 
adjacent to $x$
and $x'$ respectively. Let $z\in K$ be an 8-vertex 
adjacent to $x$ and such that $zxy$ is a
segment of type 878. Define $z'$ analogously. 
Then $d(x',z)$ or $d(x,z')=\arccos(-\frac{1}{\sqrt{3}})$, which is not possible.

If $C$ is 3-dimensional (see Section~\ref{sec:E8coxeter} for a description of $C$),
there is an 8-vertex 
$m\in C$, such that the segments $mx$ and $mx'$
are of type 867 and $\angle_m(x,x')=\arccos(-\frac{3}{4})$.
Since $x,x'$ are interior vertices, there exist 2-vertices $u,u'\in K$,
such that $mxu$ and $mx'u'$ are segments of length $\pihalf$ and of type 8672.
The fact that $\angle_m(x,x')=\arccos(-\frac{3}{4})$ 
implies that $\pi > d(u,u')\geq \arccos(-\frac{3}{4})$.
Hence $d(u,u')=\arccos(-\frac{3}{4})$.

\parpic{\includegraphics[scale=0.5]{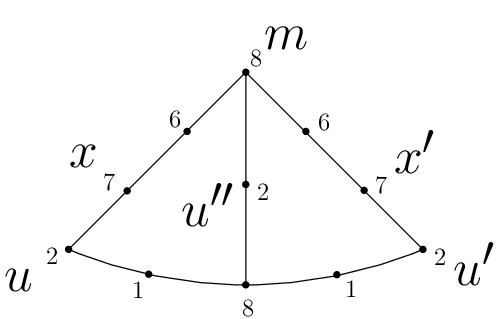}}
The segment $uu'$ is of type 21812 and $CH(m,u,u')$ is a 
(non-simplicial) spherical triangle with a 2-vertex
 $u'':=m(m,m(u,u'))$ in its interior.
This implies that the segment $xu''$ can be extended in $K$ 
beyond $u''$. In particular $u''$
is an interior 2-vertex by the propagation of property $I$ along segments
(Lemma~\ref{lem:propagationI})
and because $x$ is an interior vertex. Since $u''$ is interior and 
$K$ is top-dimensional, $\Si_{u''}K$ contains an apartment.
This contradicts Lemma~\ref{lem:6sph2pt}.
\qed

\begin{lem}\label{lem:5sph8pt}
$K$ contains no $8$-vertices $x$, such that $\Si_x K$ 
contains a singular $5$-sphere, i.e.\ a wall.
\end{lem}
\proof
Let $x_1$ be an $8$-vertex, such that $\Si_{x_1} K$ contains a wall $S_1$.
By Lemma~\ref{lem:codim1sph}, $x_1$ is an $8A$-vertex. 
Let $x_2\in G\cdot x_1$
be at distance $\2pithird$ to $x_1$, which exists by the properties of a
counterexample. 
The link $\Si_{x_2} K$ contains also a wall $S_2$.

If $\ora{x_ix_{3-i}}$ has an antipode in $S_i$ for $i=1,2$, 
then there are 7-vertices $y_i\in K$
adjacent to $x_i$, such that $y_1x_1x_2y_2$ is a 
geodesic of length $\pi$ and type 7878787. 
The direction $\ora{x_iy_i}$ lies on the 5-sphere $S_i$,
hence $\Si_{\ora{y_ix_i}}\Si_{y_i}K$ contains a 4-sphere. 
Note that $\Si_{y_i} B$
is a building of type $E_6\circ A_1$ and with Dynkin diagram
\hpic{\includegraphics[scale=0.4]{E8link7dynk}}. 
It follows that $\Si_{y_i} K$ contains
a 5-dimensional hemisphere centered at $\ora{y_ix_i}$.
Then by Lemma~\ref{lem:sphinlink}, $K$ contains a 6-dimensional
hemisphere $h\subset K$.
The midpoint $z:=m(x_1,x_2)$ is again an $8A$-vertex and it is the center of  
the hemisphere $h$.
In particular, $\Si_z K$ contains the 5-sphere $\Si_z h$ and the 7-vertices in 
this sphere are all $\pihalf$-extendable.
Let $z' \in G\cdot z$ be at distance $\2pithird$ to $z$. 
Since $z'$ is an $8A$-vertex and 
the 7-vertices in $\Si_z h$ are $\pihalf$-extendable, we deduce that
$\ora{zz'}$ has no antipodes in $\Si_z h$. 
It follows from Lemma~\ref{lem:E7_5sph} that $\Si_{\ora{zz'}}\Si_{z}K$ 
contains an apartment and, in particular, that 
for the 8A-vertex $w:=m(z,z')$ the link
$\Si_{\ora{wz}}\Si_{w}K\cong \Si_{\ora{zz'}}\Si_{z}K$ contains an apartment. 
It follows that $\Si_w K$ contains also an apartment of $\Si_w B$
because the direction $\ora{wz'}$ is antipodal to $\ora{wz}$
 (see Lemma~\ref{lem:sphinlink}). 
We obtain a contradiction to Lemma~\ref{lem:int8pts}.

We may therefore assume
w.l.o.g.\ that $\ora{x_1x_2}$ has no antipodes in $S_1$.
Using again Lemma~\ref{lem:E7_5sph} we 
conclude that $\Si_{\ora{x_1x_2}}\Si_{x_1}K$ contains an apartment
and this implies that
$\Si_z K$ contains an apartment for $z=m(x_1,x_2)$. 
Again a contradiction to Lemma~\ref{lem:int8pts}.
\qed

\begin{lem}\label{lem:5sph2pt}
$K$ contains no $2$-vertices $x$, such that 
$\Si_x K$ contains a singular $5$-sphere $S$, i.e.\ a wall.
\end{lem}
\proof
Suppose there is such an $x\in K$. Lemma~\ref{lem:codim1sph} implies that
$x$ is a $2A$-vertex.
Let $y\in K$ be an 8-vertex. If $d(x,y)=\frac{3\pi}{4}$, 
then the segment $xy$ is of type
2828. Let $y'$ be the 8-vertex between $x$ and $y$. 
The link $\Si_x K$ is of type $D_7$ and contains a wall.
If the direction $\ora{xy'}$ has an antipode in this wall, then
$\Si_{\ora{xy'}}\Si_x K$ contains a 4-sphere (i.e.\ a wall).
On the other hand, if $\ora{xy'}$ has no antipode in this wall
then Lemma~\ref{lem:Dn_wall} implies that 
$\Si_{\ora{xy'}}\Si_x K$ contains an apartment. Hence,
by Lemma~\ref{lem:sphinlink}
$\Si_{y'} K$ contains at least a singular 5-sphere, 
contradicting Lemma~\ref{lem:5sph8pt}. 
So $d(x,y)\leq\arccos(-\frac{1}{2\sqrt{2}})$
for all 8-vertices $y\in K$.

Let  $x'\in G\cdot x$ with  $d(x,x')>\pihalf$.
It also holds $d(x',y)\leq\arccos(-\frac{1}{2\sqrt{2}})$
for all 8-vertices $y\in K$.

\parpic[r]{\includegraphics[scale=0.45]{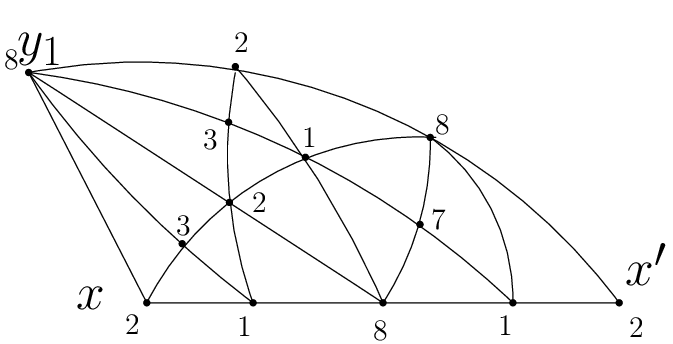}}
Case 1: $d(x,x')=\arccos(-\frac{3}{4})$. The segment $xx'$ is of type $21812$. 
Let $y_1,y_2\in K$ be 8-vertices adjacent to $x$, such that $y_1xy_2$ is a segment of
type 828. These vertices can be found, because $\Si_x K$ contains a wall.
We may assume that $\angle_x(y_1,x')\geq\pihalf$. This implies that the angle
$\angle_x(y_1,x')=\ac(-\frac{1}{\sqrt{7}})$ 
and this angle is of type 831, 
because $\Si_x B$ is a building of type $D_7$. 
$CH(y,x,x')$ is a spherical triangle, therefore we can compute that 
$d(y_1,x')=\frac{3\pi}{4}$.
A contradiction to the observation above.

Case 2: $d(x,x')=\2pithird$. 
If there is an $8$-vertex $y$ adjacent to $x$, such that 
$d(\ora{xx'},\ora{xy})>\pihalf$, then $d(x',y)=\3piquart$
(see the figure in the proof of Lemma~\ref{lem:6sph2pt}, Case 2),
which is not possible.
Therefore we see that the distance between $\ora{xx'}$ and all the 8-vertices
in $S$ must be $\pihalf$ because they come in pairs of antipodes.
This implies that $d(\ora{xx'},S')\equiv \pihalf$, where $S'\subset S$ 
is the 4-sphere spanned by the 8-vertices in $S$. 
Hence the convex hull $CH(S',\ora{xx'})$ is a 5-dimensional hemisphere
(the spherical join of $S'$ and $\ora{xx'}$).

The segments in $\Si_x K$ of length $\pihalf$ connecting the 6-vertex $\ora{xx'}$
and an 8-vertex $\in S'$ are of type 658. 
This implies that $\Si_{\ora{xx'}}\Si_x K$ contains
a 4-sphere spanned by five pairwise orthogonal 5-vertices, 
but this is impossible in a building
of type $D_4\circ A_2$ with Dynkin diagram 
\hpic{\includegraphics[scale=0.4]{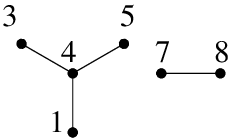}}.

Case 3: $d(x,x')=\arccos(-\quart)$. 
Let $y$ be the 8-vertex adjacent to $x$ contained in the
simplicial convex hull of $xx'$ (see table
of distances between two 2-vertices in Section~\ref{sec:E8coxeter}). 
The direction $\ora{xy}$ cannot have antipodes in $\Si_x K$. 
Otherwise there is an 8-vertex $z\in K$, such that $zxy$ is a segment of type 828
and as in Lemma~\ref{lem:6sph2pt} (Case 3), we see that $d(x',z)=\3piquart$.
It follows then from Lemma~\ref{lem:Dn_wall}, that $K$ is top-dimensional and
$\ora{xy}$ is interior in $\Si_x K$ (i.e.\ its link contains an apartment).
This implies that the direction $\ora{xx'}$ is also interior in $\Si_x K$
by Lemma~\ref{lem:propagationI}.
Then the 7-vertex $m(x,x')$ must be interior in $K$
(its link $\Si_{m(x,x')} K$ contains an apartment) by Lemma~\ref{lem:sphinlink}.
A contradiction to Lemma~\ref{lem:6sph7pt}.
\qed

\begin{lem}\label{lem:5sph7pt}
$K$ contains no $7$-vertices $x$, such that $\Si_x K$ contains a wall
 $S$ of type $1$, that is, a wall containing a pair of antipodal 8-vertices.
\end{lem}
\proof
We proceed exactly as in the proof of Lemma~\ref{lem:6sph7pt}. Recall that
$\Si_x B$ is of type $E_6\circ A_1$.
Suppose there is such an $x\in K$ and let $y\in K$ be an 8-vertex.
Again, by Lemmata~\ref{lem:sphinlink} and
\ref{lem:codim1sph}, $x$ is a $7A$-vertex (see also the properties of a 
counterexample).
If $d(x,y)=\frac{5\pi}{6}$, then the segment
$xy$ is of type 787878. The direction $\ora{xy}$ has an antipode in $S$, 
therefore the link $\Si_{\ora{xy}}\Si_x K$ contains a 
wall by Lemma~\ref{lem:sphinlink}.
This implies again by Lemma~\ref{lem:sphinlink} that
the link $\Si_{y'} K$ of the $8$-vertex $y'$ on the segment 
$xy$ adjacent to $x$ contains 
a wall, contradicting Lemma~\ref{lem:5sph8pt}. 
If $d(x,y)=\arccos(-\frac{1}{\sqrt{3}})$, then the segment
$xy$ is of type 72768. Lemma~\ref{lem:E6_5sph} implies
that $\ora{xy}$ has an antipode in $S$ or $\Si_{\ora{xy}}\Si_x K$ 
contains an apartment.
 In both
cases the link $\Si_z K$ of the $2$-vertex $z$ on the segment 
$xy$ adjacent to $x$ contains 
a wall. A contradiction to Lemma~\ref{lem:5sph2pt}.
So, $d(x,y)\leq\arccos(-\frac{1}{2\sqrt{3}})$.

Let $y_1,y_2\in K$ be 8-vertices adjacent to $x$, such that $y_1xy_2$ is a segment of
type 878. Let $x'\in G\cdot x$ with $d(x,x')>\pihalf$. Then 
$d(x',y_i)\leq\arccos(-\frac{1}{2\sqrt{3}})$ and triangle comparison
with the triangle $(x',y_1,y_2)$ implies that 
$d(x,x')\leq\arccos(-\frac{1}{3})$.

Case 1: $d(x,x')=\arccos(-\frac{1}{3})$.
If the segment $xx'$ is singular of type 76867, then Lemma~\ref{lem:E6_5sph} implies
that the 6-vertex $\ora{xx'}$ has an antipode in $S$
or $\Si_{\ora{xx'}}\Si_x K$ contains an apartment.
Either way, the link in $K$ of the 8-vertex $m(x,x')$ contains a wall, 
which is not possible by Lemma~\ref{lem:5sph8pt}.
The case, where the segment $xx'$ has 2-dimensional simplicial convex hull $C$, 
is dealt with as in the
proof of Lemma~\ref{lem:6sph7pt} (Case 1).

Case 2: $d(x,x')=\arccos(-\frac{1}{6})$. 
Let $C$ be the simplicial convex hull of the segment
$xx'$. If $C$ is 2-dimensional, we argue as in the
proof of Lemma~\ref{lem:6sph7pt} (Case 2).

If $C$ is 3-dimensional (see Section~\ref{sec:E8coxeter} for a description of $C$),
there is an 8-vertex $m\in C$, such that the segments $mx$ and $mx'$
are of type 867 and $\angle_m(x,x')=\arccos(-\frac{3}{4})$.
$C$ contains also 8-vertices $y_1,y_1'$ adjacent to $x,x'$ respectively.
Let $y_2\in K$ be an 8-vertex adjacent to 
$x$ and such that $y_2xy_1$ is a segment of type 878.
Define $y_2'$ analogously. 
Then the angle $\angle_m(x,y_2')$ is of type 6727 (compare with
$\Si_m C'$ in Section~\ref{sec:E8coxeter}). This implies that
$d(x,y_2')=\arccos(-\frac{1}{2\sqrt{3}})$.

\parpic{\includegraphics[scale=0.65]{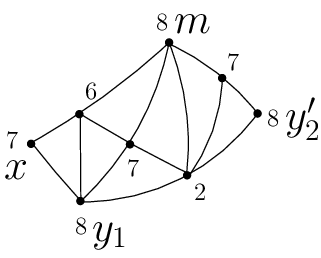}}
 If the 6-vertex $\ora{xm}$ has no antipodes in
$S$, then it follows from Lemma~\ref{lem:E6_5sph} that $\Si_{\ora{xm}}\Si_x K$ contains
an apartment, i.e.\ $\ora{xm}$ is interior in $\Si_x K$. In particular the link $\Si_w K$ of the
7-vertex $w$ in the interior of the simplicial convex hull
of $xy_2'$ contains an apartment. A contradiction to Lemma~\ref{lem:6sph7pt}.
It follows that $\ora{xm}$, and analogously,
$\ora{x'm}$ have antipodes in the walls $S\subset \Si_x K$,
respectively $S'\subset \Si_{x'} K$. Therefore, there exist 2-vertices $u,u'\in K$,
such that $mxu$ and $mx'u'$ are segments of length $\pihalf$ and of type 8672.
$\angle_m(x,x')=\arccos(-\frac{3}{4})$ implies that $\pi > d(u,u')\geq \arccos(-\frac{3}{4})$.
Hence $d(u,u')=\arccos(-\frac{3}{4})$. 

\parpic[r]{\includegraphics[scale=0.45]{lem6sph7pt}}
It follows that the
segment $uu'$ is of type 21812 and $CH(m,u,u')$ is a 
(non-simplicial) spherical triangle.
The segment $m\,m(u,u')$ has length $\pihalf$ and therefore it has type 828.
The 2-vertex $u'':=m(m,m(u,u'))$ lies in the interior of the spherical triangle 
$CH(m,u,u')$.

\parpic{\includegraphics[scale=0.4]{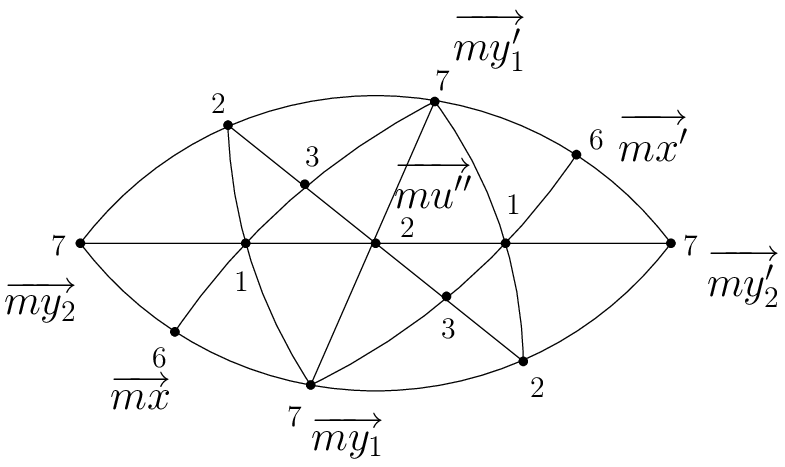}}
Consider the link of $m$. Since $\ora{xm}$ 
has an antipode in the wall $S\subset \Si_x K$,
it follows that $\Si_{\ora{mx}}\Si_m K$ contains a wall.
The link $\Si_{\ora{mx}}\Si_m B$ is of type $D_5\circ A_1$.
The wall in $\Si_{\ora{mx}}\Si_m K$ contains a wall in the $D_5$-factor.
The direction $\xi:=\ora{\ora{mx}{\scriptstyle\ora{my_1'}}}$ 
is a 1-vertex in $\Si_{\ora{mx}}\Si_m K$. 
By Lemma~\ref{lem:Dn_wall} 
we conclude that the $A_4$-factor of $\Si_{\xi}\Si_{\ora{mx}}\Si_m K$
contains at least a wall. Taking spherical join with the directions to the
7-vertices $\ora{my_2}$ and $\ora{my_1}$ 
we find a wall in $\Si_{\xi}\Si_{\ora{mx}}\Si_m K$.
This implies that $\Si_{\ora{mu''}}\Si_m K$ contains at least a 
wall. Since $\ora{mu''}$ is extendable, it follows that $\Si_{u''}K$ contains a wall. 
But this contradicts Lemma~\ref{lem:5sph2pt}.
\qed

Before we continue,
we describe first some configurations of points of $K$, which will be used 
several times during
the rest of the argument.

Let $P$ be a $G$-invariant property of 8-vertices implying $A$ 
(the property of not having antipodes in $K$)
and suppose there are $8P$-vertices in $K$. 
The following configurations appear when we try to determine
how the set of $8P$-vertices in $K$ looks like. This set has two
main restrictions. Its vertices cannot be too close to each other, otherwise
we obtain a $G$-invariant set of radius $\leq\pihalf$ contradicting the properties
of a counterexample. Its vertices can also not be too far from each other
in the sense that they do not have antipodes.

\parpic{\includegraphics[scale=0.61]{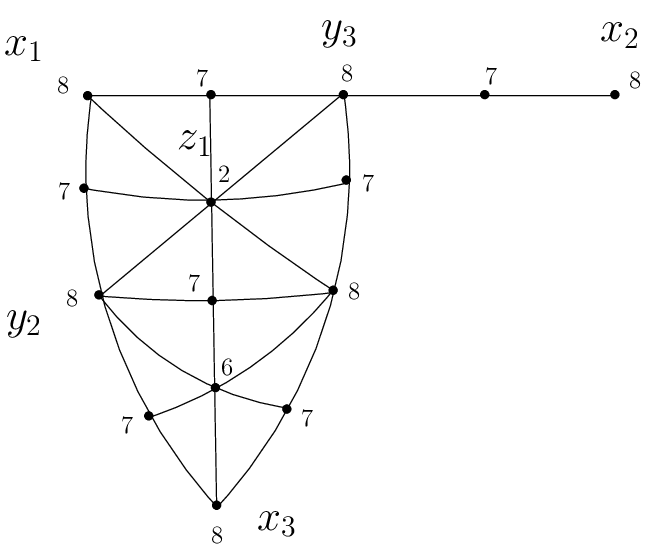}}
Let us begin with an $8P$-vertex $x_1\in K$.
Since $K$ is a counterexample, there is another $8P$-vertex $x_2\in K$ at
distance $>\pihalf$ to $x_1$. Since they do not have antipodes, it follows that 
$d(x_1,x_2)=\2pithird$. Let $y_3:=m(x_1,x_2)$, it is an $8A$-vertex
by Lemma~\ref{lem:antipodes}. 
Again there is an $8P$-vertex $x_3\in K$, such that $d(y_3,x_3)=\2pithird$
because $y_3$ lies in the convex hull of the $8P$-vertices in $K$. 
Notice that, since
$x_i$ are $8A$-vertices, $0<\angle_{y_3}(x_3,x_i)<\pi$ for $i=1,2$,
because otherwise $x_1$ or $x_2$ would be an antipode of $x_3$.
We may assume 
w.l.o.g.\ that $\angle_{y_3}(x_3,x_1)\geq \pihalf$. 
The link $\Si_{y_3}B$ is a building
of type $E_7$ and with Dynkin diagram 
\hpic{\includegraphics[scale=0.4]{E7dynk}}.
It follows that $\angle_{y_3}(x_3,x_1)=
\linebreak[3]
\arccos(-\frac{1}{3})$ and this
angle is of type 727, i.e.\ 
the segment $\ora{y_3x_1}\ora{y_3x_3}\subset \Si_{y_3} K$
is singular of type 727.
The convex hulls $CH(x_3,y_3,m(x_1,y_3))$ and 
$CH(x_3, x_1, m(x_1,y_3))$ are spherical triangles,
because $y_3$ and $m(x_1,y_3)$ ($x_1$ and $m(x_1,y_3)$, 
respectively) are contained in a 
common Weyl chamber and therefore $x_3$, $y_3$ and $m(x_1,y_3)$
($x_3$, $x_1$ and $m(x_1,y_3)$, respectively) lie in a common apartment.
This implies that $d(m(x_1,y_3),x_3)=\frac{5\pi}{6}$ and 
the segment $m(x_1,y_3)x_3$ is of type 72768 
(cf. tables in Section~\ref{sec:E8coxeter}). 
Since $\Si_{m(x_1,y_3)} B$ 
is of type $E_6\circ A_1$ with Dynkin diagram 
\hpic{\includegraphics[scale=0.4]{E8link7dynk}}, it follows that 
$\angle_{m(x_1,y_3)}(x_1,x_3)=\angle_{m(x_1,y_3)}(x_1,x_3)=\pihalf$.
Hence, the convex hull 
$CH(x_1,y_3,x_3)$ is the union of $CH(x_3,y_3,m(x_1,y_3))$ 
and $CH(x_3,x_1,m(x_1,y_3))$,
and it is an isosceles spherical triangle with sides of type
878, 87878 and 87878.
Let $y_2:=m(x_1,x_3)$ and $z_1:=m(y_2,y_3)$.

\begin{frameenv}\begin{center}
We refer to this configuration of $8P$-vertices as configuration $\ast$.
\end{center}\end{frameenv}

\bigskip
\parpic{\includegraphics[scale=0.4]{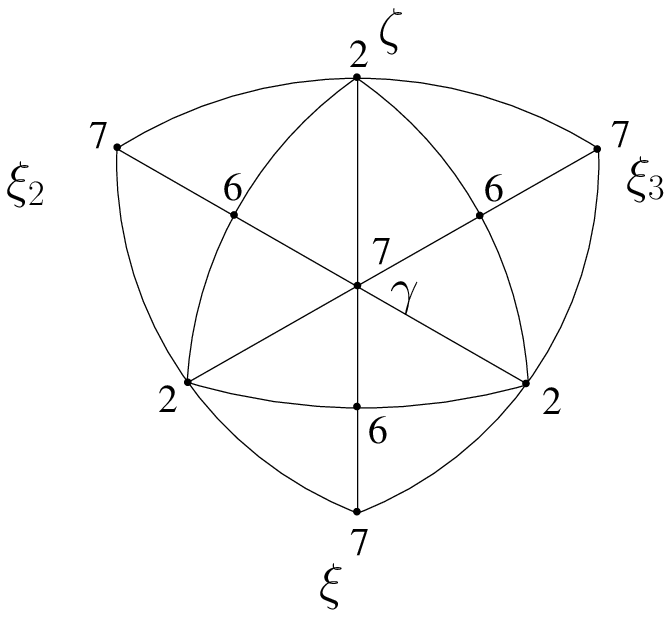}}
Let $\xi_i:=\ora{x_1x_i}$ for $i=2,3$ and 
$\zeta:=\ora{x_1z_1}=m(\xi_2,\xi_3)$.
We will add now a new $8$-vertex to our configuration $\ast$ with the
following properties and study the new configuration that this new
vertex forces, in particular,
we describe a subset of the convex hull $CH(x_1,x_2,x_3, x)$
(see the figure in the next paragraph).
Suppose there is an 8-vertex $x\in K$ at distance $\pithird$ to $x_1$
and let
$\xi:=\ora{x_1x}$. 
Assume furthermore that $d(\zeta,\xi)=\arccos(-\frac{1}{\sqrt{3}})$,
then the segment $\xi\zeta$ is of type 7672 
(see table in Section~\ref{sec:E7coxeter}).

Recall that $\xi_i$ is $\2pithird$-extendable 
to $8A$-vertices and $\xi$ is $\pithird$-extendable. 
Thus, $d(\xi,\xi_i)<\pi$ for $i=2,3$. 
It follows that $\angle_\zeta(\xi,\xi_i)=\pihalf$ 
because of the possible distances between 7-vertices in buildings
of type \hpic{\includegraphics[scale=0.4]{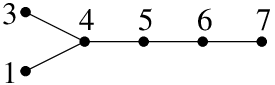}}.
The convex hulls $CH(\xi,\zeta,\xi_i)$ are spherical triangles and one
can compute $d(\xi,\xi_i)=\arccos(-\frac{1}{3})$ for $i=2,3$.
The convex hull
$CH(\xi,\xi_2,\xi_3)$ is the union of the spherical triangles
$CH(\xi,\zeta,\xi_i)$ for $i=2,3$.
It is an equilateral spherical triangle with sides of type 727.
Let $\gamma$ be the 7-vertex at the center of this triangle.
\parpic[r]{\includegraphics[scale=0.75]{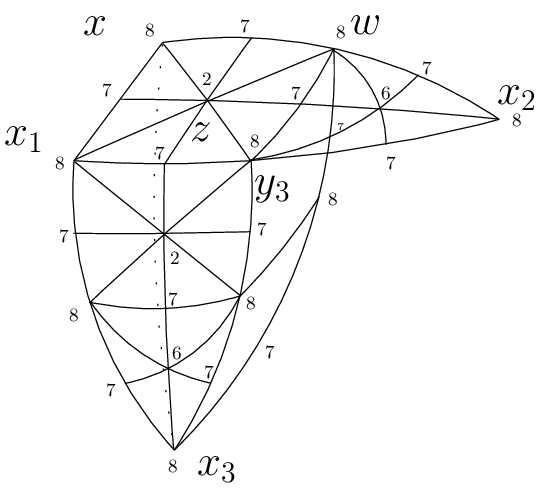}}
From $d(\xi,\xi_i)=\arccos(-\frac{1}{3})$, it follows for $i=2,3$ 
that $d(x,x_i)=\2pithird$
and the convex hulls $CH(x_1,x,x_i)$ are isosceles spherical triangles
(compare with the spherical triangle $CH(x_1,y_3,x_3)$ above).
Let $w:=m(x,x_2)$ be the $8A$-vertex between $x$ and $x_2$. 
Then by considering the triangle $CH(x_1,x,x_2)$, we see that
$\omega:=\ora{x_1w}=m(\xi,\xi_2)$. 
Let $z:=m(x_1,w)$ be the 2-vertex between $x_1$ and
$w$, then $\ora{x_1z}=m(\xi,\xi_2)$.
The angle $\angle_{x_1}(z,x_3)=\arccos(-\frac{1}{\sqrt{3}})$
is of type 2767 (compare with the triangle $CH(\xi,\xi_2,\xi_3)$).
Notice that $CH(z,x_1,x_3)$ is a spherical triangle, 
this implies that $d(z,x_3)=\3piquart$. 

\parpic{\includegraphics[scale=0.5]{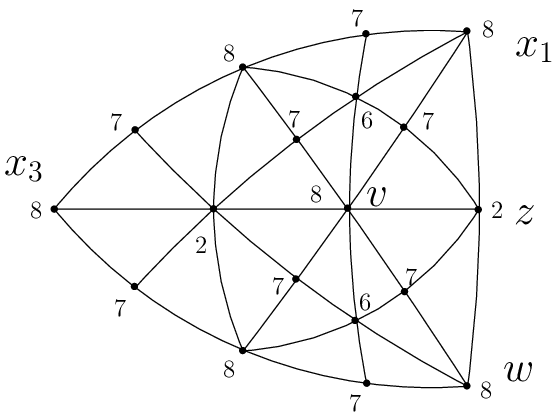}}
The segment $zx_3$ is of type 2828.
Let $v$ be the $8A$-vertex on the segment $zx_3$ adjacent to $z$.
Recall that $x_3$ is an $8A$-vertex. Then $x_3$ cannot be antipodal to
$w$, thus $d(x_3,w)=\2pithird$ and
$\angle_z(x_3,x_1)=\angle_z(x_3,w)=\pihalf$.
Recall also that $d(x_3,y_3)=d(x_3,x)=\2pithird$, therefore
$\angle_z(x_3,y_3)=\angle_z(x_3,x)=\pihalf$.
The convex hulls $CH(x_3,x_1,w)$ and $CH(x_3,y_3,x)$ are isosceles spherical
triangles with sides of type 87878, 87878 and 828.

The convex hull in $\Si_z K$ of the 8-vertices
$\ora{zx}$, $\ora{zx_1}$ $\ora{zy_3}$ $\ora{zw}$ and $\ora{zv}$ is a 2-dimensional
singular hemisphere $h$ centered at $\ora{zv}$. 
Let $s\subset \Si_z B$ be a singular 2-sphere containing $h$
and let $\wht{x_3}$ be an 8-vertex in $B$, such that it is adjacent to $z$ and
$\ora{z\wht {x_3}}$ is the antipode of $\ora{zv}$ in $s$.
It follows that $\wht {x_3}$ is antipodal to $x_3$ in $B$. 
The convex hull in $B$
of $x_3,\wht{x_3},x,x_1,y_3,w$ is a 3-dimensional spherical bigon 
(i.e.\ a 3-dimensional spherical polyhedron with just two (antipodal) vertices)
connecting
$x_3$ and $\wht{x_3}$, with edges
$x_3\alpha \wht{x_3}$ for $\alpha\in\{x,x_1,y_3,w\}$ of type 8787878. 
It follows that
the convex hull $CH(x_1,x,w,y_3,x_3)$ is a 
(3-dimensional) spherical convex polyhedron in $K$ obtained
by truncating this spherical bigon.
Notice that the 7-vertex $\gamma$ at the center of the triangle 
$CH(\xi,\xi_2,\xi_3)\subset \Si_{x_1} K$ is $\2pithird$-extendable
in $K$ to the $8A$-vertex $m(x_3,w)$.

\begin{frameenv}\begin{center}
We refer to this configuration in $K$ as configuration $\aast$.
\end{center}\end{frameenv}

Now we continue with the proof of our main result. 
Recall that our strategy is to prove first that $K$ cannot contain
$2$- or $8$-vertices, whose links contain spheres of large dimension.
In a special case we can also exclude 8-vertices, whose links contain a 3-sphere:

\begin{lem}\label{lem:3sph8pt}
$K$ contains no $8A$-vertices $x$, such that $\Si_x K$ contains a 
singular $3$-sphere $S$ with the following properties:
$S$ contains a pair of antipodal $2$-vertices 
$\xi_1, \xi_2$, such that $\Si_{\xi_i}S$ is a singular
$2$-sphere containing three pairwise orthogonal
$7$-vertices.
Furthermore, all $7$-vertices in $S$ are 
$\pithird$-extendable to $8A$-vertices. 
\end{lem}

Notice that all $7$-vertices in $S$ are adjacent to $\xi_i$
for some $i=1,2$. Indeed, a segment in $\Si_x K$ (of type $E_7$) 
connecting a 2- and a 7-vertex at distance $\leq\pihalf$
is of type 27 or 217. This last segment cannot occur 
between $\xi_i$ and a 7-vertex in $S$ 
because $\Si_{\xi_i} S$ does not contain 1-vertices
(the 2-sphere $\Si_{\xi_i} S$ is tessellated by simplices with
vertices of type 5,6,7).
Observe also, that the link $\Si_\lambda S$ of a 7-vertex $\lambda\in S$ contains a 
singular circle of type 2626262: suppose w.l.o.g.\ that $\lambda$ is adjacent
to $\xi_1$, then $\ora{\xi_1\lambda}$ is contained in a circle in $\Si_{\xi_1} S$
of type 767676767. In particular $\Si_{\ora{\lambda\xi_1}}\Si_\lambda S$ contains a pair
of antipodal 6-vertices. It follows that the antipodal directions 
$\ora{\lambda\xi_1}$ and $\ora{\lambda\xi_2}$ are contained in a singular circle in 
$\Si_\lambda S$ of type 2626262.

\bigskip\no
\emph{Proof of Lemma~\ref{lem:3sph8pt}.}
Suppose there are such $8A$-vertices.
Let $x_1,x_2,x_3\in K$ be such $8A$-vertices as in configuration $\ast$,
and let $S_{x_i}\subset\Sigma_{x_i}K$ denote 
the corresponding 3-spheres in their links. 
Let $y_3,z_1\in K$ be as in the notation of the configuration $\ast$.
Suppose that there is a 7-vertex $\xi\in S_{x_1}\subset\Si_{x_1}{K}$, such that
$d(\xi,\zeta)=\arccos(-\frac{1}{\sqrt{3}})$ for $\zeta:=\ora{x_1z_1}$. The segment
$\xi\zeta$ is of type 7672. 
By assumption, there exists an 8A-vertex $x\in K$,
such that $d(x_1,x)=\pithird$ and $\ora{x_1x}=\xi$.
Under these circumstances we obtain the configuration $\aast$.
 We use the same notation as
in the configuration $\aast$.
Let $\alpha_i\in \Si_{x_3}K$ for $i=1,\dots,4$ be the directions
$\ora{x_3x_1}$, $\ora{x_3x}$, $\ora{x_3w}$ and $\ora{x_3y_3}$. 
Let $\beta:=\ora{x_3z}$. Then the 7-vertices $\alpha_i$ are adjacent to the 
2-vertex $\beta$. And the directions $\ora{\beta\alpha_i}$ lie on a circle $\kappa$ of
type 767676767 contained in $\Si_\beta\Si_{x_3} K$.

\parpic{\includegraphics[scale=0.55]{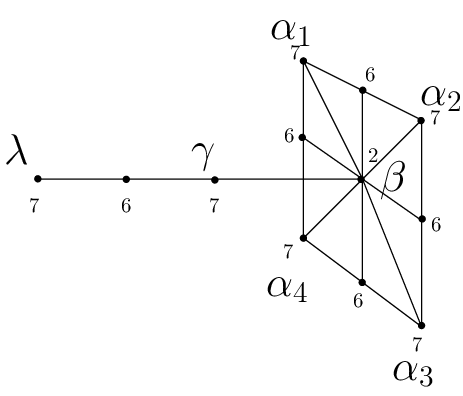}}
Suppose again that there is a 7-vertex $\lambda$
 in the 3-sphere $S_{x_3}\subset\Si_{x_3}K$,
such that $d(\beta,\lambda)=\arccos(-\frac{1}{\sqrt{3}})$. So the segment 
$\beta\lambda$ is of type 2767. 
Recall that the 7-vertices $\alpha_i$ are 
$\2pithird$-extendable and $\lambda$ is $\pithird$-extendable to an $8A$-vertex,
so they cannot be antipodal.
It follows that $\angle_\beta(\lambda,\alpha_i)=\pihalf$ and
$d(\alpha_i, \lambda)=\arccos(-\frac{1}{3})$.
The segments $\alpha_i\lambda$ are of type 727.
Let $\gamma\in \Si_{x_3}K$ be the 7-vertex
on the interior of the segment $\beta\lambda$. 
Then $\gamma$ is the center of an equilateral spherical triangle 
$CH(\lambda,\alpha_1,\alpha_3)$ with sides of type 727.
We are now in the situation of the configuration
$\aast$ (compare with the triangle $CH(\xi,\xi_2,\xi_3)$ in the definition
of the configuration $\aast$). It follows that 
$\gamma$ is $\2pithird$-extendable.

\parpic{\includegraphics[scale=0.4]{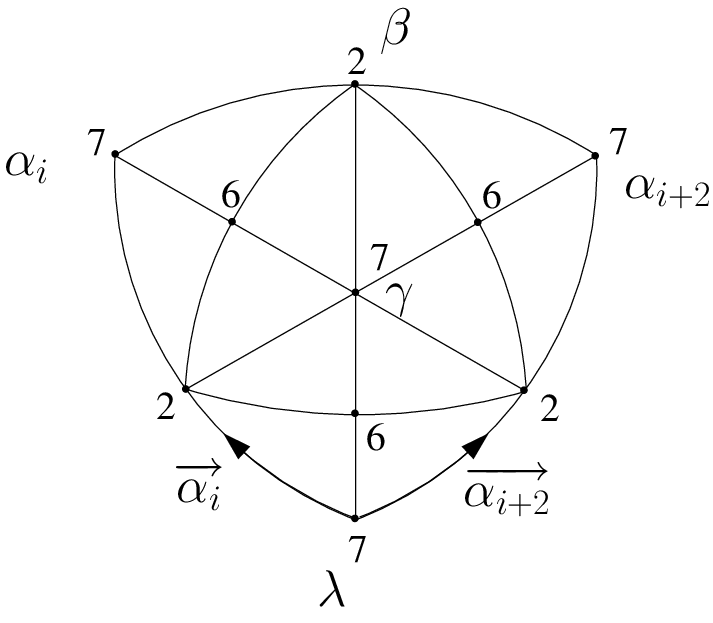}}
The convex hull $CH(\kappa,\ora{\beta\lambda})$ is a 2-dimensional hemisphere
centered at $\ora{\beta\lambda}$. Hence, 
$\Si_{\ora{\beta\lambda}}\Si_\beta \Si_{x_3}K$ 
(of type \hpic{\includegraphics[scale=0.4]{E6link2dynk}})
contains a circle of type 656565656.
This is equivalent to $\Si_{\ora{\lambda\beta}}\Si_\lambda \Si_{x_3}K$
(of type \hpic{\includegraphics[scale=0.4]{E6link6dynk}})
containing a circle of type 232323232. Note that the 2-vertices on this
circle correspond to the 2-vertices $m(\lambda,\alpha_i)\in \Si_{x_3}K$
(consider the equilateral spherical triangles 
$CH(\lambda,\alpha_i,\alpha_{i+2})$ with sides of type 727).
Let $\ora{\alpha_i}:=\ora{\lambda\alpha_i}\in \Si_\lambda \Si_{x_3}K$.

Recall that the link $\Si_\lambda S_{x_3}$ contains a circle $c$ of type 2626262
and notice that $\ora{\lambda\beta}$ cannot be antipodal to any of the 2-vertices
on this circle: otherwise, we find a 7-vertex
in the 3-sphere $S_{x_3}$ antipodal to $\gamma$. 
This cannot happen, because $\gamma$
is $\2pithird$-extendable and the 7-vertices in $S_{x_3}$ are $\pithird$-extendable
to $8A$-vertices.
It is also clear that $\ora{\lambda\beta}$ cannot have distance $<\pihalf$ to the three
6-vertices on the circle $c$, otherwise $c$ would be contained in a ball
centered at $\ora{\lambda\beta}$ with radius $<\pihalf$, but this is not possible since
$diam(c)=\pi$.

\parpic{\includegraphics[scale=0.55]{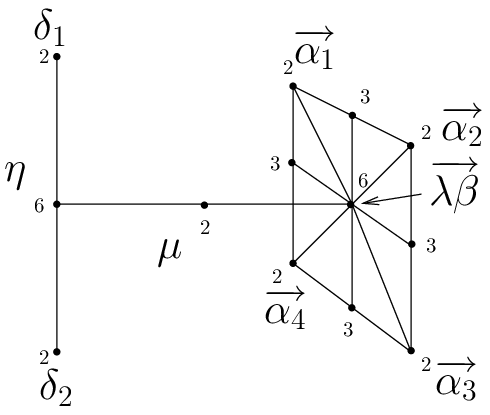}}
Therefore we can find a 6-vertex $\eta$ on the circle $c\subset \Si_\lambda S_{x_3}$, 
such that 
$d(\eta,\ora{\lambda\beta})\geq\pihalf$. Hence,
$d(\eta,\ora{\lambda\beta})=\2pithird$ and
the segment $\eta\ora{\lambda\beta}$ is of type 626.
Let $\mu:=m(\eta,\ora{\lambda\beta})$.
Let also $\delta_1,\delta_2$ be the two 2-vertices in the circle 
$c\subset\Si_\lambda \Si_{x_3}K$ adjacent to $\eta$.

We have already seen, that $\ora{\lambda\beta}$ cannot be antipodal to $\delta_i$.
This implies that $\angle_\eta(\delta_i,\mu)=\pihalf$ and these angles are of type 232.
It follows that
$\Si_{\ora{\eta\mu}}\Si_\eta\Si_\lambda\Si_{x_3}K$
(of type $D_4$: 
\hpic{\includegraphics[scale=0.4]{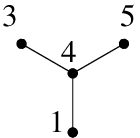}})
 contains a pair of antipodal 3-vertices.
On the other hand, if $\eta$ is antipodal to some $\ora{\alpha_i}$, then $\alpha_i\in \Si_{x_3}K$
has an antipode in $S_{x_3}$, but this cannot happen either, because $\alpha_i$ is 
$\2pithird$-extendable in $K$. Therefore $\angle_{\ora{\lambda\beta}}(\mu,\ora{\alpha_i})=\pihalf$
and these angles are of type 232.
It follows that 
$\Si_{\ora{{\scriptscriptstyle{\ora{\lambda\beta}}}\mu}}\Si_{\ora{\lambda\beta}}\Si_\lambda\Si_{x_3}K$
 contains a singular circle of type 343434343. 
This in turn implies, 
that $\Si_{\ora{\eta\mu}}\Si_\eta\Si_\lambda\Si_{x_3}K$ contains
a singular circle of type 141414141,
 because the antipode of a 3- (4)-vertex in $\Si_\mu \Si_\lambda\Si_{x_3} K$,
of type \hpic{\includegraphics[scale=0.4]{E6link2dynk}}, adjacent to
$\ora{\mu{\scriptstyle \ora{\lambda\beta}}}$ is a 1- (4)-vertex adjacent to $\ora{\mu\eta}$.
We apply now Lemma~\ref{lem:Dn_n-2sph} to conclude that 
$\Si_{\ora{\eta\mu}}\Si_\eta\Si_\lambda\Si_{x_3}K$ contains a wall. 
Hence $\Si_\mu\Si_\lambda\Si_{x_3}K$ contains a wall.

\parpic{\includegraphics[scale=0.45]{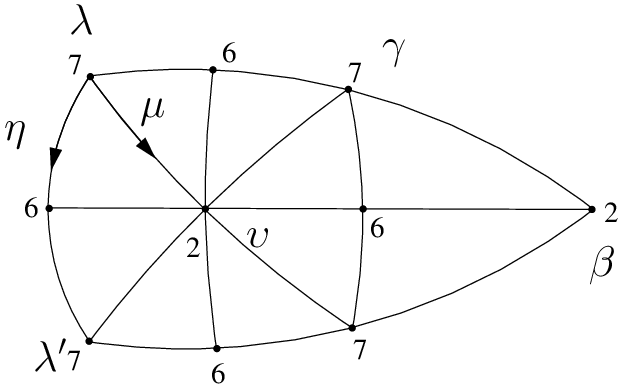}}\label{lem:3sph8ptextend}
Let $\lambda'\in S_{x_3}$ be the 7-vertex at distance $\ac(\third)$ to $\lambda$,
 so that $\ora{\lambda\lambda'}=\eta$.
By considering the spherical triangle 
$CH(\lambda,\lambda',\beta)\subset \Si_{x_3} K$
we deduce that $\mu$ is $\ac(-\third)$-extendable in $\Si_{x_3}K$.
Let $\upsilon$ be the 2-vertex in $\Si_{x_3}K$ adjacent to $\lambda$
with $\ora{\lambda\upsilon}=\mu$.
It follows that $\Si_\upsilon\Si_{x_3}K$ contains a wall.

\parpic[r]{\includegraphics[scale=0.45]{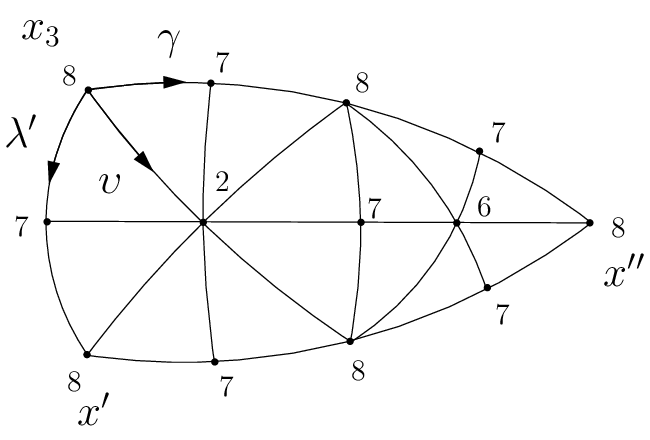}}
Recall that $\gamma$ is $\2pithird$-extendable and let $x''\in K$ be an 8-vertex 
with $d(x_3,x'')=\2pithird$ and $\ora{x_3x''}=\gamma$. Since $\lambda'\in S_{x_3}$,
it is $\pithird$-extendable. Let $x'\in K$ be an 8-vertex, so that $d(x_3,x')=\pithird$
and $\ora{x_3x'}=\lambda'$.
Consider the spherical triangle $CH(x_3,x'',x')$. 
One sees that $\upsilon$ is $\pihalf$-extendable
in $K$, thus we have found a 2-vertex in $K$, whose link contains a wall, contradicting
Lemma~\ref{lem:5sph2pt}.

So it follows that $d(\beta,\lambda)\leq\pihalf$ for all
7-vertices $\lambda\in S_{x_3}$.
Since $S_{x_3}$ is the convex hull of the 7-vertices contained in it,
this implies that $d(\beta,S_{x_3})\equiv\pihalf$
and $s:=\Si_\beta CH(\beta,S_{x_3})$ is a 3-sphere.
Let $\theta\in S_{x_3}\subset\Si_{x_3}K$ be a 2-vertex,
so that $\Si_\theta S_{x_3}$ is a 2-sphere spanned by three pairwise
orthogonal 7-vertices (compare with the description of the 3-sphere $S_{x_3}$).
The segment $\theta\beta$ is of type 262.

Notice that $d(\beta,S_{x_3})\equiv\pihalf$ implies that
$d(\ora{\theta\beta},\Si_{\theta}S_{x_3})\equiv\pihalf$.
It follows that $\Si_{\ora{\theta\beta}}\Si_{\theta} CH(\beta,S_{x_3})$ 
(subset of a building of type \hpic{\includegraphics[scale=0.4]{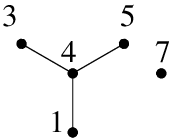}}) is a 2-sphere.
Notice that in the building $\Si_\theta \Si_{x_3} B$
of type \hpic{\includegraphics[scale=0.4]{E8link28dynk}};
two 7-, 6-vertices at distance $\pihalf$ are joined by a segment of type 756.
This implies that $\Si_{\ora{\theta\beta}}\Si_{\theta} CH(\beta,S_{x_3})$ 
is spanned by three pairwise
orthogonal 5-vertices. Such a 2-sphere in
the Coxeter complex of type \hpic{\includegraphics[scale=0.4]{E6link26dynk}}
is not a subcomplex, thus, its simplicial convex hull
is a 3-sphere. 
Therefore the 3-sphere $s\subset\Si_\beta\Si_{x_3}K$ is not a subcomplex and its
simplicial convex hull is a wall.
Recall that $\beta\in \Si_{x_3}K$ is $\pihalf$-extendable in $K$, 
hence, there are 2-vertices
in $K$, with links containing a wall. We have now a 
contradiction to Lemma~\ref{lem:5sph2pt}.

It follows that our first assumption, that 
there is a 7-vertex $\xi\in S_{x_1}\subset\Si_{x_1}{K}$, such that
$d(\xi,\zeta)=\arccos(-\frac{1}{\sqrt{3}})$ cannot occur. Thus,
$d(\zeta,S_{x_1})\equiv\pihalf$ and repeating the previous argument, we can
see that $\Si_\zeta\Si_{x_1}K$ contains a wall. Hence, $\Si_{z_1} K$ contains a wall,
contradicting again Lemma~\ref{lem:5sph2pt}.
\qed

\begin{lem}\label{lem:4sph2pta}
 Let $x\in K$ be a $2$-vertex, 
such that $\Si_x K$ contains a singular $4$-sphere $S$ of type
$757$ or $\pithird$ (cf. Section~\ref{sec:Dncoxeter}, page~\pageref{pag:spheresDn}). 
Then,
the $8$-vertices in $S\subset \Si_x K$ are not $\pihalf$-extendable
 and there are no $8$-vertices in $K$
at distance $\frac{3\pi}{4}$ to $x$.
In particular, $x$ is a $2A$-vertex, and all $8$-vertices in $S$ are directions to
$8A$-vertices in $K$ adjacent to $x$.
\end{lem}
\proof
We show first by contradiction that there is no 8-vertex in $S$ that is
$\pihalf$-extendable.
Suppose there is an 8-vertex in $S$ that is $\pihalf$-extendable. 
This means that there
is a 2-vertex $y\in K$ at distance $\pihalf$ to $x$, such that the segment
$xy$ is of type 282 and $\ora{xy} \in S$. In particular $\Si_{\ora{xy}}\Si_x S$ is a 
singular 3-sphere. 
This implies for the 8-vertex $z:=m(x,y)$, 
that its link $\Si_z K$ contains a 4-sphere.
By Lemma~\ref{lem:int8pts}, $dim(K)\geq 6$,
otherwise $dim(K)=5$ and $z$ would be an $8I$-vertex. 
In particular, $\Si_x K$ has dimension $\geq 5$ and it must
contain a 5-dimensional
hemisphere $h$ bounded by $S$.

The hemisphere $h$ is the intersection of 
a {\em wall} and a {\em root} in a building of type $D_7$
with Dynkin diagram \hpic{\includegraphics[scale=0.4]{E8link2dynk}}.
Recall the description of hemispheres of codimension 1 in Section~\ref{sec:Dncoxeter}.
If $S$ is of type 757, then $h$ is centered at a 7-vertex $\alpha$
and $\Si_\alpha h$ is a wall of type 5. In particular $\Si_\alpha h$
contains a pair of antipodal 8-vertices.
If $S$ is of type $\pithird$, then $h$ is centered at point contained in the interior
of an edge of type $86$. In particular, the 8-vertex of this edge is contained in $h$.
In both cases $h$ contains an 8-vertex $\eta$ in its interior
(notice that this is not true for a hemisphere bounded by a singular
4-sphere of type 787).
It is clear that $d(\eta,\ora{xy})=\pihalf$ and the segment is of type 878. 
The midpoint 
$\zeta:=m(\eta,\ora{xy})$ is also in the interior of $h$, and in particular, 
$\Si_\zeta\Si_x K$
contains a wall of type 5, that is, a wall containing a pair of antipodal 8-vertices.

\parpic{\includegraphics[scale=0.7]{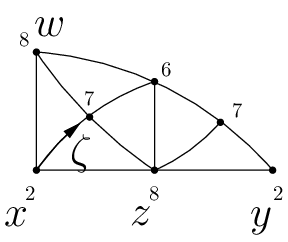}}
Let $w\in K$ be the 8-vertex in $K$ adjacent to $x$ with $\eta=\ora{xw}$. 
Then if we consider
the spherical triangle $CH(x,y,w)$, we see that $\zeta$ is extendable 
to a segment of type 276 in $K$. Therefore we find 
7-vertices in $K$, whose links in $K$ contain a wall of type 1.
A contradiction to Lemma~\ref{lem:5sph7pt}.
This finishes te proof that 
the 8-vertices in $S\subset \Si_x K$ are not $\pihalf$-extendable.

We prove now that $x$ is a $2A$-vertex. Otherwise we find an antipode 
$\wht x\in K$ of $x$ and
a segment connecting $x$ and $\wht x$ with initial direction an 
8-vertex in $S$ is of type
28282. It follows that the 8-vertices in $S$ are $\3piquart$-extendable. 
A contradiction.

Now we show the last statement of the lemma.
Let $u\in K$ be an 8-vertex adjacent to $x$, such that $\ora{xu}\in S$
and suppose that $u$ has an
antipode $\wht u \in K$. Let $c$ be the segment connecting $u$ and $\wht u$ through
$x$. It is of type 82828. Since the direction 
$\ora{x\wht u}$ has an antipode in $S$, namely
$\ora{xu}$, it follows that the 8-vertex 
$\ora{x\wht u}$ lies in a sphere $S'\subset \Si_x K$
of the same type as $S$. Hence $\ora{x\wht u}$ cannot be $\pihalf$-extendable, but
the segment $x\wht u$ is of type 2828. A contradiction.
Thus, all $8$-vertices in $S$ are directions to $8A$-vertices in $K$ adjacent to $x$.

\parpic{\includegraphics[scale=0.55]{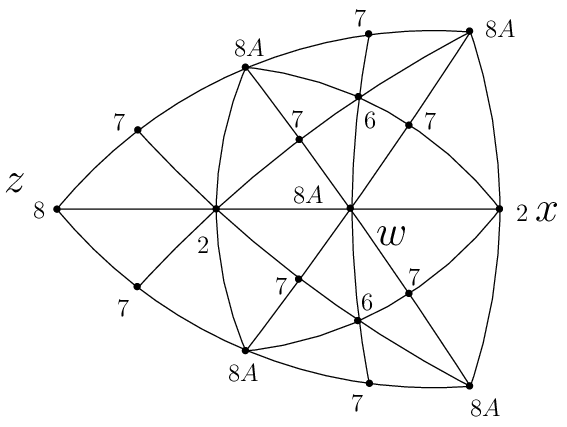}}
For the second assertion, 
suppose there is an 8-vertex $z\in K$ with $d(x,z)=\3piquart$. Since
all 8-vertices in $S$ correspond to $8A$-vertices in $K$, 
the 8-vertex $\ora{xz}$ must be orthogonal
to the 8-vertices in $S$.
In both cases (of type 757 or $\pithird$), 
$S$ contains a singular 2-sphere spanned
by three pairwise orthogonal 8-vertices (cf. Section~\ref{sec:Dncoxeter}).
This implies that $\Si_{\ora{xz}}\Si_x K$ contains a 2-sphere
spanned by three pairwise orthogonal 7-vertices. 
Let $w$ be the 8-vertex in $xz$ adjacent to $x$.
Recall that $x$ is a $2A$-vertex, therefore $w$ is an $8A$-vertex.
Then $\Si_w K$ contains a 3-sphere as described
in the statement of Lemma~\ref{lem:3sph8pt}
by Lemma~\ref{lem:sphinlink} and
the fact that $\Si_{\ora{xz}}\Si_x K$ contains a 2-sphere
spanned by three pairwise orthogonal 7-vertices. 
It also follows that the 7-vertices
in this sphere are $\pithird$-extendable to $8A$-vertices
(see figure), 
contradicting Lemma~\ref{lem:3sph8pt}. 
\qed

\begin{lem}\label{lem:4sph2ptb}
$K$ contains no $2$-vertices $x$, such that $\Si_x K$ contains a singular
$4$-sphere $S$ of type $\pithird$.
\end{lem}
\proof
Let $x$ be such a $2$-vertex. 
It follows from Lemma~\ref{lem:4sph2pta} that $x$ is a $2A$-vertex
and $rad(x,8\text{-vert. in } K)\leq\arccos(-\frac{1}{2\sqrt{2}})$.
Let $x'\in G\cdot x$ with $d(x,x')=diam(G\cdot x)$.

\parpic[r]{\includegraphics[scale=0.45]{lem6sph2pt}}
Case 1: $diam(G\cdot x)=\arccos(-\frac{3}{4})$.
The segment $xx'$ is of type $21812$. 
Let $y_1,y_2\in K$ be 8-vertices adjacent to $x$, such that $y_1xy_2$ is a segment of
type 828. These vertices can be found, because $\Si_x K$ contains a singular
$4$-sphere $S$ of type $\pithird$, which contains antipodal 8-vertices.
We may assume that $\angle_x(y_1,x')\geq\pihalf$. This implies that the angle
$\angle_x(y_1,x')=\ac(-\frac{1}{\sqrt{7}})$ 
and this angle is of type 831, 
because $\Si_x B$ is a building of type $D_7$. 
$CH(y,x,x')$ is a spherical triangle, therefore we can compute that 
$d(y_1,x')=\frac{3\pi}{4}$.
A contradiction to $rad(x,8\text{-vert. in } K)\leq\arccos(-\frac{1}{2\sqrt{2}})$.

Case 2: $diam(G\cdot x)=\2pithird$.
The segment $xx'$ is of type 26262.
If there is an $8$-vertex $y$ adjacent to $x$, such that 
$d(\ora{xx'},\ora{xy})>\pihalf$, then $d(x',y)=\3piquart$
(see the figure in the proof of Lemma~\ref{lem:6sph2pt}, Case 2),
which is not possible.
Therefore we see that the distance between 
the 6-vertex $\ora{xx'}$ and all the 8-vertices
in $S$ must be $\pihalf$ because they come in pairs of antipodes.
If $S'\subset S$ is
the 3-sphere spanned by the 8-vertices in $S$, 
then $d(\ora{xx'},S')\equiv \pihalf$.
It follows that $\Si_{\ora{xx'}}\Si_x K$ 
(of type \hpic{\includegraphics[scale=0.4]{E8link26dynk}})
contains a 3-sphere containing four pairwise orthogonal 5-vertices
(the directions of the segments from $\ora{xx'}$ to the 8-vertices
in $S'$),
this sphere is an apartment in the $D_4$-factor. 
Let $y:=m(x,x')$. Then the link $\Si_{\ora{yx}} \Si_y K$ 
(again of type \hpic{\includegraphics[scale=0.4]{E8link26dynk}})
contains also an apartment in the $D_4$-factor. This is a 3-sphere
containing a simplex of type 1345.
This implies that the link $\Si_y K$ contains a singular 4-sphere $S_y$
containing a simplex of type 13456.
A singular 4-sphere containing a simplex of type 13456 in a building of type
\hpic{\includegraphics[scale=0.4]{E8link2dynk}} is precisely a sphere
of type $\pithird$.
Hence $S_y$ is of type $\pithird$ 
and the 6-vertices $\ora{yx}$, $\ora{yx'}$ are orthogonal to the 3-sphere
$S_y'\subset S_y$ spanned by the 8-vertices in $S_y$.
To see this consider the vector space model of the Coxeter 
complex of type $D_7$ introduced
in the Appendix~\ref{app:coxeter}.
The sphere $S_y$ can be identified with the sphere 
$\{x_5=x_6=x_7\}\cap S^6 \subset\R^7$ and $S_y'$, with the sphere
$\{x_5=x_6=x_7=0\}\cap S^6$.
A 5-vertex in $S_y'$ is of the form $(\pm1,\dots,\pm1,0,0,0)$ and a 6-vertex
orthogonal to this sphere must be of the form $(0,\dots,0,\pm1,\pm1,\pm1)$.
Hence, a 5-vertex in $S_y'$ and a 6-vertex orthogonal to $S_y'$
are connected by a segment of type 536 or 516.

We have seen that $y$ is a 2-vertex as in the statement of the lemma.
Therefore we can proceed as in the beginning of the proof, 
to conclude that
$rad(y,8\text{-vert. in } K)\leq\arccos(-\frac{1}{2\sqrt{2}})$ 
and $diam(G\cdot y)\leq\2pithird$.
If $diam(G\cdot y)=\arccos(-\quart)$, then we are in the situation
of the next case below (Case 3).
Hence, we may assume again that $diam(G\cdot y)=\2pithird$.
Let $y'\in G\cdot y$ have distance 
$\2pithird$ to $y$. It follows as above, that $\Si_{\ora{yy'}}\Si_y K$ contains
an apartment in the $D_4$-factor.

Let $\xi,\xi'\in S_y'$ be antipodal 5-vertices. 
The vertices $\ora{yx}$, $\ora{yx'}$, $\xi$ and $\xi'$ lie on a singular circle
of type 635161536 contained in $S_y$.
The link $\Si_\xi\Si_y B$ is of type
$A_3\circ A_3$ and has Dynkin diagram {\includegraphics[scale=0.4]{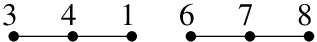}.
Notice that $\Si_\xi S_y'$ is an apartment in the second $A_3$-factor. 
Therefore the second factor in the spherical join splitting of $\Si_\xi\Si_y K$ is
a subbuilding.
Since
$rad(y,8\text{-vert. in } K)\leq\arccos(-\frac{1}{2\sqrt{2}})$,
this implies as above that $d(\ora{yy'},S_y')\equiv \pihalf$. In particular,
$d(\ora{yy'},\xi)=\pihalf$ and the direction $\ora{\xi\scriptstyle{\ora{yy'}}}$
must be orthogonal to the 2-sphere $\Si_\xi S_y'$. 
Recall that this sphere is an apartment
in the second $A_3$-factor. Thus $\ora{\xi\scriptstyle{\ora{yy'}}}$
must lie on the \includegraphics[scale=0.4]{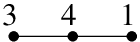} -factor of $\Si_\xi\Si_y K$.

\parpic{\includegraphics[scale=0.6]{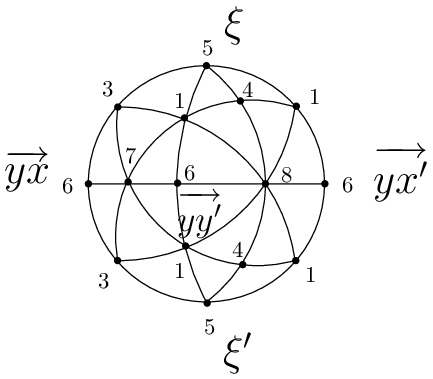}}
It follows from this that the segments $\xi\ora{yy'}$ and $\xi'\ora{yy'}$ must be of type 536
or 516. Further, since $d(\xi,\ora{yy'})+d(\xi',\ora{yy'})=d(\xi,\xi')=\pi$, the segments
are of the same type. 
Observe also, that $\ora{yy'}$ cannot be antipodal to $\ora{yx}$
or $\ora{yx'}$, otherwise the $2A$-vertex $y'$ 
would be antipodal to $x$ or $x'$.
Suppose w.l.o.g.\ that the segments $\xi\ora{yy'}\xi'$ and
$\xi\ora{yx'}\xi'$ are of type 51615.
This implies that the segment $\xi\ora{yx}\xi'$ is of type 53635.
Since $\ora{yy'}$ is not antipodal to $\ora{yx}$, then the directions
$\ora{\xi\ora{yx}}$ and $\ora{\xi\ora{yy'}}$ of type 3 and 1, respectively,
cannot be antipodal, thus, they are adjacent
(recall that these directions lie in a building of type 
\includegraphics[scale=0.4]{E8link25678dynk}).
 This implies that the segment
$\ora{yx}\ora{yy'}$ has length $\ac(\third)$ and is of type 676.
It also follows that $\ora{yy'}$ lies on a segment of length $\pi$
and type 67686 connecting $\ora{yx}$ and $\ora{yx'}$.
Therefore, the segment $\ora{yx'}\ora{yy'}$ has length $\ac(-\third)$ and is of type 686.
Hence, $\Si_{\ora{yy'}}\Si_y K$ contains antipodal 7- and 8-vertices, that is,
it contains a wall in the $A_2$-factor.
Together with the apartment in the $D_4$-factor (compare with the beginning of Case 1), 
this implies that 
the link $\Si_{\ora{yy'}}\Si_y K$
contains a wall. It follows that the link in $K$ of the 2-vertex
$m(y,y')$ contains a wall, contradicting Lemma~\ref{lem:5sph2pt}.

Thus, $diam(G\cdot y)=\ac(-\quart)$ and by relabelling $y$ by $x$ 
we have reduced the possibilities to the following case.

Case 3: $diam(G\cdot x)=\arccos(-\quart)$. 
The simplicial convex hull $C$ of $xx'$ is
2-dimensional. Let $y,y'\in C$ be the 8-vertices adjacent to $x$ and $x'$, 
respectively (see table of distances between two 2-vertices
in Section~\ref{sec:E8coxeter}).
If $\ora{xy}$ has an antipode in $\Si_x K$, then there would be an 8-vertex in $K$
at distance $\3piquart$ to $x'$, but this is not possible 
(cf. proof of Lemma~\ref{lem:6sph2pt}, Case 3).
It follows that $\ora{xy}$ has distance $\pihalf$ to all 8-vertices in $S'$,
where $S'\subset S$ is the 3-sphere
spanned by the 8-vertices in $S$,
thus, $d(\ora{xy},S')\equiv \pihalf$. 
The link $\Si_{\ora{xy}}CH(\ora{xy},S')$ is a 3-sphere containing
four pairwise orthogonal 7-vertices. 

Let $w\in C$ be the 7-vertex $m(x,x')$
(see table of distances between two 2-vertices
in Section~\ref{sec:E8coxeter} or the figure below).
Let $x''\in G\cdot x$ with $d(w,x'')>\pihalf$,
which exists by the first property of a counterexample explained
at the beginning of Section~\ref{sec:centerconj}. 
The possible distances between 2- and 7-vertices in the Coxeter complex
of type $E_8$ are of the form $\ac(-\frac{k}{2\sqrt{6}})$ for $k$ an integer 
(this can be deduced from the table of 2- and 7-vertices 
in Appendix~\ref{app:E8coxeter}).
Notice that $d(x,w)=d(w,x')=\ac(\frac{3}{2\sqrt{6}})$.
Triangle comparison for the triangle $(x,x',x'')$ 
and $diam(G\cdot x)\leq \ac(-\quart)$
imply that $d(x'',w)=d(x'',m(x,x'))\leq\ac(-\frac{1}{\sqrt{6}})$. 
If $d(w,x'')=\ac(-\frac{1}{\sqrt{6}})$, then by rigidity,
$CH(x,x',x'')$ is an equilateral spherical triangle with side lengths $\ac(-\quart)$.
In particular $d(x,x'')=\ac(-\quart)$ and $\angle_x(x',x'')>\pihalf$.

If $d(w,x'')=\ac(-\frac{1}{2\sqrt{6}})$, we may assume w.l.o.g.\ that
$\angle_w(x,x'')\geq\pihalf$. 
This implies that $d(x,x'')\geq\ac(-\frac{1}{8})$, i.e.\ $d(x,x'')=\ac(-\quart)$. 
Again by triangle comparison and
$\angle_w(x,x'')\geq\pihalf$ we want to
see that $CH(x,w,x'')$ must be a spherical triangle:
let $\tilde x,\tilde{x''}$ be 2-vertices and let $\tilde w$ be a 7-vertex in the Coxeter 
complex of type $E_8$, such that $d(\tilde x,\tilde w)=d(x,w)=\ac(\frac{3}{2\sqrt{6}})$, 
$d(\tilde w,\tilde{ x''})=d(w,x'')=\ac(-\frac{1}{2\sqrt{6}})$ and 
$\angle_w(x,x'')=\angle_{\tilde w}(\tilde x,\tilde{x''})$.
By triangle comparison, $d(\tilde x, \tilde{x''})\leq d(x,x'')=\ac(-\quart)$,
but since the angle $\angle_{\tilde w}(\tilde x,\tilde{x''})=\angle_w(x,x'')\geq\pihalf$, 
then $d(\tilde x, \tilde{x''})>\pihalf$.
It follows that $d(\tilde x, \tilde{x''})=\ac(-\quart)= d(x,x'')$ and by rigidity
$CH(x,w,x'')$ is a spherical triangle.
We can now compute that $\angle_x(x',x'')=\ac(-\frac{1}{15})>\pihalf$.

\parpic{\includegraphics[scale=0.6]{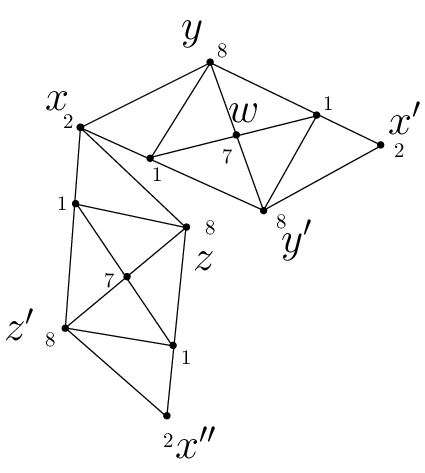}}
Let $C'$ be the 2-dimensional simplicial convex hull of $xx''$ and let
$z,z'\in C'$ be the 8-vertices adjacent to $x$ and $x''$.
By considering the spherical triangle $CH(x,x',y)$, we can compute
$\angle_x(y,x')=\ac(\frac{3}{\sqrt{15}})<\piquart$.
Then we can see that, if $\ora{xy}=\ora{xz}$, 
it follows $\angle_x(x',x'')<\pihalf$, thus $\ora{xy}\neq\ora{xz}$.
They cannot be antipodal either, because
$\ora{xy}$ has no antipodes in $\Si_x K$ (compare with the beginning of Case 2). 
Hence, the segment $\ora{xy}\ora{xz}$ has length $\pihalf$ and 
is of type 878. 

Let $\xi\in \Si_x K$ be the 7-vertex $m(\ora{xy},\ora{xz})$.
Notice that as for $\ora{xy}$, it also holds $d(\ora{xz},S')\equiv \pihalf$.
This implies that the convex hull of $S'$ and the segment $\ora{xy}\ora{xz}$
is isometric to the spherical join $S'\circ \ora{xy}\ora{xz}$.
In particular, $d(\xi,S')\equiv \pihalf$.
Notice that in a building of type $D_7$ with Dynkin diagram
\hpic{\includegraphics[scale=0.4]{E8link2dynk}},
a 7- and an 8-vertex at distance $\pihalf$ are joined by a segment of
type 768.
It follows that 
$\Si_\xi\Si_x K$ 
(of type \hpic{\includegraphics[scale=0.4]{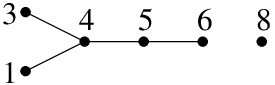}})
contains a 3-sphere spanned
by four pairwise orthogonal 6-vertices. This 3-sphere is not simplicial,
and its simplicial convex hull is an apartment in the $D_5$-factor
of $\Si_\xi\Si_x K$.
Since $\{\ora{xy},\ora{xz}\}$ is an apartment in 
the $A_1$-factor of $\Si_\xi\Si_x K$, 
it follows
that $\Si_\xi\Si_x K$ contains an apartment. In particular $\xi$ is an interior
7-vertex in $\Si_x K$.

We can also see, that if both 1-vertices $\ora{xy'}$ and $\ora{xz'}$
are adjacent to $\xi$, then $\angle_x(x',x'')<\pihalf$,
because in this case 
$d(\xi,\ora{xw})=d(\xi,\ora{xx''})=
\ac(\frac{2\sqrt{2}}{\sqrt{15}})<\piquart$
(just consider the spherical triangle $CH(\ora{xy},\ora{xw},\xi)$
with sides $d(\ora{xy},\ora{xw})=\ac(\frac{3}{\sqrt{15}})$, $d(\ora{xy},\xi)=\piquart$
and angle $\angle_{\ora{xy}}(\ora{xw},\xi)=\ac(\frac{1}{\sqrt{6}})$).

\parpic{\includegraphics[scale=0.65]{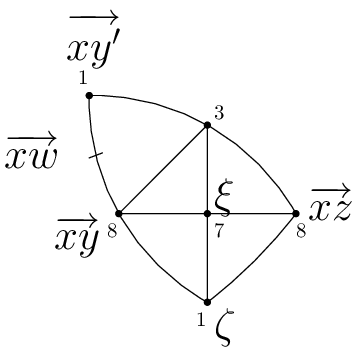}}
Therefore
w.l.o.g.\ $\ora{xy'}$ is not adjacent to $\xi$, but since both are
adjacent to $\ora{xy}$, the angle $\angle_{\ora{xy}}(\xi,\ora{xy'})$
must be of type 731, because $\Si_{\ora{xy}}\Si_x B$ is of type
$D_6$ with Dynkin diagram \hpic{\includegraphics[scale=0.4]{E8link28dynk}}.
Now recall that $\xi$ is an interior vertex in $\Si_x K$,
this implies that we can find a 1-vertex $\zeta\in \Si_x K$, so that
$\ora{\scriptstyle{xy'}}\ora{xy}\zeta$ is a segment of type 181.
Thus, the link $\Si_{\ora{xy}}\Si_x K$ (of type $D_6$)
contains a pair of antipodal 1-vertices and
a 3-sphere spanned by four pairwise orthogonal 7-vertices
(compare with the beginning of Case 2). 
We can apply Lemma~\ref{lem:Dn_n-2sph} to see that $\Si_{\ora{xy}}\Si_x K$
contains a wall.
By Lemma~\ref{lem:Dn_wall},
$\Si_{\ora{{\scriptscriptstyle \ora{xy}\ora{xw}}}}\Si_{\ora{xy}}\Si_x K$
contains at least a wall. This implies
that $\Si_{\ora{xw}}\Si_x K$ contains a wall and
$\Si_w K$ contains a wall of type 1, 
contradicting Lemma~\ref{lem:5sph7pt}.
\qed

We introduce now two new $G$-invariant properties of $8A$-vertices in $K$

\parpic(6cm,3cm)[r]{\includegraphics[scale=0.5]{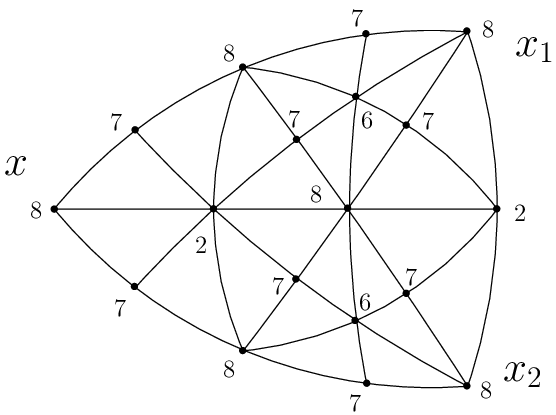}}
Let $x\in K$ be an $8A$-vertex. We say that 
\textbf{$x$ has the {\em property $T$}}, if there is no
spherical triangle in $K$ with $8A$-vertices $x, x_1$ and $8$-vertex $x_2$,
with side lengths $d(x,x_i)=\2pithird$,
$d(x_1,x_2)=\pihalf$, and such that the direction
$\ora{xx_2}$ is $\2pithird$-extendable to an $8A$-vertex in $K$.
This last assumption is fulfilled if e.g.\ $x_2$ is also an $8A$-vertex.
\bigskip

We analyze further the configuration $\ast$ in the case where we
have $8T$-vertices.
Let $x_1,x_2,x_3\in K$ be $8T$-vertices as in configuration $\ast$. 
If $\angle_{y_3}(x_3,x_2)=\arccos(\third)$,
then the convex hull $CH(x_3,y_3,m(y_3,x_2))$ is a
spherical triangle but it is not simplicial.
Its simplicial convex hull is a spherical
triangle with vertices $x_3,y_3$ and an 8-vertex $x_2'$
which is not necessarily equal to $x_2$.
This triangle has
sides $y_3x_3$, $x_3x_2'$ and $x_2'y_3$
of type 87878, 828 and 878, respectively, and $m(y_3,x_2)=m(y_3,x_2')$.
It follows that the simplicial convex hull of $x_1,m(y_3,x_2),x_3$ 
is a spherical triangle in $K$ 
as ruled out by the property $T$,
hence the property $T$ implies that 
$\angle_{y_3}(x_3,x_i)=\arccos(-\third)$ and $d(x_3,x_i)=\2pithird$ for $i=1,2$.
Thus, $\angle_{x_i}(x_{i-1},x_{i+1})=\arccos(-\third))$ for $i=1,2$ 
(the indices to be understood modulo 3) and
these angles are of type 727. 
Let $y_1:=m(x_2,x_3)$ and $y_2:=m(x_1,x_3)$. Then it also follows that 
$d(x_i,y_i)=\2pithird$ for $i=1,2$.
Consider the vertices $x_1,x_3,x_2,y_2$, then we are again in the situation
of the configuration $\ast$ (just exchange the indices $2\leftrightarrow3$).
It follows as above that $\angle_{y_2}(x_2,x_3)=\ac(-\third)$ because $x_1$ is
an $8T$-vertex. This implies that $\angle_{x_3}(x_1,x_2)=\ac(-\third)$ as well,
and this angle is of type 727.

The convex hulls $CH(x_i,y_j,x_j)$ for distinct $i,j=1,2,3$ are isosceles
spherical triangles with sides of type 87878, 87878 and 878. 
This implies $d(y_i,y_{i+1})=\pihalf$ and the segments $y_iy_{i+1}$
are of type 828.
The intersection $CH(x_i,y_{i+1},x_{i+1})\cap CH(x_i,y_{i-1},x_{i-1})$ is the
spherical triangle $CH(x_i,y_{i-1},y_{i+1})$ with sides of type 878, 878 and 828.
In particular the 8-vertices $m(x_i,y_i)$ are pairwise distinct.

Observe that the 2-vertices $\ora{y_3y_2},\ora{y_3y_1}\in \Si_{y_3} K$ are
adjacent to the antipodal 7-vertices $\ora{y_3x_1},\ora{y_3x_2}$, respectively.
This implies that $d(\ora{y_3y_2},\ora{y_3y_1})\geq \ac(\third)>\pithird$, thus
$d(\ora{y_3y_2},\ora{y_3y_1})\geq\pihalf$.
On the other hand, triangle comparison for the triangle $(y_1,y_2,y_3)$ implies 
$d(\ora{y_3y_2},\ora{y_3y_1})\leq\pihalf$ and it follows that this triangle 
is rigid, i.e.\ the convex hull
$CH(y_1,y_2,y_3)$ is an equilateral spherical triangle with sides of type 828.
Let $z_i:=m(y_i,y_{i-1})$. Notice that $z_i$ does not lie on the 
segment $x_iy_i$ of type 87878.
Let $w$ be the 7-vertex at the center of the triangle $CH(y_1,y_2,y_3)$
and consider the spherical triangles $CH(x_i,z_i,y_i)$ for $i=1,2,3$ 
with sides of type 82, 2768 and 87878. 
Notice that $w$ is the 7-vertex on the segments $z_iy_i$. It follows that
$w$ is adjacent to the $8A$-vertices $m(x_i,y_i)$ for $i=1,2,3$ and in particular,
$\Si_w K$ contains three pairwise antipodal 8-vertices.
\begin{center}
\hpic{\includegraphics[scale=0.55]{propertyTb}}
\hspace{2cm}
\hpic{\includegraphics[scale=0.55]{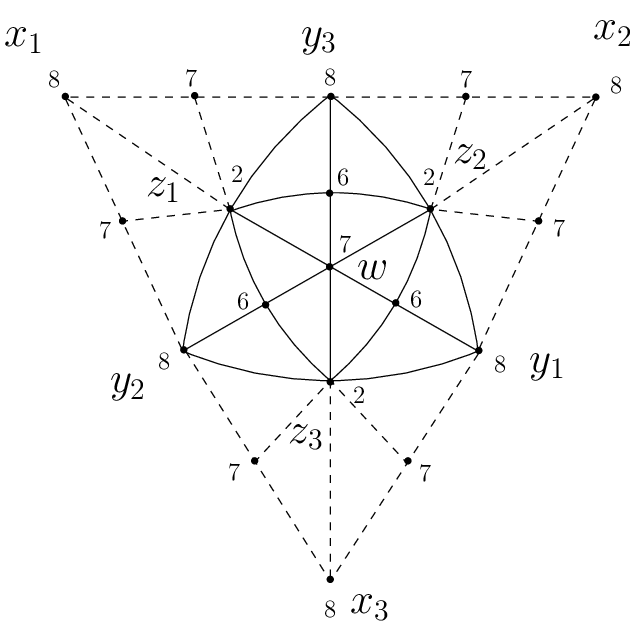}}
\end{center}

We say that an $8T$-vertex \textbf{$x\in K$ has the {\em property $T'$}}, if 
$rad(z_i,\{8\text{-vert. in } K\})\leq \arccos(-\frac{1}{2\sqrt{2}})$ 
for $i=1,2,3$ and for any such configuration of vertices $x_1,x_2,x_3\in G\cdot x$.

\begin{lem}\label{lem:8T'}
$K$ contains no $8T'$-vertices.
\end{lem}
\proof
Suppose there are $8T'$-vertices. We use the notation as in the definition
of the property $T'$. Let $w$ be the center of the triangle $CH(y_1,y_2,y_3)$.

Let $u\in K$ be an $8$-vertex. Then for some $i=1,2,3$, $\angle_w(z_i,u)\geq\pihalf$.
Suppose w.l.o.g.\ that it holds for $i=1$.
If $d(w,u)=\frac{5\pi}{6}$, then
$\ora{wu}$ is an 8-vertex and $\angle_w(z_1,u)=\pihalf$. It follows that
$d(u,z_1)=\3piquart$, but this contradicts the definition of
the property $T'$.
If $d(w,u)=\arccos(-\frac{1}{\sqrt{3}})$, then $\ora{wu}$ is a 2-vertex
and $\angle_w(z_1,u)=\2pithird$. It follows again that $d(u,z_1)=\3piquart$. 
Hence, $d(w,u)\leq\arccos(-\frac{1}{2\sqrt{3}})$ for all 8-vertices $u\in K$
and as in the beginning of the proof of Lemma~\ref{lem:6sph7pt} 
we deduce by triangle comparison that if $w'\in G\cdot w$, then
$d(w,w')\leq\arccos(-\third)$. We may also choose $w'$, so that $d(w,w')>\pihalf$.

Case 1: $d(w,w')=\arccos(-\third)$. If the segment $ww'$ is singular of type $76867$, then
for some $i=1,2,3$, $\angle_w(y_i,w')=\2pithird$ and this angle is of type 626. It follows
that $d(w',y_i)=\arccos(-\frac{1}{\sqrt{3}})$, a contradiction. If the simplicial convex
hull of $ww'$ is 2-dimensional, we can argue as in 
the proof of Lemma~\ref{lem:6sph7pt} (Case 1) 
to see that this case is not possible either.

Case 2:  $d(w,w')=\arccos(-\frac{1}{6})$. The argument in the proof
of Lemma~\ref{lem:6sph7pt} (Case 2)
rules out the case where $ww'$ has a 2-dimensional simplicial convex hull.

\parpic{\includegraphics[scale=0.6]{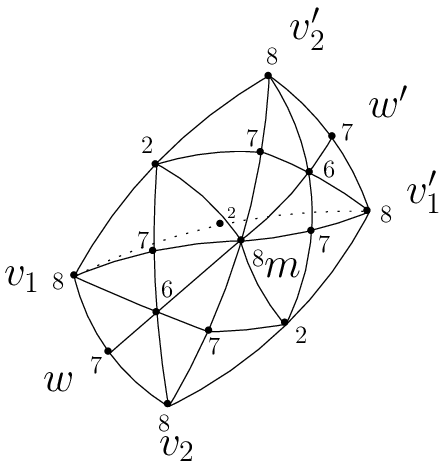}}
It remains to show that the case where the simplicial 
convex hull $C$ of $ww'$ is 3-dimensional
is not possible either. Let $v_1,v_1'\in C$ be the 8-vertices adjacent to $w$ and $w'$, 
respectively. Notice that they are $8A$-vertices, 
otherwise an antipode of e.g.\ $v_1$ in $K$ would
have distance $\frac{5\pi}{6}$ to $w$; but this cannot happen.
Recall that there is an 8-vertex $m\in C$, such that $mw$ and $mw'$ are segments of type
867 and $\angle_m(w,w')=\ac(-\frac{3}{4})$.
Let $v_2\in K$ be an $8A$-vertex adjacent to $w$ and so that $v_1wv_2$ is a segment of type 878.
We can choose $v_2$ to be one of the $8A$-vertices $m(x_i,y_i)$.
Define $v_2'$ analogously. 
Then the convex hulls $CH(m,v_1,v_2)$ and $CH(m,v_1',v_2')$ are equilateral spherical
triangles with sides of type 878.

\parpic[r]{\includegraphics[scale=0.4]{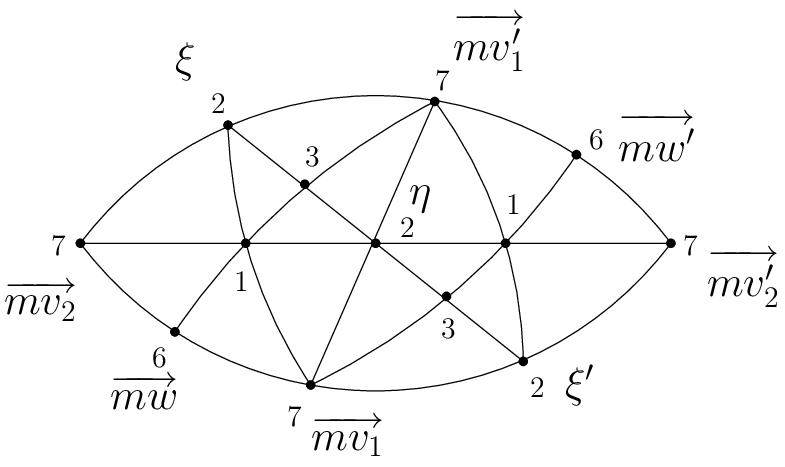}}
We want now to consider the convex hull $C':=CH(C,v_2,v_2')$.
The link $\Si_m C$ is a 2-dimensional spherical quadrilateral with vertices
$\ora{mw}$, $\ora{mv_1}$, $\ora{mw'}$ and $\ora{mv_1'}$.
Notice that $\ora{mv_2}\ora{mw}\ora{mv_1}$
and $\ora{mv_2'}\ora{mw'}\ora{mv_1'}$ are segments of type 767. 
It follows that
$CH(\Si_m C, \ora{mv_2},\ora{mv_2'})$ is a bigon connecting the antipodal 7-vertices
$\ora{mv_2}$ and $\ora{mv_2'}$. 
Then $d(v_2,v_2')=\2pithird$ and $m=m(v_2,v_2')$,
in particular, $m$ is an $8A$-vertex. 
Let $\xi,\xi'\in\Si_m C'$ be the 2-vertices
$m(\ora{mv_1'},\ora{mv_2})$ and $m(\ora{mv_1},\ora{mv_2'})$.
Let $\eta$ be the 2-vertex $m(\ora{mv_1},\ora{mv_1'})$.
The convex hulls $CH(v_1,v_2,v_2')$ and $CH(v_1',v_2',v_2)$
are spherical triangles with sides of type 878, 87878 and 828.

Since $m\in K$ is contained in the convex hull of the $8T'$-vertices, it is 
also contained in the convex hull of the $8T$-vertices. We can find another
$8T$-vertex $u_1\in K$, such that $d(m,u_1)=\2pithird$. Notice that the $8A$-vertex
$u_1$ cannot be antipodal to $v_2$ or $v_2'$, in particular,
$\angle_m(u_1,v_2),\angle_m(u_1,v_2')<\pi$.
Suppose w.l.o.g.\ that $\angle_m(u_1,v_2)\geq \pihalf$. 
Then $\angle_m(u_1,v_2)=\arccos(-\third)$ and $d(u_1,v_2)=\2pithird$.
$CH(v_2,m,u_1)$ is an isosceles spherical triangle 
(as in the configuration $\ast$) with a
2-vertex $z$ in its interior.
Recall that $d(w,u_1)\leq\arccos(-\frac{1}{2\sqrt{3}})$. This implies that
$\angle_{v_2}(w,u_1)\leq\ac(\third)$. This angle cannot be 0, because
$\angle_{v_2}(m,w)=\ac(\third)$ and $\angle_{v_2}(m,u_1)=\ac(-\third)$.
Thus $\angle_{v_2}(w,u_1)=\ac(\third)$ and it is of type 767.
$CH(\ora{v_2w},\ora{v_2m},\ora{v_2u_1})$ is then
a spherical triangle with sides of type 767, 767 and 727. In particular $w$ is 
adjacent to the 2-vertex $z$.

\parpic{\includegraphics[scale=0.55]{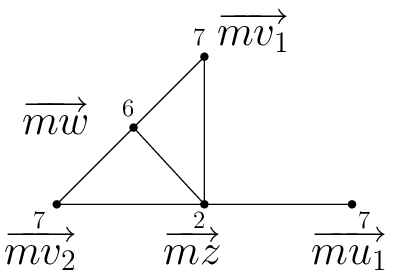}}
This consideration implies in the link $\Si_m K$ that $\ora{mz}$ and $\ora{mw}$ are
adjacent. Suppose that the segment $\ora{mv_1}\ora{mu_1}$ is of type 727.
This implies that the angle $\angle_{\ora{mv_1}}(\ora{mu_1},\xi')$ is of type 262.
It follows that the segment $\ora{mu_1}\xi'$ is of type 7672. Hence, 
$d(u_1,m(v_2',v_1))=\3piquart$ and $CH(v_1,v_2',u_1)$ is a spherical triangle with sides
87878, 87878 and 828. But this contradicts the definition of the property $T$ for $u_1$.
Therefore the segment $\ora{mv_1}\ora{mu_1}$ is of type 767.

If $\angle_m(u_1,v_2')=\arccos(-\third)$ we argue analogously and conclude that the
segment $\ora{mv_1'}\ora{mu_1}$ is of type 767. If $\angle_m(u_1,v_2')=\ac(\third)$ we
see as above that $d(\ora{mu_1},\xi)\leq\pihalf$, otherwise we violate the property $T$ for $u_1$.
Using triangle comparison with the triangle $(\xi, \ora{mu_1}, \ora{mv_2'})$ 
(or using the convexity of the ball centered at $\ora{mu_1}$ with radius $\pihalf$)
we see that $d(\ora{mv_1'},\ora{mu_1})\leq\ac(\third)$. 
Since $\ora{mv_1}\ora{mu_1}$ is of type 767, then $\ora{mu_1}\neq\ora{mv_1'}$.
Thus, $d(\ora{mv_1'},\ora{mu_1})=\ac(\third)$ and the 
segment $\ora{mv_1'}\ora{mu_1}$ is of type 767 also in this case. 
It follows that $CH(\ora{mv_1'},\ora{mv_1},\ora{mu_1})$ is a spherical triangle
with sides 767, 767 and 727. In particular $\ora{mu_1}$ is adjacent to $\eta$.

We have shown so far that any 7-vertex in $\Si_m K$ that is $\2pithird$-extendable to
an $8T$-vertex in $K$ must be adjacent to $\eta$ and the segments connecting it
with $\ora{mv_1}$ and $\ora{mv_1'}$ are of type 767.

Let $r_1:=m(m,u_1)\in K$ and let $u_2'\in K$ be an $8T$-vertex with
$d(r_1,u_2')=\2pithird$.
Since $u_1$ is an $8T$-vertex, the angle $\angle_{r_1}(m,u_2')$ cannot
be of type 767. Hence, it is of type 727. If the angle $\angle_{r_1}(u_1,u_2')$
is also of type 727, then set $u_2:=u_2'$.

\parpic{\includegraphics[scale=0.65]{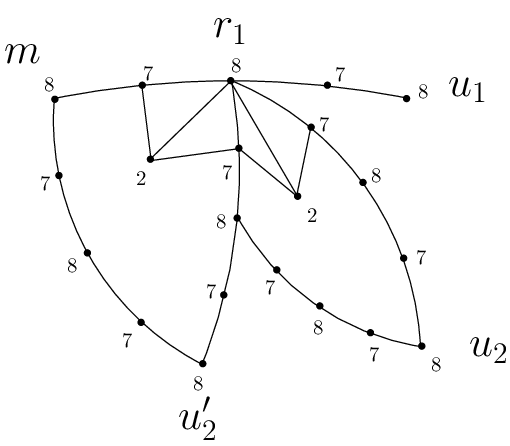}}
 Otherwise, let $u_2\in K$ be another $8T$-vertex, 
so that $d(u_2,m(r_1,u_2'))=\2pithird$.
Again, because $u_2'$ is an $8T$-vertex, the angle $\angle_{m(r_1,u_2')}(r_1,u_2)$
is of type 727.
In particular $d(r_1,u_2)=\2pithird$ and again $\angle_{r_1}(m,u_2)$ is of type 727.
We want to see now, that $\angle_{r_1}(u_1,u_2)$ is also of type 727.
Suppose that $\angle_{r_1}(u_1,u_2)$ is of type 767. Then 
$CH(\ora{r_1u_2},\ora{r_1u_1},\ora{r_1u_2'})$ is a spherical triangle with sides
of type 767, 767 and 727. In particular $\ora{r_1u_1}$ is adjacent to 
$\delta:=m(\ora{r_1u_2},\ora{r_1u_2'})$, 
this means that the segment $\delta\ora{r_1m}$ is of type 2767. 
Notice that this is 
the configuration $\aast$ for the vertices $r_1,u_2',u_2,m$.
This implies that $CH(r_1,u_2,m(m,u_2'))$
is a spherical triangle with vertices of type $8A$ and sides 
of type 87878, 87878 and 828 and $u_2$ could not
be an $8T$-vertex, a contradiction.

\parpic{\includegraphics[scale=0.51]{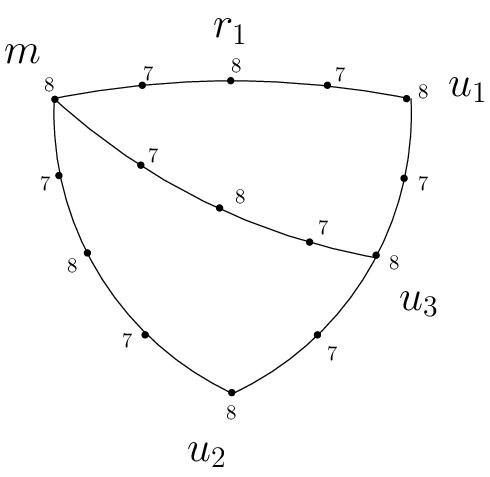}}
Thus $\angle_{r_1}(u_1,u_2)$ is of type 727.
This implies that $d(u_1,u_2)=\2pithird$
and $\angle_{u_1}(m,u_2)$ is of type 727. Let $u_3\in K$ be the $8A$-vertex
$m(u_1,u_2)$, then $\angle_{u_1}(m,u_3)$ is of type 727 and this implies
that $d(m,u_3)=\2pithird$. Observe that $u_3$ is not necessarily an $8T$-vertex.
Notice that $\ora{mu_1}\ora{mu_2}$ is of type 727 and recall that
$\ora{mu_i}$ is adjacent to $\eta$ for $i=1,2$. It follows
that $\eta=m(\ora{mu_1},\ora{mu_2})$. In particular $\eta$ is
$\pihalf$-extendable in $K$.
Consider the triangles $(m,u_1,u_3)$ and $(m,u_2,u_3)$, then
by triangle comparison, it follows that 
$\angle_m(u_1,u_3), \angle_m(u_2,u_3)\leq \ac(\third)$
and since $\angle_m(u_1,u_2)=\ac(-\third)$, this implies that,
$\angle_m(u_1,u_i)=\ac(\third)$ 
and $\ora{mu_3}\ora{mu_i}$ is of type 767 for $i=1,2$. 
Hence, $CH(\ora{mu_1},\ora{mu_2},\ora{mu_3})$ is a spherical triangle with sides
of type 767, 767 and 727. In particular,
$\ora{mu_3}$ is adjacent to $\eta$ as well.

Write $\ora{\eta\star}:=\ora{\eta\ora{m\star}}\in \Si_\eta\Si_m K$, 
where $\star$ is
any vertex in $K$ adjacent to $m$, so that 
$\ora{m\star}\in\Si_m K$ is adjacent to $\eta$.

\parpic{\includegraphics[scale=0.5]{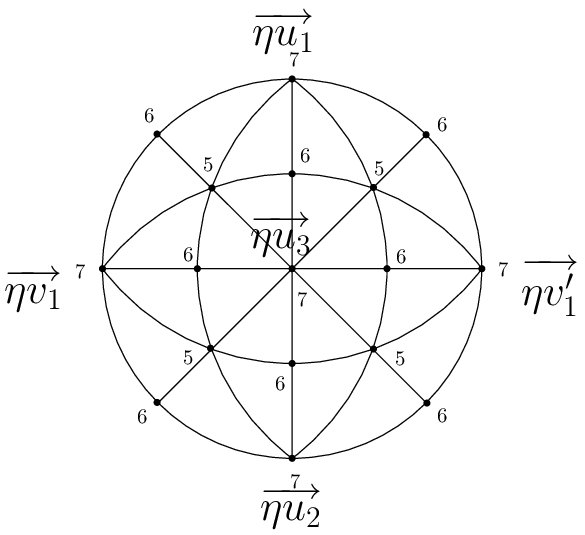}}
The 7-vertices $\ora{\eta v_1}$, $\ora{\eta v_1'}$, $\ora{\eta u_1}$ and $\ora{\eta u_2}$
are the 7-vertices of a circle $c\subset \Si_\eta\Si_m K$ of type 767676767,
because as seen above, $\ora{\eta u_i}$ for $i=1,2$ is the midpoint of a geodesic
of length $\pi$ connecting $\ora{\eta v_1}$ and $\ora{\eta v_1'}$,
and $\ora{\eta u_i}$ are antipodal for $i=1,2$.
From the construction above we see that $d(\ora{\eta u_3},\ora{\eta u_i})=\pihalf$
for $i=1,2$ (the segments $\ora{\eta u_3}\ora{\eta u_i}$ are of type 767).
Suppose $\ora{\eta u_3}$ is antipodal to $\ora{\eta v_1}$.
This would imply that the segment $\ora{mu_3}\ora{mw}\subset \Si_m K$
is of type 7316 and therefore $d(\ora{mu_3},\ora{mw})>\pihalf$
(compare with the figure for $\Si_m C'$ above).
Consider now the triangle $(w,m,u_3)$, it has sides
$d(m,w)=\ac(\frac{1}{\sqrt{3}})$, $d(m,u_3)=\2pithird$ and angle 
$\angle_m(w,u_3)>\pihalf$. It follows that 
$d(w,u_3)>\ac(-\frac{1}{2\sqrt{3}})$, which is not possible.
Hence, $d(\ora{\eta u_3},\ora{\eta v_1})=d(\ora{\eta u_3},\ora{\eta v_1'})=\pihalf$.
Therefore $\ora{\eta u_3}$ is the center of a 2-dimensional hemisphere 
in $\Si_\eta\Si_m K$ bounded by $c$.

Let $r_3:=(m,u_3)\in K$ and let $u_4'\in K$ be another $8T$-vertex, 
so that $d(r_3,u_4')=\2pithird$. 
Recall that $u_3$ is not necessarily an $8T$-vertex, therefore we cannot
conclude directly that $\angle_{r_3}(m,u_4')$ is of type 727.
If $\angle_{r_3}(m,u_4')$ is actually of type 727, then
set $u_4:=u_4'$. 

\parpic{\includegraphics[scale=0.65]{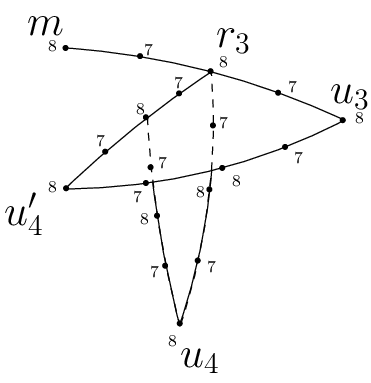}}
Otherwise (i.e.\ if $\angle_{r_3}(m,u_4')$ is of type 767),
 let $u_4\in K$ be an $8T$-vertex, 
so that $d(u_4,m(r_3,u_4'))=\2pithird$.
Then, since $u_4'$ is an $8T$-vertex, the angle $\angle_{m(r_3,u_4')}(r_3,u_4)$
must be of type 727. This implies that $d(r_3,u_4)=\2pithird$ and the angle 
$\angle_{r_3}(u_4,u_4')$ is of type 727.
It follows that $\angle_{r_3}(m,u_4)$ is of type 727, otherwise (as in
the argument above for $u_2$) we find the configuration $\aast$ and 
$CH(u_3,u_4,m(r_3,u_4'))$ is a spherical triangle with sides of type 87878, 87878
and 828, contradicting the property $T$ for $u_4$.
From this we conclude that $d(m,u_4)=\2pithird$ and $\angle_m(u_3,u_4)$ is of type 727.
Recall that $\ora{mu_4}$ must be adjacent to $\eta$. This implies that
$\ora{\eta u_4}$ is antipodal to $\ora{\eta u_3}$.

Thus, $\Si_\eta\Si_m K$ (of type $D_6$)
contains a singular 2-sphere spanned by 3 pairwise orthogonal
7-vertices. Recall that it also contains a pair of antipodal 3-vertices $\ora{\eta\xi}$
and $\ora{\eta\xi'}$. Lemma~\ref{lem:Dn_n-2sph} implies that
$\Si_\eta\Si_m K$ contains a 3-sphere spanned by a simplex of type $1567$.
Since $\eta$ is $\pihalf$-extendable in $K$, 
we have found a 2-vertex in $K$, whose link contains
a 4-sphere spanned by a simplex of type $15678$. 
This 4-sphere is of type $\pithird$ (this can be easily seen in the vector space
realization of the Coxeter complex of type $D_n$ presented in Appendix~\ref{app:coxeter}).
A contradiction to Lemma~\ref{lem:4sph2ptb}.
\qed

\parpic[r]{\includegraphics[scale=0.3]{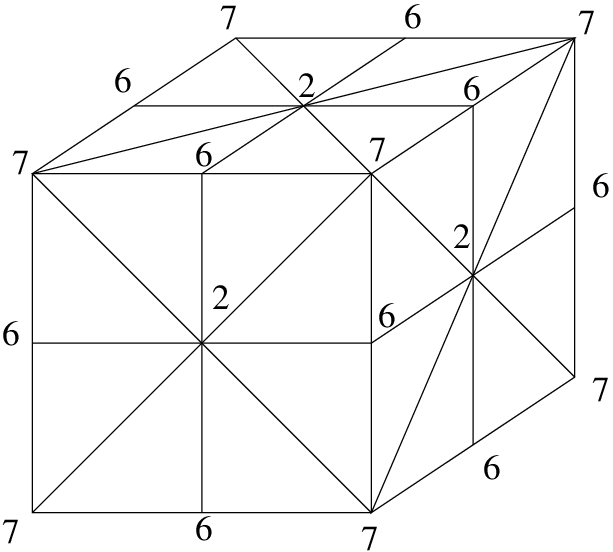}}
\picskip{3}
\vspace{1cm}
Let \textbf{$B_3$ be the property of an $8A$-vertex $x\in K$}, 
such that $\Si_x K$ contains a singular
2-sphere with $B_3$-geometry 
\includegraphics[scale=0.5]{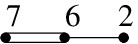}, and such that all the
$7$-vertices in this sphere are $\pithird$-extendable.
\bigskip

Consider the configuration $\aast$ and notice that the 8-vertex $v$
on the segment $zx_3$
(of type 2828) adjacent to $z$ is an $8B_3$-vertex. 
Indeed, the convex hull of the 8-vertices 
$x,x_1,w,y_3$, $m(x,x_3),m(x_1,x_3),m(w,x_3),m(y_3,x_3)$ is a spherical cube
with these 8-vertices as vertices. The 8-vertex $v$ is the center of this
spherical cube and therefore it has a 2-sphere in its link. This sphere
is spanned by three pairwise orthogonal pairs of antipodal 2-vertices 
(the directions of the segments from $v$ to the 
centers of the faces of the spherical cube).
Such a sphere has  $B_3$-geometry 
\includegraphics[scale=0.5]{B3_762dynk}.

Another similar way of finding $8B_3$-vertices is the following. 
Let $x_1,x_2,x_3,x_4\in K$ be $8A$-vertices adjacent to a 2-vertex $y$, so
that $CH(x_i)$ is a 2-dimensional spherical quadrilateral with sides $x_ix_{i+1}$ 
of type 878. Let $x\in K$ be an 8-vertex at distance $\3piquart$ to $y$.
Since the $x_i$ are $8A$-vertices, it follows that $\angle_y(x,x_i)=\pihalf$.
This implies that $\Si_{\ora{yx}}\Si_y K$ contains a singular circle of type
767676767. Let $z$ be the 8-vertex in $yx$ adjacent to $y$. Then
$\Si_z K$ contains a 2-sphere with $B_3$-geometry 
\includegraphics[scale=0.5]{B3_762dynk}.
Considering the spherical triangles $CH(x,x_i,x_{i+2})$, 
we see that the 7-vertices in this 2-sphere
are $\pithird$-extendable. Hence $z$ is an $8B_3$-vertex.
\begin{center}
 \hpic{\includegraphics[scale=0.6]{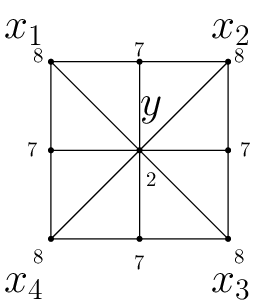}}\hspace{3cm}
 \hpic{\includegraphics[scale=0.45]{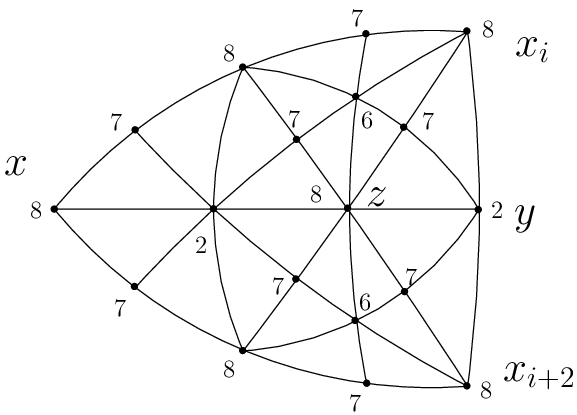}}
\end{center}

Consider now the definition of the property $T'$. The 2-vertices $z_i$ are centers
of 2-dimensional spherical quadrilaterals as described above.
In particular, if there are no $8B_3$-vertices in $K$, then it follows from
the observation above, that
$rad(z_i,\{8\text{-vert. in } K\})\leq \arccos(-\frac{1}{2\sqrt{2}})$
for $i=1,2,3$.
Hence, if $K$ contains no $8B_3$-vertices, it follows that Property $T$ implies
Property $T'$.

Recall that our strategy is to find spheres of large dimension in the links of
vertices of type 2 or 8. Notice that we have made the first step in this
direction:

\begin{cor}\label{cor:8Aptwithcircle}
 If $K$ contains $8A$-vertices, then it contains 8-vertices, whose links 
in $K$ contain a singular circle.
\end{cor}
\proof
If $K$ contains $8B_3$-vertices, we are done. Otherwise, 
$8T\Rightarrow 8T'$, and Lemma~\ref{lem:8T'} implies that there are no $8T$-vertices in $K$.
In particular, we find a spherical triangle in $K$ with sides of type 87878, 87878 and 828.
The link in $K$ of the 8-vertex in the interior of this triangle contains a singular circle. 
\qed

Now we find 8-vertices, such that their links contain singular 2-spheres.

\begin{lem}\label{lem:exist8B3pt}
 If $K$ contains $8A$-vertices, then it also contains $8B_3$-vertices.
\end{lem}
\proof
Suppose that $K$ contains $8A$-vertices but no $8B_3$-vertices. Then, $8T\Rightarrow 8T'$
and Lemma~\ref{lem:8T'} implies that there are no $8T$-vertices in $K$.

\parpic{\includegraphics[scale=0.5]{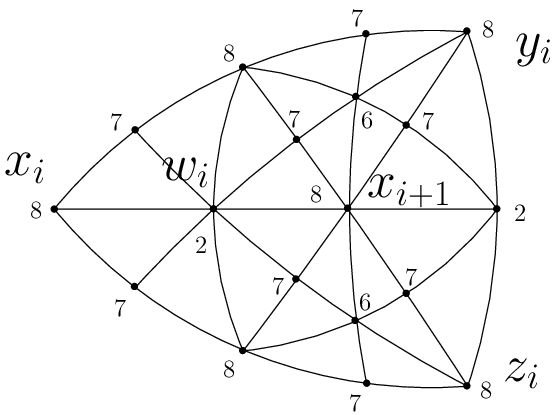}}
Hence, there are $8A$-vertices $x_0,y_0\in K$ and an 8-vertex $z_0\in K$, so that
$T_0:=CH(x_0,y_0,z_0)$ is a spherical triangle with sides of type 87878, 87878 and 828; where
$y_0z_0$ is the side of type 828 (as in the definition of the property T).
Let $x_1\in K$ be the $8A$-vertex on the segment $x_0m(y_0,z_0)$ (of type 8282)
adjacent to the 2-vertex $m(y_0,z_0)$. Since $x_1$ is not an $8T$-vertex, we can
find 8-vertices $y_1,z_1\in K$ as vertices of a spherical triangle
$T_1:=CH(x_1,y_1,z_1)$ as above. Define $x_i,y_i,z_i\in K$ and $T_i\subset K$ inductively.
Let $w_i$ be the $2A$-vertex $m(x_i,x_{i+1})$.

If $\xi\in \Si_{x_i} K$ is a $\pithird$-extendable 7-vertex and 
$d(\xi,\ora{x_ix_{i+1}})=\ac(-\frac{1}{\sqrt{3}})$, then we are in the
setting of the configuration $\aast$
because $\ora{x_iy_i}$ and $\ora{x_iz_i}$ are both 
$\2pithird$-extendable to $8A$-vertices (definition of the property $T$).
This implies that there are $8B_3$-vertices in $K$, contradicting our assumption.
Hence, $\ora{x_ix_{i+1}}$ has distance $\leq \pihalf$ to all $\pithird$-extendable
7-vertices in $\Si_{x_i} K$. Notice also that $d(\ora{x_ix_{i-1}},\ora{x_iy_i})$
and $d(\ora{x_ix_{i-1}},\ora{x_iz_i})$ are both $\leq\pihalf$, otherwise $w_{i-1}$
would have distance $\3piquart$ to the 8-vertex $y_i$ or $z_i$ and we would find
an $8B_3$-vertex on the segment $w_{i-1}y_i$ ($w_{i-1}z_i$).

\parpic{\includegraphics[scale=0.53]{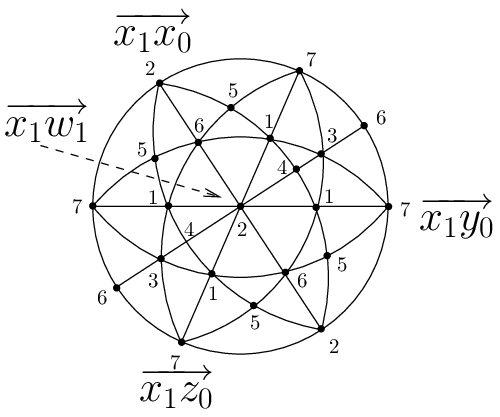}}\vspace{0.5cm}
From these observations it follows, that $\ora{x_1w_1}$ has distance $\equiv\pihalf$
to the circle $\Si_{x_1}T_0$ of type 727672767.
This implies that $\Si_{\ora{x_1w_1}}\Si_{x_1}K$ 
(of type \hpic{\includegraphics[scale=0.4]{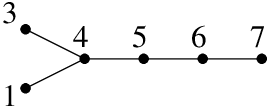}})
contains a singular circle of type 161416141.
It also contains the pair of antipodal 7-vertices 
$\xi:=\ora{\scriptstyle{\ora{x_1w_1}\ora{x_1y_1}}}$ and 
$\xi':=\ora{\scriptstyle{\ora{x_1w_1}\ora{x_1z_1}}}$.
\bigskip

Since $d(\ora{x_1x_0},\ora{x_1y_1})$, $d(\ora{x_1x_0},\ora{x_1z_1})\leq\pihalf$
and $d(\ora{x_1x_0},\ora{x_1w_1})=\pihalf$, it follows from triangle comparison that
$d(\ora{x_1x_0},\ora{x_1y_1})=d(\ora{x_1x_0},\ora{x_1z_1})=\pihalf$, because
the triangle $(\ora{x_1x_0},\ora{x_1y_1},\ora{x_1z_1})$ must be rigid.
Let $\zeta:=\ora{\scriptstyle{\ora{x_1w_1}\ora{x_1x_0}}}$. Then the segments
$\zeta\xi$ and $\zeta\xi'$ have length $\pihalf$ and are of type 657.
\begin{center}
\hpic{\includegraphics[scale=0.5]{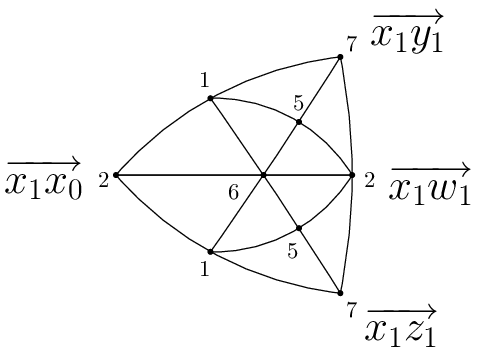}}\hspace{3cm}
\hpic{\includegraphics[scale=0.5]{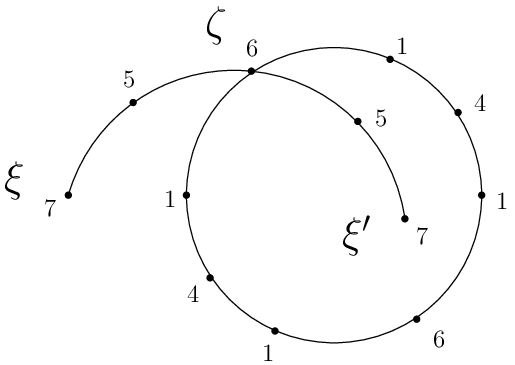}}
\end{center}

\begin{sublem}\label{sublem:exist8B4pt}
 $\Si_{\ora{x_1w_1}}\Si_{x_1}K$ contains a singular circle of type
$756575657$. This circle contains the vertices $\xi$, $\xi'$ and $\zeta$.
\end{sublem}
\proof
Let $\zeta'\in \Si_{\ora{x_1w_1}}\Si_{x_1}K$ be the 6-vertex in the circle of type
161416141 antipodal to $\zeta$.
If $d(\xi,\zeta')=\pihalf$, then $\zeta\xi\zeta'$ is a geodesic of type 65756. 
In particular, $\ora{\xi\zeta}$ has an antipode in 
$\Si_\xi\Si_{\ora{x_1w_1}}\Si_{x_1}K$ and we find the desired circle.
If $d(\xi,\zeta')>\pihalf$, then the segment $\xi\zeta'$
is of type 7676.

\parpic{\includegraphics[scale=0.6]{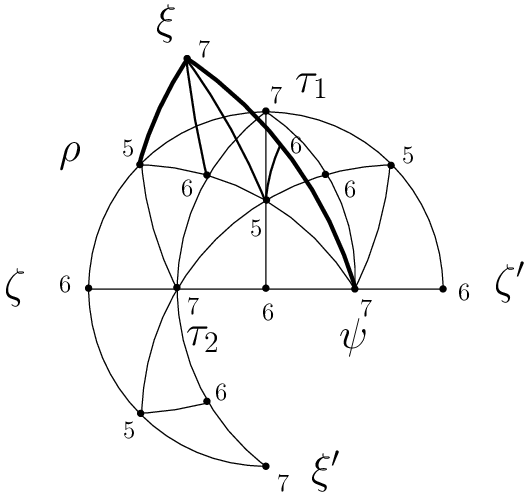}}
Let $\rho$ be the 5-vertex on the segment $\zeta\xi$ and let $\psi$ be the 7-vertex
on the segment $\xi\zeta'$ adjacent to $\zeta'$. 
Consider the geodesics $c_\rho$ and $c_\psi$ of length $\pi$
connecting $\zeta$ and $\zeta'$ through $\rho$ and $\psi$.
Let $\tau_1$ be the 7-vertex at the center of $c_\rho$ and $\tau_2$ be the 7-vertex
in $c_\psi$ adjacent to $\zeta$. 
Then $\rho$ and $\tau_2$ are adjacent
because $\Si_\zeta\Si_{\ora{x_1w_1}}\Si_{x_1}K$ is of type 
\hpic{\includegraphics[scale=0.4]{E8link268dynk}}.
$\xi$ cannot be adjacent to the 6-vertex at the center of $c_\psi$, otherwise
it would have distance $\3piquart$ to $\zeta$. Thus, the intersection of the segments
$\xi\zeta'$ and $c_\psi$ is the segment $\psi\zeta'$.
Considering the spherical triangle $CH(\rho, \xi, \psi)$ with
sides of type 57, 767 and 7565,
it follows that $\xi$ is adjacent to the 6-vertex $m(\tau_1,\tau_2)$
on the segment $\rho\psi$. 
In particular, $\xi'$ must be antipodal to 
at least one of $\tau_1$ or $\tau_2$. Since $\tau_2$ is adjacent
to $\zeta$ and $d(\zeta,\xi')=\pihalf$, then $\xi'$ cannot be antipodal to $\tau_2$.
It follows that $\xi'$ and $\tau_1$ are antipodal.
Let finally $c$ be the geodesic connecting $\tau_1$ and $\xi'$, so that the initial
direction coincides with $\ora{\tau_1\zeta'}$.
Then the initial direction of $c$ at $\xi'$ is antipodal 
to $\ora{\xi'\zeta}$ and we can find
the desired circle.
\qed

\bigskip\no
\emph{Continuation of proof of Lemma~\ref{lem:exist8B3pt}.}
\parpic{\includegraphics[scale=0.5]{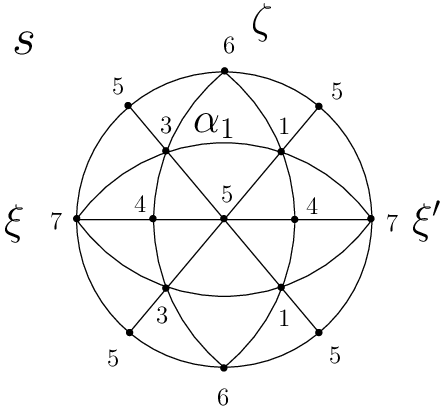}}
The link $\Si_\zeta\Si_{\ora{x_1w_1}}\Si_{x_1}K$ (of type 
\hpic{\includegraphics[scale=0.5]{E8link268dynk}})
 contains a pair of antipodal
5-vertices $\ora{\zeta\xi}$ and $\ora{\zeta\xi'}$ and a pair of antipodal 1-vertices.
We apply Corollary~\ref{lem:D4_circle} 
to conclude that 
$\Si_\zeta\Si_{\ora{x_1w_1}}\Si_{x_1}K$ contains a singular circle of type 5135135 with
$\ora{\zeta\xi}$ and $\ora{\zeta\xi'}$ on it. It follows now from
Sublemma~\ref{sublem:exist8B4pt} that $\Si_{\ora{x_1w_1}}\Si_{x_1}K$
contains a singular 2-sphere $s$ containing the vertices $\zeta$, $\xi$ and $\xi'$.
Therefore $\Si_{x_2} K$ contains a singular 3-sphere $S$ containing the singular
circle $\Si_{x_2}T_1$. 

We investigate below which 7-vertices in $S$ are $\pithird$-extendable. 
Clearly the 7-vertices in $\Si_{x_2}T_1\subset S$ are $\pithird$-extendable.

Let $\alpha_1,\alpha_2\in s$ be the 3-vertices adjacent to $\zeta$ and recall 
that $\zeta$ is $\pihalf$-extendable (to $\ora{x_1x_0}$) in $\Si_{x_1} K$. This implies
that $\alpha_i$ is $\pithird$-extendable to a segment of type 232 in $\Si_{x_1} K$.
Therefore, we find 7-vertices $\beta_1,\beta_2\in S$ 
at distance $\pihalf$ to $\ora{x_2w_1}$ which are $\pithird$-extendable in $K$
(compare with the figure below).
\begin{center}
 \hpic{\includegraphics[scale=0.5]{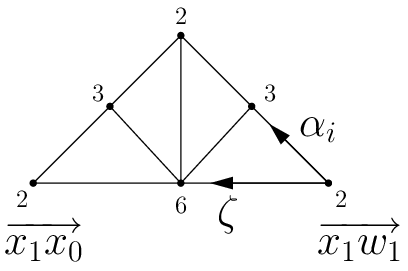}}\hspace{2cm}
 \hpic{\includegraphics[scale=0.6]{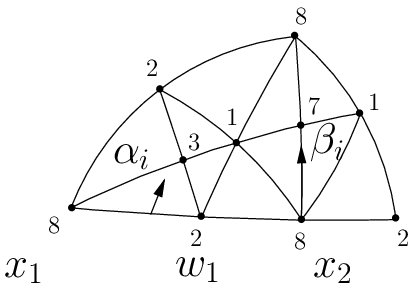}}
\end{center}

The segment $\alpha_1\alpha_2\subset \Si_{x_1w_1}\Si_{x_1}K$ is of type 363
with midpoint the 6-vertex $\zeta$, this
implies that the angle $\angle_{\ora{x_2w_1}}(\beta_1,\beta_2)$ is of type 161
and this implies in turn, that the 
segment $\beta_1\beta_2\subset\Si_{x_2}K$ is of type 727.
Let $\gamma\in S$  be the 2-vertex $m(\beta_1,\beta_2)$.
\begin{center}
 \hpic{\includegraphics[scale=0.7]{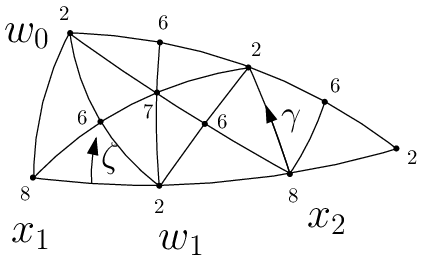}}\hspace{2cm}
 \hpic{\includegraphics[scale=0.6]{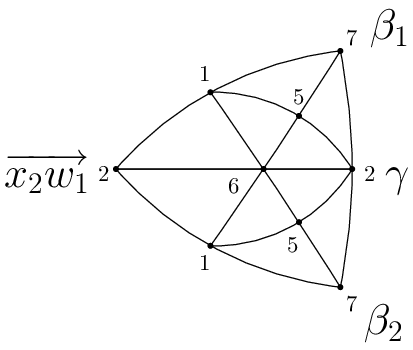}}
\end{center}

Let $\zeta_2:=\ora{\scriptstyle{\ora{x_2w_2}\ora{x_2x_1}}}$.
We can use the same argument as above to see that 
$\Si_{\zeta_2}\Si_{\ora{x_2w_2}}\Si_{x_2}K$ contains a singular 
circle of type 5135135. We want to prove next
that it also contains a pair of antipodal 7-vertices. 

\begin{sublem}\label{sublem:exist8B3ptb}
The link $\Si_{\zeta_2}\Si_{\ora{x_2w_2}}\Si_{x_2}K$ contains
a pair of antipodal 7-vertices.
\end{sublem}
\proof
Notice again that $d(\ora{x_2w_2},\Si_{x_2}T_1)\equiv\pihalf$, in particular,
$d(\ora{x_2w_2},\ora{x_2w_1})=\pihalf$.

\parpic{\includegraphics[scale=0.75]{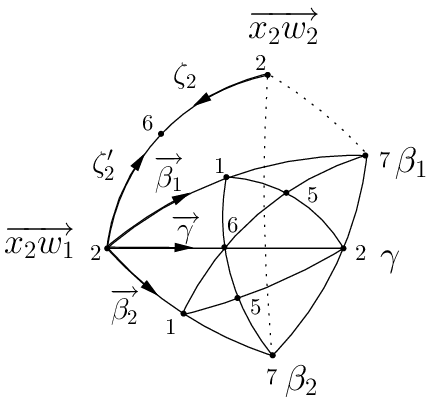}}
Recall also from the beginning of the proof Lemma~\ref{lem:exist8B3pt} that 
$d(\ora{x_2w_2},\beta_i)\leq\pihalf$, because $\beta_i$ is $\pithird$-extendable.
This implies that $d(\ora{x_2w_2},\gamma)\leq\pihalf$ and
$\angle_{\ora{x_2w_1}}(\ora{x_2w_2},\gamma)\leq\pihalf$.
Let $\zeta_2':=\ora{\scriptstyle{\ora{x_2w_1}\ora{x_2w_2}}}$,
$\ora{\beta_i}:=\ora{\ora{x_2w_1}\beta_i}$ and 
$\ora\gamma:=\ora{\ora{x_2w_1}\gamma}$.
We have already seen that $d(\zeta_2',\ora{\beta_i})\leq\pihalf$ and
$d(\zeta_2',\ora{\gamma})\leq\pihalf$. Furthermore, it follows from triangle comparison,
that if  $d(\zeta_2',\ora{\gamma})=\pihalf$, then 
$d(\zeta_2',\ora{\beta_i})=\pihalf$ for $i=1,2$. 

Notice that the link $\Si_{\zeta_2}\Si_{\ora{x_2w_2}}\Si_{x_2}K$ 
contains a pair of antipodal 7-vertices
if and only if $\Si_{\zeta_2'}\Si_{\ora{x_2w_1}}\Si_{x_2}K$
contains a pair of antipodal 7-vertices. The latter is what we will show.

\parpic{\includegraphics[scale=0.55]{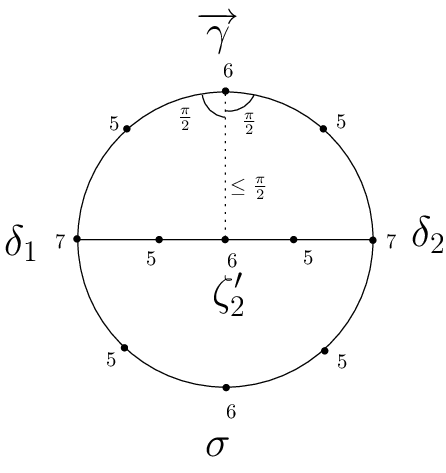}}
Let $\delta_1,\delta_2\in \Si_{\ora{x_2w_1}}\Si_{x_2}K$ be the two 7-vertices
in $\Si_{\ora{x_2w_1}}\Si_{x_2}T_1$
and recall that the 2-sphere
$\Si_{\ora{x_2w_1}} S$ contains a singular circle of type 756575657 containing the 
vertices $\delta_1$, $\delta_2$ and $\ora{\gamma}$
(this is just the circle in $\Si_{\ora{x_2w_1}} \Si_{x_2}K$
corresponding to the circle in $\Si_{\ora{x_1w_1}}\Si_{x_1} K$ from the
Sublemma~\ref{sublem:exist8B4pt} containing
$\xi$, $\xi'$ and $\zeta$).
Let $\sigma$ be the 6-vertex in this circle antipodal to $\ora{\gamma}$.
Further, we know that $d(\zeta_2',\delta_i)=\pihalf$, because 
$d(\ora{x_2w_2},\Si_{x_2}T_1)\equiv\pihalf$.
If $\zeta_2'$ has an antipode in the 2-sphere
$\Si_{\ora{x_2w_1}} S$, then $\ora{x_2w_2}$ has an antipode in $S$.
But this is impossible, since $\ora{x_2w_2}=\ora{x_2m(y_2,z_2)}$ and 
$m(y_2,z_2)\in K$ is a $2A$-vertex at distance $\3piquart$ to $x_2$. 
Hence $\pihalf\geq d(\zeta_2',\ora\gamma)> 0$ and $d(\zeta_2',\sigma)<\pi$.

Notice that $\Si_{\ora{x_2w_1}} \Si_{x_2} B$ is a building of type $D_6$ and
Dynkin diagram \hpic{\includegraphics[scale=0.4]{E7link2dynk}}. The distances
between 6-vertices are $0$, $\pithird$, $\pihalf$, $\2pithird$ and $\pi$.
The link $\Si_{\zeta_2'}\Si_{\ora{x_2w_1}}\Si_{x_2}K$ is of type 
\hpic{\includegraphics[scale=0.4]{E8link268dynk}}, thus two distinct 7-vertices
in this link must be antipodal.

\parpic[r]{\includegraphics[scale=0.5]{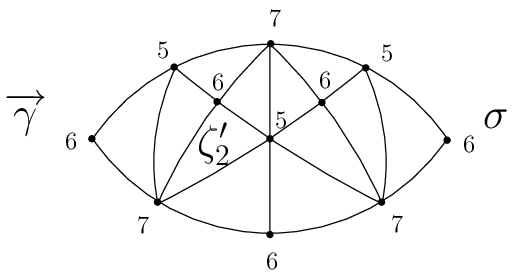}}
Case 1: $d(\zeta_2',\ora\gamma)=\pithird$. Since $d(\zeta_2',\sigma)<\pi$,
it follows that
$\ora\gamma\zeta_2'\sigma$ is a geodesic of length $\pi$. Its simplicial convex
hull is 2-dimensional and contains two 7-vertices adjacent to $\zeta_2'$.
It follows that $\Si_{\zeta_2'}\Si_{\ora{x_2w_1}}\Si_{x_2}K$
contains a pair of antipodal 7-vertices.

\parpic[r]{\includegraphics[scale=0.5]{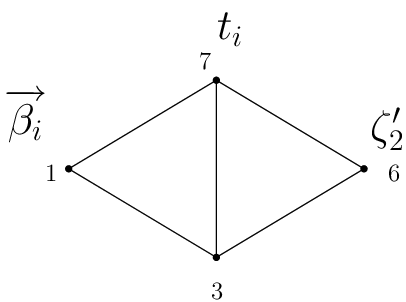}}
Case 2: $d(\zeta_2',\ora\gamma)=\pihalf$ and the segment $\zeta_2'\ora\gamma$
is of type 646. In this case, we know that
$d(\zeta_2',\ora{\beta_i})=\pihalf$ for $i=1,2$. Thus, 
$CH(\ora{\beta_1},\ora{\beta_2},\zeta_2')$ is an isosceles
spherical triangle with side lengths $\pihalf$, $\pihalf$ and $\ac(-\third)$. 
The simplicial convex hull of the segment
$\zeta_2'\ora{\beta_i}$ contains a 7-vertex $t_i$ adjacent to $\zeta_2'$ and
to $\ora{\beta_i}$ for $i=1,2$.
If $t_1=t_2$, then $t_1$ is adjacent to $\ora{\beta_i}$ for $i=1,2$. It follows 
that $t_1$ is also adjacent to $\ora{\gamma}=m(\ora{\beta_1},\ora{\beta_2})$. 
This means that 
$d(\ora\gamma,t_1)=d(t_1,\zeta_2')=\piquart$. 
Since $d(\zeta_2',\ora\gamma)=\pihalf$, $\zeta_2't_1\ora\gamma$ must be a geodesic.
This contradicts the fact that the segment $\zeta_2'\ora\gamma$ is of type 646.
Hence, $t_1\neq t_2$ and $\Si_{\zeta_2'}\Si_{\ora{x_2w_1}}\Si_{x_2}K$
contains a pair of antipodal 7-vertices.

Case 3: $d(\zeta_2',\ora\gamma)=\pihalf$ and the segment $\zeta_2'\ora\gamma$
is of type 676. If $d(\zeta_2',\sigma)=\pihalf$ then $\ora\gamma\zeta_2'\sigma$
is a geodesic of length $\pi$ and of type 67676.
If $d(\zeta_2',\sigma)=\2pithird$, then the segment $\zeta_2'\ora\gamma$
contains a 7-vertex adjacent to $\zeta_2'$ 
at distance $\piquart$ to $\ora\gamma$ and the simplicial
convex hull of the segment $\zeta_2'\sigma$ contains a 7-vertex 
adjacent to $\zeta_2'$ at distance
$\pihalf$ to $\sigma$. It follows that $\zeta_2'$ is adjacent to two different
7-vertices. Thus, $\Si_{\zeta_2'}\Si_{\ora{x_2w_1}}\Si_{x_2}K$
contains a pair of antipodal 7-vertices.
\qed 

\bigskip\no
\emph{End of proof of Lemma~\ref{lem:exist8B3pt}.}
We know now that
$\Si_{\zeta_2}\Si_{\ora{x_2w_2}}\Si_{x_2}K$
(of type \hpic{\includegraphics[scale=0.4]{E8link268dynk}})
contains a singular circle of type 5135135 and a pair of antipodal
7-vertices. Hence, it
contains a singular 2-sphere (the spherical join of the singular circle
and the pair of antipodal 7-vertices). 
Since $\zeta_2$ has an antipode in $\Si_{\ora{x_2w_2}}\Si_{x_2}K$, 
this implies that
$\Si_{\ora{x_2w_2}}\Si_{x_2}K$ contains a 3-sphere spanned by a
simplex of type 1567. This in turn implies that $\Si_{w_2}K$
contains a singular 4-sphere spanned by a simplex of type 15678.
This sphere is of type $\pithird$ as can be verified by considering
the vector space realization of the Coxeter complex of type $D_n$ 
presented in Appendix~\ref{app:coxeter}.
We get a contradiction to 
Lemma~\ref{lem:4sph2ptb} finishing the proof of the lemma.
\qed

\begin{lem}\label{lem:no8B3pt}
 $K$ contains no $8B_3$-vertices.
\end{lem}
\proof
We want to show first that an $8B_3$-vertex has the property $T$.
Suppose $x_1\in K$ is an $8B_3$-vertex and let $x_2,x_3\in K$ 
be $8$-vertices as in the 
configuration $\ast$.
Suppose further, that $x_3$ is an $8A$-vertex and that $\ora{x_1x_2}$ is 
$\2pithird$-extendable to an $8A$-vertex.
To prove that $x_1$ has the property $T$,
we have to show that $CH(x_1,x_2,x_3)$ is not a spherical triangle. 
Let $S\subset\Si_{x_1}K$ be the singular 2-sphere from
the definition of the property $B_3$. Let $\zeta:=\ora{x_1z_1}$ and 
$\xi_i:=\ora{x_1x_i}$ for $i=2,3$, as in the notation of the 
configuration $\ast$.

Suppose there is a 7-vertex $\xi \in S$, such that 
$d(\zeta,\xi)=\arccos(-\frac{1}{\sqrt{3}})$.
The segment $\zeta\xi$ is of type 2767. 
Since $\xi$ is $\pithird$-extendable in $K$ and $\xi_i$ is $\2pithird$-extendable to 
an $8A$-vertex, $\xi$ is not antipodal to $\xi_i$ for $i=1,2$. 
It follows that 
$CH(\xi,\xi_2,\xi_3)$ is an equilateral spherical triangle sides of type 727.
Let $\gamma$ be the 7-vertex in $\zeta\xi$
adjacent to $\zeta$. The 7-vertex $\gamma$ is the center of the spherical triangle
 $CH(\xi,\xi_2,\xi_3)$.
It follows from the configuration $\aast$, that 
$\gamma$ is $\2pithird$-extendable to an $8A$-vertex in $K$.

$\Si_\xi S$ is a singular circle of type 2626262. 
Notice that $\ora{\xi\zeta}=\ora{\xi\gamma}$ is not antipodal to any 
2-vertex in this circle, otherwise we could
find in $S$ an antipodal 7-vertex to $\gamma$, but this is not possible, since
$\gamma$ is $\2pithird$-extendable to an $8A$-vertex in $K$.
On the other hand, $\ora{\xi\zeta}$ cannot have distance $<\pihalf$
to all the 6-vertices in this circle, so
let $\eta$ be a 6-vertex in
$\Si_\xi S$, so that $d(\eta,\ora{\xi\zeta})=\2pithird$ and let
$\delta_i\in \Si_\xi S$ be the 2-vertices adjacent to $\eta$.
Let $\mu:=m(\eta,\ora{\xi\zeta})$.
(Compare with the configuration in the proof of Lemma~\ref{lem:3sph8pt}.)

\parpic{\includegraphics[scale=0.5]{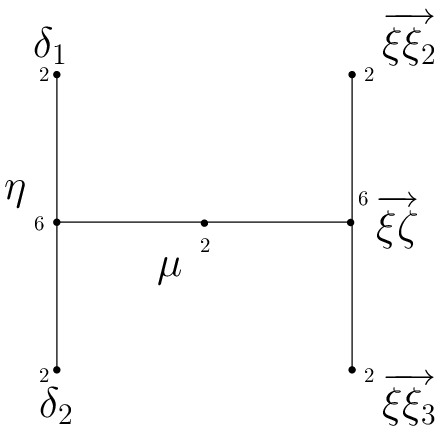}}
Since $\ora{\xi\zeta}$ is not antipodal to $\delta_i$,
it follows that $\angle_\eta(\delta_i,\ora{\xi\zeta})=\pihalf$ and
these angles are of type 232.
Therefore,
$\Si_{\ora{\mu\eta}}\Si_\mu\Si_\xi\Si_{x_1}K$ contains a pair of antipodal 3-vertices.
Similarly, we see that $\eta$ cannot be antipodal to $\ora{\xi\xi_i}$
because $\xi_i$ has no antipodes in $S$. Thus,
$\Si_{\ora{\mu{\scriptscriptstyle \ora{\xi\zeta}}}}\Si_\mu\Si_\xi\Si_{x_1}K$ 
contains a pair of antipodal 3-vertices.
This implies in turn, that $\Si_{\ora{\mu\eta}}\Si_\mu\Si_\xi\Si_{x_1}K$
contains a pair of antipodal 1-vertices.
We apply now Corollary~\ref{lem:D4_circle} to the building 
$\Si_{\ora{\mu\eta}}\Si_\mu\Si_\xi\Si_{x_1} B$ of type $D_4$
and conclude that $\Si_{\ora{\mu\eta}}\Si_\mu\Si_\xi\Si_{x_1}K$ contains a singular circle
of type 1351351. Therefore 
$\Si_\mu\Si_\xi\Si_{x_1}K$ contains a singular 2-sphere spanned by a simplex of type 156.
The same argument as in the proof of Lemma~\ref{lem:3sph8pt} 
(p. \pageref{lem:3sph8ptextend})
 shows that $\mu$ is extendable
in $\Si_{x_1}K$ to a segment of type 727
and the 2-vertex on this segment is extendable
in $K$ to a segment of type 828
(this uses that $\gamma\in \Si_{x_1}K$ is $\2pithird$-extendable 
and the 7-vertices in $S$ are $\pithird$-extendable). 
This produces a 2-vertex in $K$, whose link contains a 
4-sphere spanned by a simplex of type $15678$. This singular 4-sphere is of type $\pithird$,
a contradiction to Lemma~\ref{lem:4sph2ptb}.

From this, it follows that $\zeta$ has distance $\leq\pihalf$ to all the 7-vertices in $S$.
Since $S$ is the convex hull of its 7-vertices, 
it follows that $d(\zeta,S)\equiv \pihalf$. Hence $\Si_\zeta\Si_{x_1}K$
contains the 2-sphere $s:=\Si_\zeta CH(\zeta,S)$. The segments connecting $\zeta$
with the 2-vertices of $S$ are of type 262, the segments connecting $\zeta$
with the 7-vertices of $S$ are of type 217 and since the 6-vertices in $S$ are midpoints
of segments of type 767 in $S$, this implies that the segments connecting $\zeta$
with the 6-vertices of $S$ are of type 2436. Since the sphere $S$ has $B_3$-geometry
\includegraphics[scale=0.4]{B3_762dynk}, it follows that $s$
has $B_3$-geometry \includegraphics[scale=0.4]{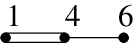}.
$\Si_\zeta\Si_{x_1}K$ also contains the two antipodal 7-vertices 
$\ora{\zeta\ora{x_1x_i}}$ for $i=2,3$.

\begin{sublem}\label{sublem:no8B3pt}
 Let $L\subset B$ be a convex subcomplex of a building of type $D_6$ with Dynkin
diagram \hpic{\includegraphics[scale=0.4]{E7link2dynk}}. Suppose $L$ contains
a singular $2$-sphere $S$ with $B_3$-geometry \includegraphics[scale=0.4]{B3_146dynk}
and also a pair of antipodal $7$-vertices. Then $L$ contains a $3$-sphere
spanned by a simplex of type $1467$.
\end{sublem}
\proof
Let $a,a'\in L$ be the antipodal 7-vertices and let $b,b'$ be antipodal 
1-vertices in $S\subset L$. By Lemma~\ref{lem:Dn_n-2sph} 
it follows that
$L$ contains a circle of type 7317317 through $b$ and $b'$. In particular
$\Si_b L$ contains a pair of antipodal 7- and 3-vertices.
$\Si_b S$ is a singular circle of type 6464646.
$\Si_b B$ is a building of type $A_5$.
So, it will suffice to show that under these circumstances $\Si_b L$ contains a 2-sphere
spanned by a simplex of type 467 (notice that such a sphere is also spanned by
a simplex of type 346):
\begin{center}
\hpic{\includegraphics[scale=0.5]{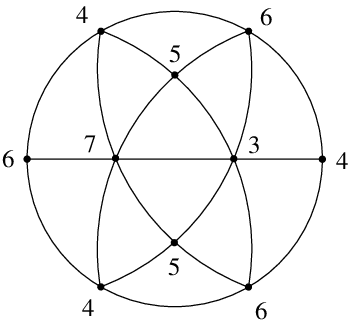}}
\end{center}

Let $d,d'\in \Si_b L$ be the antipodal 3- and 7-vertices, respectively.
Let $c\in \Si_b S$ be a 4-vertex and 
let $c'$ the 6-vertex in $\Si_b S$ antipodal to $c$.
If $c$ is adjacent to $d$, then $dcd'$ is a geodesic of type 3437 and
$\Si_c\Si_b L$ contains a pair of antipodal 3-vertices.
If $c$ is adjacent to $d'$, then $dcd'$ is a geodesic of type 3547 and
$\Si_c\Si_b L$ contains a pair of antipodal 5- and 7-vertices.

\parpic{\includegraphics[scale=0.55]{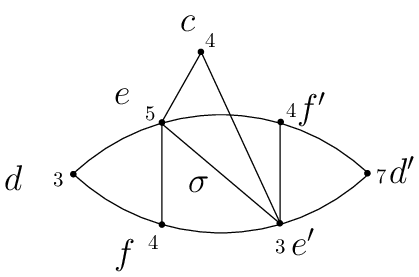}}
Otherwise the segments $cd$ and $cd'$ are of type 453 and 437 respectively.
Let $e$ be the 5-vertex in $cd$ and let $e'$ be the 3-vertex in $cd'$.
Let also $f$ be the 4-vertex on the segment $e'd$ (of type 343)
 and let $f'$ be the 4-vertex on the segment $ed'$ (of type 547).
Notice that since $e,e'$ are adjacent to $c$, then $e$ is adjacent to $e'$.
It follows that $e$ is adjacent to $f$ and $e'$ is adjacent to $f'$.
Let $\sigma$ be the edge $ee'$. The link $\Si_\sigma\Si_b B$ is of type 
\includegraphics[scale=0.4]{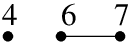}; and the direction 
$\ora{\sigma c'}$ is of type 4. It follows that $c'$ is antipodal to $f$ or $f'$
and $c'$ is contained in a circle in $\Si_b L$ of type 7673437 or 3657453.
This implies that $\Si_{c'}\Si_b L$ contains a pair of antipodal 7-vertices
or a pair of antipodal 3- and 5-vertices.
This means for $\Si_c\Si_b L$, that it contains a pair of antipodal 3-vertices
or a pair of antipodal 5- and 7-vertices.

Recall that $\Si_c\Si_b S$ consists of a pair of antipodal 6-vertices.
$\Si_c\Si_b B$ is of type $A_1\circ A_3$ with Dynkin diagram 
\includegraphics[scale=0.4]{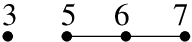}.
If $\Si_c\Si_b L$ contains a pair of antipodal 3-vertices, then it contains 
a circle of type 63636. This implies that $\Si_b L$ contains a 2-sphere 
spanned by a simplex of type 364 as desired.
If $\Si_c\Si_b L$ contains a pair of antipodal 5- and 7-vertices, 
then we apply Lemma~\ref{lem:Dn_n-2sph} (for $n=k=3$) to the $A_3$-factor of
$\Si_c\Si_b B$ and conclude that $\Si_c\Si_b L$ contains a circle of type 7675657. 
We get again the 2-sphere in $\Si_b L$ as we wanted.
\qed

\bigskip\no
\emph{End of proof of Lemma~\ref{lem:no8B3pt}.}
Sublemma~\ref{sublem:no8B3pt} implies that $\Si_\zeta\Si_{x_1}K$
contains a 3-sphere spanned by a simplex of type 1467.
Recall the notation of the configuration $\ast$. 
Let $u$ be the 8-vertex $m(x_3,y_3)$. The segment $x_1z_1u$ is of type 828.
Then, it follows that
$\Si_{z_1} K$ contains a singular 4-sphere spanned by a simplex of type
14678. This sphere is of type 757 (to verify this, one can consider
the vector space realization of $D_n$ in Appendix~\ref{app:coxeter}).
Lemma~\ref{lem:4sph2pta} implies that the segment $x_1u$ cannot be extended
beyond $u$ in $K$. This implies in turn, that $CH(x_1,x_2,x_3)$
cannot be a spherical triangle. In particular $x_1$ must be an $8T$-vertex.
i.e.\ $8B_3\Rightarrow 8T$.

Let now $x_1,x_2,x_3$ be $8B_3$-vertices as in the definition of the property
$T'$. Our argument above shows that $\Si_{z_i} K$ contains a 4-sphere of type
757 for $i=1,2,3$. We apply again Lemma~\ref{lem:4sph2pta} and see that 
$rad(z_i,\{8\text{-vert. in } K\})\leq \arccos(-\frac{1}{2\sqrt{2}})$
for $i=1,2,3$.
Hence, $x_1$ is an $8T'$-vertex. A contradiction to Lemma~\ref{lem:8T'}.
\qed

If we combine the Lemmata \ref{lem:exist8B3pt} and \ref{lem:no8B3pt} 
we obtain the following result, which is the main step towards the proof of
Theorem~\ref{thm:ccE8}.

\begin{cor}\label{lem:no8Apt}
All $8$-vertices in $K$ have antipodes in $K$.
\end{cor}

Now we proceed to prove that the other vertices in $K$ must also have
antipodes in $K$. We will make constant use of 
the information about types of segments between
vertices in the Coxeter complex of type $E_8$ listed in 
Section~\ref{sec:E8coxeter}.

\begin{lem}\label{lem:no2Apt}
 All $2$-vertices in $K$ have antipodes in $K$.
\end{lem}
\proof
First note that a $2A$-vertex $x\in K$ cannot be adjacent to an 8-vertex in $K$.
Otherwise let $y\in K$ be an antipode of the 8-vertex adjacent to $x$. The segment
$xy$ is of type 2828. This in not possible due to Lemma~\ref{lem:antipodes} and
\ref{lem:no8Apt}.

Suppose there is a $2A$-vertex $x\in K$.
There exists $x'\in G\cdot x$ with $d(x,x')>\pihalf$. 
From the observation above it follows,
that $d(x,x')\neq\ac(-\quart)$. 
$d(x,x')$ cannot be $\ac(-\frac{3}{4})$ either, because in this case
the midpoint of the segment $xx'$ is an 8-vertex. 
It follows that $d(x,x')=\2pithird$ and
the segment $xx'$ is of type 26262. Let $y:=m(x,x')$, it is
also a $2A$-vertex. Therefore we can find $y'\in G\cdot y$ with $d(y,y')=\2pithird$.
Suppose w.l.o.g.\ that $\angle_y(x,y')\geq\pihalf$. Then $d(x,y')>\pihalf$, thus
$d(x,y')=\2pithird$. This implies by triangle comparison, that 
$\angle_y(x,y')\leq\ac(-\third)$.

\parpic{\includegraphics[scale=0.6]{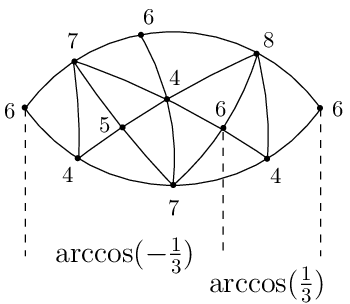}}
If $\angle_y(x,y')=\ac(-\third)$, then either this angle is of type 686, 
which is not possible because $K$ contains no 8-vertex adjacent to $y$;
or the simplicial convex hull of the segment $\ora{yx}\ora{yy'}$ contains
a 7-vertex adjacent to $\ora{yx}$. The segment connecting $\ora{yx}$ and $\ora{yx'}$
through this 7-vertex is of type 67686. This cannot happen either. 
Hence, $\angle_y(x,y')=\pihalf$. 
It follows that $\angle_y(x',y')\geq\pihalf$ and we conclude
analogously that $\angle_y(x',y')=\pihalf$.

\parpic{\includegraphics[scale=0.6]{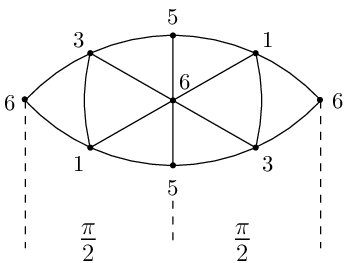}}
Let $\gamma\subset \Si_y K$ be the geodesic
connecting $\ora{yx}$ and $\ora{yx'}$ through $\ora{yy'}$. The simplicial convex hull
of $\gamma$ is either 3-dimensional, in which case the direction
$\ora{\ora{yx}\ora{yy'}}$ spans a simplex of type 578 and in particular, 
$\Si_y K$ contains 8-vertices, but this is not possible; 
or it is 2-dimensional and it contains a pair of 1-vertices adjacent to $\ora{yy'}$.
Let $z:=m(y,y')$ and let $w$ be the 6-vertex $m(y,z)$. The segment joining 
$\ora{wy}$ and $\ora{wz}$ through the 1-vertex adjacent to $\ora{wy}$ is of type 2152. 
It follows that $\ora{zy}$ is adjacent to a 5-vertex. The geodesic
connecting $\ora{zy}$ and $\ora{zy'}$ through this 5-vertex is of type 65856, but $z$
is a $2A$-vertex and $\Si_z K$ cannot contain 8-vertices. A contradiction.
\qed

\begin{lem}
 All $7$-vertices in $K$ have antipodes in $K$.
\end{lem}
\proof
Considering the singular circles in $E_8$,
we observe again that a $7A$-vertex cannot be adjacent to 2- or 8-vertices in $K$
by Lemma~\ref{lem:antipodes}.
Suppose $K$ contains $7A$-vertices, then there exist $7A$-vertices
$x_1,x_2\in K$ at distance $>\pihalf$. There are only two types of
segments $x_1x_2$ of length $>\pihalf$ and such that the simplices containing
$\ora{x_ix_{3-i}}$ in their interiors have no 2- or 8-vertices
(see the table in Appendix~\ref{app:E8coxeter}, page~\pageref{page:7vertices}). 
They are of type 76867 and 7342437. These
segments have a vertex of type 2 or 8 in their interiors, 
which yields a contradiction to Lemma~\ref{lem:antipodes}.
\qed

\begin{lem}
 All $1$-vertices in $K$ have antipodes in $K$.
\end{lem}
\proof
Suppose $x$ is an $1A$-vertex in $K$. 
By an argument similar to the one at the beginning of the proof 
of Lemma~\ref{lem:no2Apt}, we see that 
 $x$ cannot be adjacent to 2-, 7- or 8-vertices in $K$
by Lemma~\ref{lem:antipodes}.
Let $x'\in G\cdot x$ be another $1A$-vertex at distance $>\pihalf$ to $x$.
It follows that the simplex spanned by the direction $\ora{xx'}$ has no 2-, 7- or
8-vertices.
There are four possible types of segments $xx'$ 
(see table in Section~\ref{sec:E8coxeter}, page~\pageref{page:1vertices}).
If $d(x,x')=\ac(-\frac{3}{8})$, then the simplicial convex hull of $xx'$ contains an
8-vertex adjacent to $x'$. 
If $d(x,x')=\2pithird$ or $\ac(-\frac{7}{8})$, then the midpoint of $xx'$ is an 8-vertex.
If $d(x,x')=\ac(-\frac{5}{8})$, then the midpoint of $xx'$ is a 7-vertex.
This is not possible by Lemma~\ref{lem:antipodes}.
Hence, there are no $1A$-vertices in $K$.
\qed

\begin{lem}
 All $6$-vertices in $K$ have antipodes in $K$.
\end{lem}
\proof
Let $x$ be a $6A$-vertex. By the previous lemmata 
and according to the list of singular 
1-spheres in the Coxeter complex of type $E_8$, $x$ cannot be adjacent to
vertices of type 1, 2, 7 or 8. 
There exists another $6A$-vertex $x'\in K$ at distance
$>\pihalf$ to $x$. It follows that the direction $\ora{xx'}$ span
a simplex with no 1, 2, 7 or 8-vertices. Hence
$d(x,x')\in\{\ac(-\quart), \2pithird, \ac(-\frac{3}{4})\}$
(see table in Section~\ref{sec:E8coxeter}, page~\pageref{page:6vertices}).
In the first case the midpoint of $xx'$ is an 8-vertex and in the third
case, it is a 7-vertex.  In the second case the simplicial convex hull
of $xx'$ contains an 8-vertex adjacent to $x'$. 
A contradiction.
\qed

\begin{lem}
 All $3$-vertices in $K$ have antipodes in $K$.
\end{lem}
\proof
Observe first,
that a $3A$-vertex cannot be adjacent to a vertex of type 1, 2, 6, 7 or 8. 
\parpic{\includegraphics[scale=0.7]{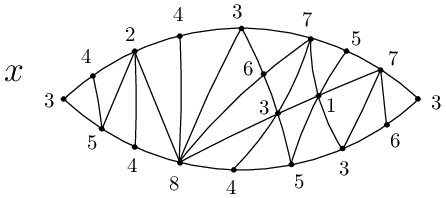}}
If $K$ contains $3A$-vertices, then it contains at least
two distinct $3A$-vertices $x,x'$.
Then $\ora{xx'}$ is contained in an
edge of type 45, otherwise $CH(x,x')$ would contain vertices of type 1, 2, 6, 7 or 8
adjacent to $x$. 
Consider the bigon in the Coxeter complex of type $E_8$, which
is the convex hull of a simplex of type 345 and the antipode of the 3-vertex
of this simplex (cf. Remark~\ref{rem:1-2-spheres}). 
We see from this bigon (see figure) that 
there are only three possibilities 
for the type of the segment $xx'$. 
In one of them, the midpoint of $xx'$ is a 2-vertex; and
in another possibility, it is an 8-vertex. The simplicial convex hull of $xx'$ 
for the last possibility contains an 8-vertex adjacent to $x'$.
We obtain again a contradiction to Lemma~\ref{lem:antipodes}.
\qed

\begin{lem}
 All $4$- and $5$-vertices in $K$ have antipodes in $K$.
\end{lem}
\proof
A vertex in $K$ of type 4 or 5 without antipodes in $K$ cannot have vertices of type
1, 2, 3, 6, 7 or 8 in $K$ adjacent to it. It follows that, 
if $K$ contains $4A$- or $5A$-vertices,
then it has dimension $\leq 1$. A contradiction.
\qed

We have shown in the previous lemmata that all vertices of a counterexample 
$K$ have antipodes in $K$, 
contradicting Proposition~\ref{lem:Ksubbuild}. This proves our main result:

\begin{thm}\label{thm:ccE8}
 The Center Conjecture~\ref{centerconj} holds for spherical buildings of type $E_8$.
\end{thm}

\begin{rem}\label{rem:ccE8}
Our proof actually shows that $K$ is a subbuilding or the action of the group 
$Aut_B(K)\acts K$ fixes a point, where $Aut_B(K)\supseteq Stab_{Aut(B)}(K)$
denotes the possibly larger
group of isometries of $K$ preserving the polyhedral structure of $K$ induced 
by the polyhedral structure of $B$ and such that the
changes of vertex types of $K$ correspond to a symmetry of the Dynkin diagram of $B$. 
(Compare with \cite[Sec.\ 3.5]{LeebRamos}.)

Notice that the automorphisms in $Aut_B(K)$ are not necessarily extendable to
automorphisms of $B$, as the following example shows. Let $\sigma\subset B$
be a panel (i.e.\ a face of codimension one) and let $K_\sigma$ 
be the union of the Weyl chambers in $B$ containing $\sigma$.
It is a convex subcomplex of $B$ and $Aut_B(K_\sigma)$ is the group of permutations
of the set of Weyl chambers containing $\sigma$. This group is very large
if e.g. the set of Weyl chambers containing $\sigma$ is uncountable.

Although the automorphisms in $Aut_B(K)$ must not be extendable
to automorphisms of $B$, the group $Aut_B(K)$ depends on the
ambient building $B$, or more precisely,
on the Coxeter complex on which $B$ is modelled,
 in the sense illustrated by the following example.
Let $B$ be a building of type $F_4$ and 
let $K\subset B$ be a convex subcomplex.
Suppose we can embed $B$ in a building $\widetilde B$ of type $E_8$, such
that the polyhedral structure of $B$ coincides with the structure
induced by the polyhedral structure of $\widetilde B$
(the Coxeter complex of type $F_4$ can be viewed as a subcomplex
of the Coxeter complex of type $E_8$).
The image of $K$ under this embedding
is a convex subcomplex of $\widetilde B$.
Then $Aut_{\widetilde B}(K)$ is the possible smaller
subgroup of $Aut_B(K)$ of
type preserving automorphisms. This is a consequence of the fact
that the Dynkin diagram of $F_4$ has one symmetry whereas 
the Dynkin diagram of $E_8$ has none.
\end{rem}

\subsection{The $E_6$-case}

\begin{thm}\label{thm:ccE6}
The Center Conjecture~\ref{centerconj} holds for spherical buildings of type $E_6$.
\end{thm}
\proof
Let $K$ be a convex subcomplex
of a spherical building $B$ of type $E_6$
and suppose that $K$ is not a subbuilding.

Let $\kappa$ be a circle of radius 1 with 
the structure of the spherical Coxeter complex of type $I_2(6)$ 
with labelling of its Dynkin diagram \includegraphics[scale=0.4]{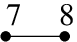}.
Let $p_1,p_2$ be a pair of antipodal 7-vertices in $\kappa$.

\begin{center}
 \includegraphics[scale=0.7]{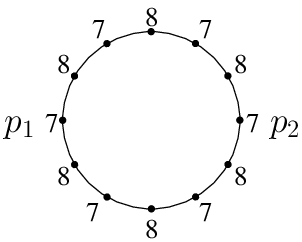}
\end{center}

Consider the spherical join $\wht B := B \circ \kappa$. 
There is a natural embedding $B\hookrightarrow \wht B$ so
we can regard $B$ as a subset of $\wht B$.
Let $\wht K:= K\circ \kappa \subset \wht B$.

Let $B_{p_i}:=\Si_{p_i} (B\circ\{p_i\})\subset \Si_{p_i}\wht B$ for $i=1,2$.
Then $\Si_{p_i}\wht B=B_{p_i}\circ\Si_{p_i}\kappa$ and we have
isometries $B \stackrel{\rho_i}{\rightarrow} B_{p_i}$ defined by
$\rho_i(v):=\ora{p_iv}$.
Let $K_{p_i}:=\Si_{p_i}(K\circ\{p_i\})\cong K$.

Let $B_1\stackrel{\rho}{\rightarrow} B_2$ be the isometry that sends a direction
$\xi\in\Si_{p_1}(B\circ\{p_i\})$ to the initial direction at $p_2$ of the geodesic
connecting $p_1$ and $p_2$ with initial direction $\xi$ at $p_1$.
Then $\rho=\rho_1^{-1}\circ\rho_2$.

Recall that the link of a 7-vertex in the Coxeter complex
of type $E_8$ is a Coxeter complex of type $E_6\circ A_1$.
We consider the building $B_{p_1}$ of type 
$E_6$ with the labelling of vertices
induced by the labelling of $B$ and the isometry $\rho_1$.
With this labelling, a chart 
$(S^5,W_{E_6}) \hookrightarrow B_{p_1}$ of the building 
$B_{p_1}$ induces a chart $(S^7,W_{E_8}) \hookrightarrow \wht B$, giving
$\wht B$ a structure of spherical building of type $E_8$,
where the induced polyhedral structure of $\kappa$ coincides with
its structure as Coxeter complex of type $I_2(6)$.
The labelling of the vertices of $B_{p_2}$ induced by this building
structure in $\wht B$ can be obtained from the one induced by $\rho_2$
by exchanging the labels $2\leftrightarrow 6$, $3\leftrightarrow 5$
and fixing 1 and 4, that is, by the symmetry of the Dynkin diagram
of type $E_6$.

As an example, we present the $E_8$-structure of a 2-dimensional hemisphere 
$\{v\}\circ\kappa\subset \wht B$, where $v$ is a 2-vertex of $B$:
\begin{center}
 \hpic{\includegraphics[scale=0.9]{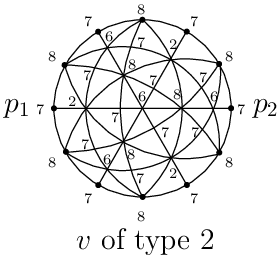}}
\end{center}

Notice that there is a natural isomorphism  
$Aut_0(B_{p_1}) \cong Fixator_{Aut(\wht B)}(\kappa)$
between the type preserving automorphisms of $B_{p_1}$
and the automorphisms of $\wht B$ fixing $\kappa$ pointwise.
It extends to an embedding 
$Aut(B_{p_1})\stackrel{\iota}\hookrightarrow Stab_{Aut(\wht B)}(\kappa)$,
where the image of a non type preserving automorphism of $B_{p_1}$
restricted to $\kappa$ is the antipodal involution $ant_\kappa$
of $\kappa$.
In particular, $\iota(\varphi)(p_1)=p_2$
for a non type preserving automorphism $\varphi\in Aut(B_{p_1})$, 
hence
$\iota(\varphi)$ induces an isometry $B_{p_1}\rightarrow B_{p_2}$.
This isometry is type preserving and
coincides with $\rho\circ\varphi:B_{p_1}\rightarrow B_{p_2}$.
This means that the image $\iota(Aut(B_{p_1}))$
acts on $B_{p_1}$ via $\iota(\varphi)|_{B_{p_1}}$,
if $\varphi\in Aut(B_{p_1})$ is type preserving; 
and via $\rho^{-1}\!\circ(\iota(\varphi)|_{B_{p_1}})$,
if $\varphi\in Aut(B_{p_1})$ is not type preserving.
The embedding $\iota$ restricts to an embedding
$G:=Stab_{Aut(B_{p_1})}(K_{p_1})\stackrel{\iota}{\hookrightarrow} G':=Stab_{Aut(\wht B)}(\wht K)$.

There is an isometry $\phi_0$ of $\wht B$ that rotates $\kappa$ an angle
of $\2pithird$ and preserves every 2-dimensional hemisphere bounded by $\kappa$.
The restriction of $\phi_0$ to an apartment $a\subset \wht B$ is the
composition of the reflections on two walls orthogonal to two 8-vertices
in $\kappa$ at distance $\pithird$. It follows that $\phi_0$ is an automorphism
of $\wht B$ and $\phi_0\in G'$.

We apply now the Center Conjecture for buildings of type $E_8$ (Theorem~\ref{thm:ccE8})
to $\wht K\subset \wht B$. 
Since $K$ is not a subbuilding, then $\wht K$ cannot be a subbuilding.
It follows that
$G'$ fixes a point $x\in\wht K$. But since $\phi_0\in G'$ and $\phi_0$ has no fixed
points in $\kappa$, it follows that $x\not\in \kappa$.
This implies that $\iota(G)$ preserves the 2-dimensional hemisphere $h\subset \wht K$ 
bounded by $\kappa$ and containing $x$ and in particular, it fixes its center.
Hence, it preserves the geodesic $\gamma$ connecting $p_1$ and $p_2$ 
through the center of $h$.
It follows that $\iota(G)\acts K_{p_1}$ fixes the initial direction of 
$\gamma$ at $p_1$.
\qed

\begin{rem}\label{rem:ccE6}
Notice that the embedding $G\hookrightarrow G'$ in the proof of 
Theorem~\ref{thm:ccE6} extends to an embedding
$Aut_{B_{p_1}}(K_{p_1})\hookrightarrow Aut_{\widetilde B}(\widetilde K)$. 
Then by Remark~\ref{rem:ccE8},
the proof actually shows that $K$ is a subbuilding or the action of the group 
$Aut_B(K)\acts K$ fixes a point.
\end{rem}

\subsection{The $E_7$-case}

\begin{thm}\label{thm:ccE7}
The Center Conjecture~\ref{centerconj} holds for spherical buildings of type $E_7$.
\end{thm}
\proof
It can be deduced from the $E_8$-case as follows: Let $K\subset B$ be a convex subcomplex
of a spherical building of type $E_7$. Suppose that $K$ is not a subbuilding.
Let $\widetilde B$ be the suspension of $B$, 
i.e.\ the spherical join of $B$ and a 0-sphere $\{p_1,p_2\}$.
There is a natural embedding $B\hookrightarrow\widetilde B$, 
so we can consider $B$ as a subset of $\widetilde B$.
Notice that the map $v\mapsto \ora{p_iv}$
for $v\in B\subset \widetilde B$
is an isometry $B \cong B_{p_i}:=\Si_{p_i} \widetilde B$.
Let $\widetilde K\subset \widetilde B$ be the suspension of $K$ and
let $K_{p_i}:=\Si_{p_i} \widetilde K\cong K$.

Recall that the link of an 8-vertex in the Coxeter complex
of type $E_8$ is a Coxeter complex of type $E_7$. Hence
a chart $(S^6,W_{E_7}) \hookrightarrow B_{p_1}$
of the building $B_{p_1}$ induces a chart 
$(S^7,W_{E_8}) \hookrightarrow \widetilde B$, giving
$\widetilde B$ a structure of spherical building of type $E_8$,
where $p_1$ and $p_2$ are 8-vertices.

Observe that there is a natural isomorphism $Aut(B_{p_1})\cong Stab_{Aut(\widetilde B)}(p_1)$.
The embedding $Aut(B_{p_1})\hookrightarrow Aut(\widetilde B)$ restricts to an embedding
$G:=Stab_{Aut(B_{p_1})}(K_{p_1})\hookrightarrow \widetilde G:= Stab_{Aut(\widetilde B)} (\widetilde K)$.

There is an isometry $\phi_0$ of $\widetilde B$ that exchanges the points
$p_1\leftrightarrow p_2$ and preserves the geodesics connecting $p_1$ and $p_2$.
The restriction of $\phi_0$ to an apartment of $\widetilde B$ is the reflection on the wall
orthogonal to $p_1,p_2$. Hence $\phi_0$ is an automorphism of $\widetilde B$ and
$\phi_0\in \widetilde G$.

We apply now the Center Conjecture for buildings of type $E_8$ (Theorem~\ref{thm:ccE8})
to the building $\widetilde B$ and the convex subcomplex $\widetilde K$.
Since $K$ is not a subbuilding, then $\widetilde K$ is not a subbuilding either.
It follows that $\widetilde G$ fixes a point $x\in\widetilde K$ and since
$\phi_0\in \widetilde G$, this fixed point cannot be $p_1$ or $p_2$.
The image of $G$ in $\widetilde G$ 
fixes $x$, $p_1$ and $p_2$, hence it fixes pointwise the geodesic 
$\gamma\subset \widetilde K$ through $x$ connecting $p_1$
and $p_2$. 
Therefore, the action $G\acts K_{p_1}\cong K$ 
has a fixed point $\ora{p_1x} \in K_{p_1}$.
\qed

\begin{rem}\label{rem:ccE7}
Notice that the embedding $G\hookrightarrow \widetilde G$ in the proof of 
Theorem~\ref{thm:ccE7} extends to an embedding 
$Aut_{B_{p_1}}(K_{p_1})\hookrightarrow Aut_{\widetilde B}(\widetilde K)$. 
Then by Remark~\ref{rem:ccE8},
the proof actually shows that $K$ is a subbuilding or the action of the group 
$Aut_B(K)\acts K$ fixes a point.
\end{rem}

\appendix

\section{Vector space realizations of Coxeter complexes}\label{app:coxeter}

In this appendix we present a vector space realization of the
irreducible spherical Coxeter complexes. 
The information on
the root systems can be found in \cite[Ch. 5]{GroveBenson}.
The orders of the irreducible Weyl groups can be found
in \cite[p. 80]{GroveBenson}.

We consider the spherical Coxeter complex $(S^{n-1},W)$ embedded in $\R^{n}$
as the unit sphere. 
Let $\{e_i\}_{i=1}^{n}$ denote the canonical base of $\R^{n}$.

The {\em root system} of a Coxeter complex $(S,W)$ is the set of (unit) vectors
orthogonal to the hyperplanes inducing the reflections in $W$. The elements
of the root system are called root vectors.

A subset $F$ of the root system is called a base if there is a vector $v\in\R^n$
such that $\langle r, v\rangle \neq 0$ for all root vectors $r$, and $F$ is minimal 
with respect to
the property that any root vector $r$, such that $\langle r, v\rangle > 0$,
can be written as a linear combination of elements in $F$ with nonnegative
coefficients. 
The {\em fundamental root vectors} are the elements of a given base of
the root system. 
The {\em fundamental Weyl chamber} of $(S,W)$ is 
$\triangle:=\bar\triangle \cap S$, where $\bar\triangle$
is the intersection
of the half spaces $\langle r_i, \cdot \rangle\geq 0$, where 
$r_1,\dots,r_{n}$ are the fundamental root vectors. 
$\bar\triangle$ is a fundamental domain for the action of $W$ in $\R^n$.

Let $v_i$ be the vertex of $\triangle$ opposite to the face determined
by $\langle r_i, \cdot \rangle = 0$. We say that a vertex of $(S,W)$
is of type $i$, if it lies on the orbit $W\cdot v_i$.

We use the following labelling of the Dynkin diagrams of the irreducible
spherical Coxeter complexes:
\begin{center}
\begin{tabular}{C{3cm}C{4cm}C{3cm}C{4cm}}
$I_2(m)$ & \includegraphics[scale=0.5]{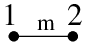} 
	& $H_3$ & \includegraphics[scale=0.5]{H3dynk}\\
$A_n$ & \includegraphics[scale=0.5]{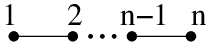}
	& $H_4$ & \includegraphics[scale=0.5]{H4dynk}\\
$B_n$, $n\geq 3$ & \includegraphics[scale=0.5]{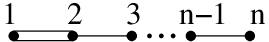}
	& $E_6$ & \hpic{\includegraphics[scale=0.5]{E6dynk}}\\
$D_n$, $n\geq 4$ &\hpic{\includegraphics[scale=0.5]{Dndynk}}
	 & $E_7$ & \hpic{\includegraphics[scale=0.5]{E7dynk}}\\
$F_4$ & \includegraphics[scale=0.5]{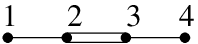}
	& $E_8$ & \hpic{\includegraphics[scale=0.5]{E8dynk}}\\
\end{tabular}
\end{center}

Recall, that the link $\Si_v S$ of a vertex $x\in S$ 
is a spherical Coxeter complex with Weyl group $Stab_W(v)$ and with
Dynkin diagram obtained from the Dynkin diagram of $(S,W)$ by deleting the
vertex with label corresponding to the type of $v$.

The antipodal involution $v\mapsto -v$ is type preserving for the
spherical Coxeter complexes of type $I_2(m)$ ($m$ even), $B_n$ ($n\geq 3$),
$D_{2n}$ ($n\geq 2$), $H_3$, $H_4$, $F_4$, $E_7$ and $E_8$.
It exchanges the types $1\leftrightarrow 2$ in $I_2(m)$ for $m$ odd;
the types $i\leftrightarrow (n+1-i)$, for $i=1,\dots,[\frac{n}{2}]$ in $A_n$;
the types $1\leftrightarrow 2$, in $D_{2n+1}$, $n\geq 2$ and
the types $2\leftrightarrow 6$ and $3\leftrightarrow 5$ in $E_6$.

In order to deduce the geometric properties of Coxeter complexes
listed in Section~\ref{sec:coxeter} one needs to do several computations in
the vector space realizations of the Coxeter complexes. 
We will not write down each computation but instead explain the general
ideas of how to do these computations in the following remarks.

\begin{rem}[Singular 1- and 2-spheres]\label{rem:1-2-spheres}
Suppose $xy$ is an edge of $S$ of type $ij$. By deleting the vertex with label $j$
from the Dynkin diagram of $(S,W)$, 
we obtain the Dynkin diagram of $(\Si_y S, Stab_W(y))$.
We can easily read off this Dynkin diagram which type the antipode of $\ora{yx}$
in $\Si_y S$ has.
Say it has type $k$, then the edge $xy$ extends to a segment of type $ijk$. 
Repeating this procedure and taking into account the lengths of the different types of
segments (which can be deduced from the description of the fundamental
Weyl chamber), we can determine the different singular 1-spheres in $S$.
A similar consideration can be used to determine the 2-dimensional singular bigons
bounded by singular segments and with it the 2-dimensional singular spheres.
\end{rem}

\begin{rem}[Types of segments between vertices]\label{rem:segmenttype}
To determine the different types of segments modulo the action of the Weyl group
connecting a vertex of type $i$ with a vertex of type $j$, it suffices to compute
the vertices of type $j$ in the {\em spherical bigon}
$\beta_i:=CH(\triangle, \wht{v_i})\subset S$, where
$\wht{v_i}$ is the vertex antipodal to $v_i$. 

The bigon $\beta_i$ can be described by the set of inequalities
$\{\langle r_l, \cdot \rangle\geq 0\}_{l\neq i}$.
\end{rem}

\begin{rem}[Types of segments between vertices with restricted
initial direction]\label{rem:segmenttypeII}
More generally, suppose we want to determine the different types 
of segments connecting a vertex $x$ 
of type $i$ and a vertex $y$ of type $j$, such that the vertices of the simplex 
in $\Si_x S$ spanned by
the direction $\ora{xy}$ are not of type $i_1,\dots,i_k \neq i$.
Then, it suffices to compute the vertices of type $j$ in the spherical bigon 
$\beta_i(i_1,\dots,i_k):=CH(\triangle(i_1,\dots,i_k),\wht{v_i})$. 
Here, $\triangle(i_1,\dots,i_k)$ denotes the face of the fundamental
Weyl chamber $\triangle$, which does not contain the vertices $v_{i_1},\dots,v_{i_k}$.

The bigon $\beta_i(i_1,\dots,i_k)$ can be described by the set of (in)equalities
$$\{\langle r_l, \cdot \rangle\geq 0\}_{l\neq i,i_1,\dots,i_k},\;\;\;
\{\langle r_l, \cdot \rangle = 0\}_{l=i_1,\dots,i_k}.$$
\end{rem}

\begin{rem}\label{rem:tables}
Given a table listing the $j$-vertices in the bigon $\beta_i$, 
this list can be verified
as follows. First, we have to check that the vertices listed indeed are of type $j$
and are contained in $\beta_i$. Next we notice that $\beta_i$ is a fundamental
domain for the action $Stab_W(v_i)\acts S$. 
For a $j$-vertex $x$ in the list, let $\sigma_x$ be the face of $\triangle$ spanned by
the initial part of the segment $v_ix$.
Then the orbit $Stab_W(v_i)\cdot x$ has cardinality $|Stab_w(v_i)|/|Stab_W(\sigma_x)|$.
Since the stabilizers are again Weyl groups of spherical Coxeter complexes, 
their orders can be found in the table in \cite[p. 80]{GroveBenson}.
It remains to verify that the union of the orbits $Stab_W(v_i)\cdot x$ exhausts all
the $j$-vertices in $S$, that is, we have to check that
$\sum\limits_{x\text{ in the list}} \frac{|Stab_W(v_i)|}{|Stab_W(\sigma_x)|} 
= \frac{|W|}{|Stab_W(v_i)|}$.
\end{rem}

\subsection{$D_n$}\label{app:Dn}

Let $n\geq 4$. The Weyl group $W_{D_n}$ of type $D_n$
is the finite group of isometries of 
$\R^n$ 
generated by the reflections at the hyperplanes orthogonal 
to the {\em fundamental root vectors}:
$$
r_1=e_1+e_2,\;\;\;\; 
r_i=e_i-e_{i-1} \text{ for } 2\leq i\leq n
$$

The {\em fundamental Weyl chamber} $\triangle$ can be described by the inequalities:
$$
-x_2\overset{\text{(1)}}{\leq} x_1
\overset{\text{(2)}}{\leq} x_2
\overset{\text{(3)}}{\leq}\dots 
\overset{\text{(n)}}{\leq}x_n.
$$

Next we exhibit an element representing the vertices of the fundamental Weyl chamber $\triangle$,
i.e.\ elements of $\R^+\!\cdot v_i$:
\begin{center}
\begin{tabular}{R{3cm}C{1cm}C{6cm}}
1-vertex: & $v_1$	&	$(\phm1,1,1,\dots,1)$\\
2-vertex: & $v_2$	&	$(-1,1,1,\dots,1)$\\
3-vertex: & $v_3$	&	$(\phm0,0,1,\dots,1)$\\
$\vdots$ \phantom{vert} & $\vdots$	&	$\vdots$\\
$(n-1)$-vertex: & $v_{n-1}$ 	&	$(\phm0,\dots,0,1,1)$\\
$n$-vertex: & $v_n$&	$(\phm0,\dots,0,0,1)$\\
\end{tabular}
\end{center}

The Weyl group $W_{D_n}$ acts on $\R^n$ by permutations of the coordinates
and change of signs in an even number of places.

\subsection{$E_6$}\label{app:E6}

The Weyl group $W_{E_6}$ of type $E_6$
is the finite group of isometries of 
$\R^6\cong\{(x_1,\dots,x_8)\in\R^8\; |\; x_6=x_7=x_8\}$ 
generated by the reflections at the hyperplanes orthogonal 
to the {\em fundamental root vectors}:
$$
r_1=\half(1,1,1,-1,-1,-1,-1,-1),\;\;\;\; r_i=e_i-e_{i-1} \text{ for } 2\leq i\leq 5;.
$$$$
\text{ and }\;\;\; r_6=\half(1,1,1,1,-1,1,1,1).
$$

This choice of fundamental root vectors differs from the one in \cite{GroveBenson}.
It is the same as the one presented in \cite[Sec. 2.2.3]{LeebRamos}, where its equivalency
is explained.

The {\em fundamental Weyl chamber} $\triangle$ can be described by the inequalities:
$$
x_4+x_5+\dots+x_8 \overset{\text{(1)}}{\leq} x_1+x_2+x_3\;;\;\;\; \;\;
x_1\overset{\text{(2)}}{\leq} x_2
\overset{\text{(3)}}{\leq}\dots 
\overset{\text{(5)}}{\leq}x_5\;;\;\;\; \;\;
x_5\overset{\text{(6)}}{\leq}x_1+\dots+x_4+x_6+x_7+x_8.
$$

Next we exhibit an element representing the vertices of the fundamental Weyl chamber $\triangle$,
i.e.\ elements of $\R^+\!\cdot v_i$:
\begin{center}
\begin{tabular}{p{3cm}p{6cm}}
1-vertex: $v_1$	&	$(\phm1,\phm1,\phm1,\phm1,\phm1,-1,-1,-1)$\\
2-vertex: $v_2$	&	$(-3,\phm3,\phm3,\phm3,\phm3,-1,-1,-1)$\\
3-vertex: $v_3$	&	$(\phm0,\phm0,\phm3,\phm3,\phm3,-1,-1,-1)$\\
4-vertex: $v_4$	&	$(\phm1,\phm1,\phm1,\phm3,\phm3,-1,-1,-1)$\\
5-vertex: $v_5$	&	$(\phm3,\phm3,\phm3,\phm3,\phm9,-1,-1,-1)$\\
6-vertex: $v_6$	&	$(\phm3,\phm3,\phm3,\phm3,\phm3,\phm1,\phm1,\phm1)$\\
\end{tabular}
\end{center}

We list now the orbits of the 2-vertices of $\triangle$
 under the action of the Weyl group
(modulo the following elements of the Weyl group: 
permutations of the first five coordinates and 
change of sign in an even number of places
in the first five coordinates). 
We give representing vectors for the vertices. The 6-vertices are just
the antipodes of the 2-vertices.

\begin{center}
\begin{tabular}{C{2cm} R{6.5cm} R{6.5cm}}
2-vertices	&	$(-3,\phm3,\phm3,\phm3,\phm3,-1,-1,-1)$,	&
				$(\phm0,\phm0,\phm0,\phm0,\phm6,\phm2,\phm2,\phm2)$,\\
		&	$(\phm0,\phm0,\phm0,\phm0,\phm0,-4,-4,-4)$. \\ \\
\end{tabular}
\end{center}

This list can be verified by checking that the vertices listed indeed lie on the orbit 
$W_{E_6}\cdot v_2$ and there are as many as $|W_{E_6}|/|W_{D_5}|=3^3$.

We describe in the following table
the 2- and 6-vertices $x$ in $\beta_2$. Let $\sigma$ be the face of
$\Si_{v_2} \triangle$ containing $\ora{v_2x}$ in its interior.
\begin{center}
\begin{tabular}{C{3cm} | C{6.5cm} | C{2.5cm} | C{2cm}}
&$x$		&	$d(x,v_2)$	&	Type of $\sigma$	\\\hline
2-vertices $x\neq v_2$	&
	$(\phm3,-3,\phm3,\phm3,\phm3,-1,-1,-1)$\hfill		\linebreak
	$(\phm6,\phm0,\phm0,\phm0,\phm0,\phm2,\phm2,\phm2)$	&
	$\ac(\quart)$\hfill					\linebreak
	$\2pithird$						&
	3\hfill							\linebreak
	6	\\
\hline
6-vertices $x\neq \wht{v_2}$	&
	$(\phm3,\phm3,\phm3,\phm3,\phm3,\phm1,\phm1,\phm1)$\hfill\linebreak
	$(\phm6,\phm0,\phm0,\phm0,\phm0,-2,-2,-2)$		&
	$\pithird$\hfill\linebreak
	$\ac(-\quart)$	&
	6\hfill\linebreak
	1	\\
\end{tabular}
\end{center}

\subsection{$E_7$}\label{app:E7}

The Weyl group $W_{E_7}$ of type $E_7$
is the finite group of isometries of 
$\R^7\cong\{(x_1,\dots,x_8)\in\R^8\; |\; x_7=x_8\}$ 
generated by the reflections at the hyperplanes orthogonal 
to the {\em fundamental root vectors}:
$$
r_1=\half(1,1,1,-1,-1,-1,-1,-1),\;\;\;\; r_i=e_i-e_{i-1} \text{ for } 2\leq i\leq 6
$$$$
\text{ and }\;\;\; r_7=\half(1,1,1,1,1,-1,1,1).
$$

This choice of fundamental root vectors differs from the one in \cite{GroveBenson}.
They are equivalent, as can be seen by an argument similar to
the one in \cite[Sec. 2.2.3]{LeebRamos} for the Coxeter complex of type $E_6$.

The {\em fundamental Weyl chamber} $\triangle$ can be described by the inequalities:
$$
x_4+x_5+\dots+x_8 \overset{\text{(1)}}{\leq} x_1+x_2+x_3\;;\;\;\; \;\;
x_1\overset{\text{(2)}}{\leq} x_2
\overset{\text{(3)}}{\leq}\dots 
\overset{\text{(6)}}{\leq}x_6\;;\;\;\; \;\;
x_6\overset{\text{(7)}}{\leq}x_1+\dots+x_5+x_7+x_8.
$$

Next we exhibit an element representing the vertices of the fundamental Weyl chamber $\triangle$,
i.e.\ elements of $\R^+\!\cdot v_i$:
\begin{center}
\begin{tabular}{p{3cm}p{6cm}}
1-vertex: $v_1$	&	$(\phm1,\phm1,\phm1,\phm1,\phm1,\phm1,-2,-2)$\\
2-vertex: $v_2$	&	$(-1,\phm1,\phm1,\phm1,\phm1,\phm1,-1,-1)$\\
3-vertex: $v_3$	&	$(\phm0,\phm0,\phm1,\phm1,\phm1,\phm1,-1,-1)$\\
4-vertex: $v_4$	&	$(\phm1,\phm1,\phm1,\phm3,\phm3,\phm3,-3,-3)$\\
5-vertex: $v_5$	&	$(\phm1,\phm1,\phm1,\phm1,\phm3,\phm3,-2,-2)$\\
6-vertex: $v_6$	&	$(\phm1,\phm1,\phm1,\phm1,\phm1,\phm3,-1,-1)$\\
7-vertex: $v_7$	&	$(\phm1,\phm1,\phm1,\phm1,\phm1,\phm1,\phm0,\phm0)$\\
\end{tabular}
\end{center}

We list now the orbits of the 2- and 7-vertices of $\triangle$
 under the action of the Weyl group
(modulo the following elements of the Weyl group:
permutations of the first six coordinates, 
change of sign in an even number of places
in the first six coordinates 
and simultaneous change of sign of the last two coordinates). 
We give representing vectors for the vertices.

\begin{center}
\begin{tabular}{C{2cm} R{6.5cm} R{6.5cm}}
2-vertices	&	$(-1,\phm1,\phm1,\phm1,\phm1,\phm1,-1,-1)$,	&
				$(\phm0,\phm0,\phm0,\phm0,\phm0,\phm0,\phm2,\phm2)$,\\
		&	$(\phm0,\phm0,\phm0,\phm0,\phm2,\phm2,\phm0,\phm0)$. \\ \\
7-vertices	&	$(\phm1,\phm1,\phm1,\phm1,\phm1,\phm1,\phm0,\phm0)$,	& 
				$(\phm0,\phm0,\phm0,\phm0,\phm0,\phm2,\phm1,\phm1)$.\\
\end{tabular}
\end{center}

This list can be verified by checking that the vertices listed indeed lie on the orbits 
$W_{E_7}\cdot v_i$ and there are as many as $|W_{E_7}|/|Stab_{W_{E_7}}(v_i)|$.

We describe in the following table
the 2-vertices $x$ in $\beta_2$. Let $\sigma$ be the face of
$\Si_{v_2} \triangle$ containing $\ora{v_2x}$ in its interior.
\begin{center}
\begin{tabular}{C{3cm} | C{6.5cm} | C{2.5cm} | C{2cm}}
&$x$		&	$d(x,v_2)$	&	Type of $\sigma$	\\\hline
2-vertices $x\neq v_2,\wht{v_2}$	&
	$(\phm1,-1,\phm1,\phm1,\phm1,\phm1,-1,-1)$\hfill		\linebreak
	$(\phm2,\phm0,\phm0,\phm0,\phm0,\phm2,\phm0,\phm0)$\hfill		\linebreak
	$(\phm2,-2,\phm0,\phm0,\phm0,\phm0,\phm0,\phm0)$	&
	$\pithird$\hfill					\linebreak
	$\pihalf$\hfill					\linebreak
	$\2pithird$						&
	3\hfill							\linebreak
	6\hfill							\linebreak
	3	\\
\end{tabular}
\end{center}

We describe in the following table
the 2- and 7-vertices $x$ in $\beta_7$. Let $\sigma$ be the face of
$\Si_{v_7} \triangle$ containing $\ora{v_7x}$ in its interior.
\begin{center}
\begin{tabular}{C{3cm} | C{6.5cm} | C{2.5cm} | C{2cm}}
&$x$		&	$d(x,v_7)$	&	Type of $\sigma$	\\\hline
2-vertices $x$	&
	$(-1,\phm1,\phm1,\phm1,\phm1,\phm1,-1,-1)$\hfill		\linebreak
	$(\phm0,\phm0,\phm0,\phm0,\phm0,\phm0,-2,-2)$\hfill		\linebreak
	$(-1,-1,-1,-1,-1,\phm1,-1,-1)$	&
	$\ac(\frac{1}{\sqrt{3}})$\hfill					\linebreak
	$\pihalf$\hfill					\linebreak
	$\ac(-\frac{1}{\sqrt{3}})$						&
	2\hfill							\linebreak
	1\hfill							\linebreak
	6	\\
\hline
7-vertices $x\neq v_7, \wht{v_7}$	&
	$(\phm0,\phm0,\phm0,\phm0,\phm0,\phm2,\phm1,\phm1)$\hfill\linebreak
	$(-2,\phm0,\phm0,\phm0,\phm0,\phm0,\phm1,\phm1)$		&
	$\ac(\third)$\hfill\linebreak
	$\ac(-\third)$	&
	6\hfill\linebreak
	2	\\
\end{tabular}
\end{center}

\subsection{$E_8$}\label{app:E8coxeter}

The Weyl group $W_{E_8}$ of type $E_8$ is the finite group of isometries of $\R^8$ 
generated by the reflections at the hyperplanes orthogonal 
to the {\em fundamental root vectors}:

$$r_1=\half(1,1,1,-1,-1,-1,-1,-1) \text{ and } r_i=e_i-e_{i-1} \text{ for } 2\leq i\leq 8.$$

The {\em fundamental Weyl chamber} $\triangle$ can be described by the inequalities:

$$
x_4+x_5+\dots+x_8 \overset{\text{(1)}}{\leq} x_1+x_2+x_3\;;\;\;\; \;\;
x_1\overset{\text{(2)}}{\leq} x_2
\overset{\text{(3)}}{\leq}x_3
\overset{\text{(4)}}{\leq}\dots 
\overset{\text{(8)}}{\leq}x_8.
$$

Next we exhibit an element representing the vertices of the fundamental Weyl chamber $\triangle$,
i.e.\ elements of $\R^+\!\cdot v_i$:
\begin{center}
\begin{tabular}{p{3cm}p{6cm}}
1-vertex: $v_1$	&	$(-1,-1,-1,-1,-1,-1,-1,-1)$\\
2-vertex: $v_2$	&	$(-3,-1,-1,-1,-1,-1,-1,-1)$\\
3-vertex: $v_3$	&	$(-2,-2,-1,-1,-1,-1,-1,-1)$\\
4-vertex: $v_4$	&	$(-5,-5,-5,-3,-3,-3,-3,-3)$\\
5-vertex: $v_5$	&	$(-2,-2,-2,-2,-1,-1,-1,-1)$\\
6-vertex: $v_6$	&	$(-3,-3,-3,-3,-3,-1,-1,-1)$\\
7-vertex: $v_7$	&	$(-1,-1,-1,-1,-1,-1,\phm0,\phm0)$\\
8-vertex: $v_8$	&	$(-1,-1,-1,-1,-1,-1,-1,\phm1)$\\
\end{tabular}
\end{center}

We list now (modulo the following elements of the Weyl group:
permutations of the coordinates and change of sign in an even number 
of places) the orbits of the vertices 
of $\triangle$ of type 1, 2, 6, 7, 8
under the action of the Weyl group. 
We give representing vectors for the vertices.
\begin{center}
\begin{tabular}{C{2cm} R{6.5cm} R{6.5cm}}
1-vertices	&	$(-1,-1,-1,-1,-1,-1,-1,-1)$,	&
				$\half(-5,\phm1,\phm1,\phm1,\phm1,\phm1,\phm1,\phm1)$,\\
		&	$\half(\phm3,\phm3,\phm3,\phm1,\phm1,\phm1,\phm1,\phm1)$,	&
				$(\phm0,\phm0,\phm0,\phm1,\phm1,\phm1,\phm1,\phm2)$. \\ \\
2-vertices	&	$(-3,-1,-1,-1,-1,-1,-1,-1)$,	& 
				$(\phm2,\phm2,\phm2,\phm2,\phm0,\phm0,\phm0,\phm0)$,\\
		&	$(\phm4,\phm0,\phm0,\phm0,\phm0,\phm0,\phm0,\phm0)$.\\\\
6-vertices	&	$(-3,-3,-3,-3,-3,-1,-1,-1)$,	&
				$(\phm6,\phm2,\phm2,\phm2,\phm0,\phm0,\phm0,\phm0)$,\\
		&	$(\phm4,\phm4,\phm4,\phm0,\phm0,\phm0,\phm0,\phm0)$,	&
				$(-5,\phm3,\phm3,\phm1,\phm1,\phm1,\phm1,\phm1)$,\\
		&	$(\phm4,\phm4,\phm2,\phm2,\phm2,\phm2,\phm0,\phm0)$. \\\\
7-vertices	&	$(-1,-1,-1,-1,-1,-1,\phm0,\phm0)$,	&
				$(\phm2,\phm1,\phm1,\phm0,\phm0,\phm0,\phm0,\phm0)$,\\
		&	$\half(-3,\phm3,\phm1,\phm1,\phm1,\phm1,\phm1,\phm1)$. \\\\
8-vertices	&	$(-1,-1,-1,-1,-1,-1,-1,\phm1)$,	&
				$(\phm2,\phm2,\phm0,\phm0,\phm0,\phm0,\phm0,\phm0)$.\\
\end{tabular}
\end{center}

This list can be verified by checking that the vertices listed indeed lie on the orbits 
$W_{E_8}\cdot v_i$ and there are as many as $|W_{E_8}|/|Stab_{W_{E_8}}(v_i)|$.

We describe in the following table
the 2- and 8-vertices $x$ in $\beta_2$. Let $\sigma$ be the face of
$\Si_{v_2} \triangle$ containing $\ora{v_2x}$ in its interior.

\begin{center}
\begin{tabular}{C{3cm} | C{6.5cm} | C{2.5cm} | C{2cm}}
&$x$		&	$d(x,v_2)$	&	Type of $\sigma$	\\\hline
2-vertices $x\neq v_2, \wht{v_2}$	&
	$(\phm1,-3,-1,-1,-1,-1,-1,\phm1)$\hfill		\linebreak
	$(\phm0,-2,-2,-2,-2,\phm0,\phm0,\phm0)$\hfill		\linebreak
	$(\phm1,-3,-1,-1,-1,-1,-1,\phm1)$\hfill		\linebreak
	$(\phm1,-1,-1,-1,-1,-1,-1,\phm3)$\hfill		\linebreak
	$(\phm2,-2,-2,-2,\phm0,\phm0,\phm0,\phm0)$\hfill		\linebreak
	$(\phm3,-1,-1,-1,-1,-1,-1,\phm1)$\hfill		\linebreak
	$(\phm3,-1,-1,-1,-1,\phm1,\phm1,\phm1)$\hfill		\linebreak
	$(\phm4,\phm0,\phm0,\phm0,\phm0,\phm0,\phm0,\phm0)$	&
	$\ac(\frac{3}{4})$\hfill					\linebreak
	$\pithird$\hfill						\linebreak
	$\ac(\quart)$\hfill						\linebreak
	$\pihalf$\hfill							\linebreak
	$\pihalf$\hfill							\linebreak
	$\ac(-\quart)$\hfill						\linebreak
	$\2pithird$\hfill						\linebreak
	$\ac(-\frac{3}{4})$						&
	3\hfill							\linebreak
	6\hfill							\linebreak
	38\hfill							\linebreak
	8\hfill							\linebreak
	5\hfill							\linebreak
	18\hfill							\linebreak
	6\hfill							\linebreak
	1	\\
\hline
8-vertices $x$	&
	$(-1,-1,-1,-1,-1,-1,-1,\phm1)$\hfill\linebreak
	$(\phm1,-1,-1,-1,-1,-1,-1,-1)$\hfill\linebreak
	$(\phm1,-1,-1,-1,-1,-1,\phm1,\phm1)$\hfill\linebreak
	$(\phm2,-2,\phm0,\phm0,\phm0,\phm0,\phm0,\phm0)$\hfill\linebreak
	$(\phm2,\phm0,\phm0,\phm0,\phm0,\phm0,\phm0,\phm2)$		&
	$\piquart$\hfill\linebreak
	$\ac(\frac{1}{2\sqrt{2}})$\hfill\linebreak
	$\pihalf$\hfill\linebreak
	$\ac(-\frac{1}{2\sqrt{2}})$\hfill\linebreak
	$\ac(\3piquart)$	&
	8\hfill\linebreak
	1\hfill\linebreak
	7\hfill\linebreak
	3\hfill\linebreak
	8	\\
\end{tabular}
\end{center}

We describe in the following table
the 7-vertices $x$ in $\beta_7$, 
such that $d(x,v_7)=\ac(-\third)$ or $\ac(-\frac{1}{6})$,
and the 8-vertices $x$ in $\beta_7$,
such that $d(x,v_7)>\pihalf$. 
Let $\sigma$ be the face of
$\Si_{v_7} \triangle$ containing $\ora{v_7x}$ in its interior.

\begin{center}
\begin{tabular}{C{3cm} | C{6.5cm} | C{2.5cm} | C{2cm}}
&$x$		&	$d(x,v_7)$	&	Type of $\sigma$	\\\hline
7-vertices $x$	&
	$(\phm0,\phm0,\phm0,\phm0,\phm0,\phm2,-1,-1)$\hfill		\linebreak
	$(\phm0,\phm0,\phm0,\phm0,\phm1,\phm1,-2,\phm0)$\hfill		\linebreak
	$\half(-1,\phm1,\phm1,\phm1,\phm1,\phm1,-3,-3)$\hfill		\linebreak
	$(\phm0,\phm0,\phm0,\phm0,\phm0,\phm1,-2,\phm1)$\hfill		\linebreak
	$\half(-3,\phm1,\phm1,\phm1,\phm1,\phm1,-3,-1)$\hfill		\linebreak
	$(\phm0,\phm0,\phm0,\phm0,\phm0,\phm1,-2,-1)$	&
	$\ac(-\third)$\hfill							\linebreak
	$\ac(-\third)$\hfill							\linebreak
	$\ac(-\third)$\hfill						\linebreak
	$\ac(-\frac{1}{6})$\hfill						\linebreak
	$\ac(-\frac{1}{6})$\hfill						\linebreak
	$\ac(-\frac{1}{6})$						&
	6\hfill							\linebreak
	58\hfill							\linebreak
	12\hfill							\linebreak
	68\hfill							\linebreak
	28\hfill							\linebreak
	168	\\
\hline
8-vertices $x$	&
	$(\phm1,\phm1,\phm1,\phm1,\phm1,\phm1,-1,\phm1)$\hfill\linebreak
	$(-1,\phm1,\phm1,\phm1,\phm1,\phm1,-1,-1)$\hfill\linebreak
	$(\phm0,\phm0,\phm0,\phm0,\phm0,\phm2,-2,\phm0)$		&
	$\ac(-\frac{\sqrt{3}}{2})$\hfill\linebreak
	$\ac(-\frac{1}{\sqrt{3}})$\hfill\linebreak
	$\ac(-\frac{1}{2\sqrt{3}})$	&
	8\hfill\linebreak
	2\hfill\linebreak
	68	\\
\end{tabular}
\end{center}

In order to make it easier to verify the table above, we present the complete table in
Appendix~\ref{app:E8}.

We want to describe the simplicial convex hull $C$ \label{app:7ptsconvhull}
 of the segment $v_7x$ for
the 7-vertex $x=(0,0,0,0,0,1,-2,-1)$, for this
we present 
first a larger 3-dimensional spherical polyhedron, namely the
tetrahedron $C':=CH(v_8,y,u_8,y')$,
where $y=(-1,-1,-1,-1,-1,-1,1,-1)$,
$u_8=(0,0,0,0,0,0,-2,-2)$ and $y'=(0,0,0,0,0,2,-2,0)$.
Notice that $v_7=m(v_8,y)$ and $x=m(y',u_8)$.
$C'$ is a subcomplex with four 2-dimensional faces: the triangles
$CH(v_8,y,y')$, $CH(z,y,y')$, $CH(y,u_8,v_8)$ and $CH(y',u_8,v_8)$.
The figures illustrate the tetrahedron $C'$ from the front and from behind.

\hpic{\includegraphics[scale=0.7]{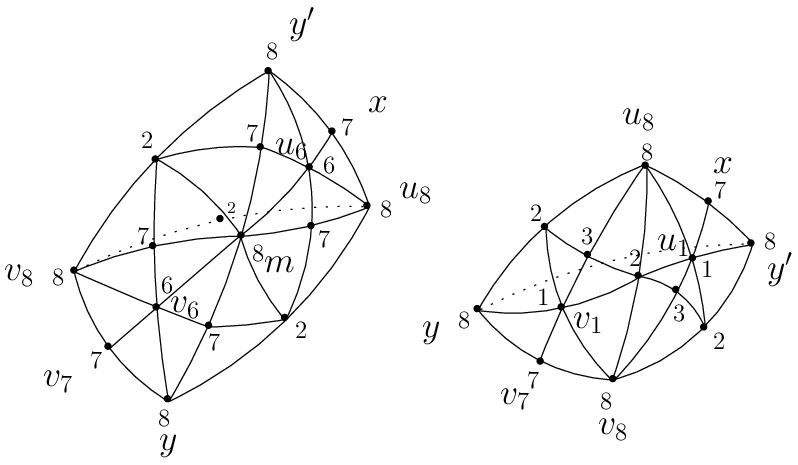}}\hspace{0.3cm}
\hpic{\parbox{7cm}{
$x=(\phm0,\phm0,\phm0,\phm0,\phm0,\phm1,-2,-1)$\\
$y=(-1,-1,-1,-1,-1,-1,\phm1,-1)$\\
$m=(-1,-1,-1,-1,-1,\phm1,-1,-1)$\\
$y'=(\phm0,\phm0,\phm0,\phm0,\phm0,\phm2,-2,\phm0)$\\
$u_1=(-1,-1,-1,-1,-1,\phm1,-5,-1)$\\
$u_6=(-1,-1,-1,-1,-1,\phm3,-5,-3)$\\
$u_8=(\phm0,\phm0,\phm0,\phm0,\phm0,\phm0,-2,-2)$\\
$m(u_8,v_8)=$\\
\hspace*{0.8cm}$(-1,-1,-1,-1,-1,-1,-3,-1)$
}}

\parpic[r]{
\hpic{\includegraphics[scale=0.65]{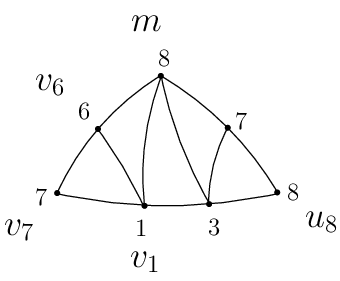}}\hspace{0.6cm}
\hpic{\includegraphics[scale=0.65]{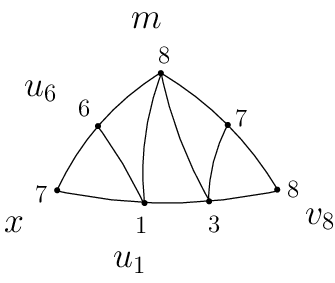}}}
The triangles $CH(v_8,m,x)$ and $CH(v_7,m,u_8)$ are 2-dimensional subcomplexes.
If we cut $C'$ along these triangles, we obtain a convex subcomplex
$C'':=CH(v_7,v_8,x,u_8,m)$. 
It has six 2-dimensional faces: the triangles\linebreak[4]
$CH(m,v_7,u_8)$, $CH(m,x,v_8)$, $CH(m,v_7,v_8)$,\linebreak[4] $CH(m,x,u_8)$, 
$CH(v_7,v_8,u_8)$ and $CH(x,u_8,v_8)$.
Recall that the direction $\ora{v_7x}$ spans the 168-face in $\Si_{v_7}\triangle$,
this implies that $v_1$, $v_6$ and $v_8$ 
are contained in the simplicial convex hull $C$ of $v_7x$. 
We can also see that the direction $\ora{xv_7}$ spans the 168-face with vertices
$\ora{xu_1}$, $\ora{xu_6}$ and $\ora{xu_8}$. In particular, $u_8\in C$. Considering
the triangle $CH(v_7,m,u_8)$ we deduce that also $m\in C$.
It follows that $C=C''$.
The next figure shows the link $\Si_m C'$.

\begin{center}\includegraphics[scale=0.5]{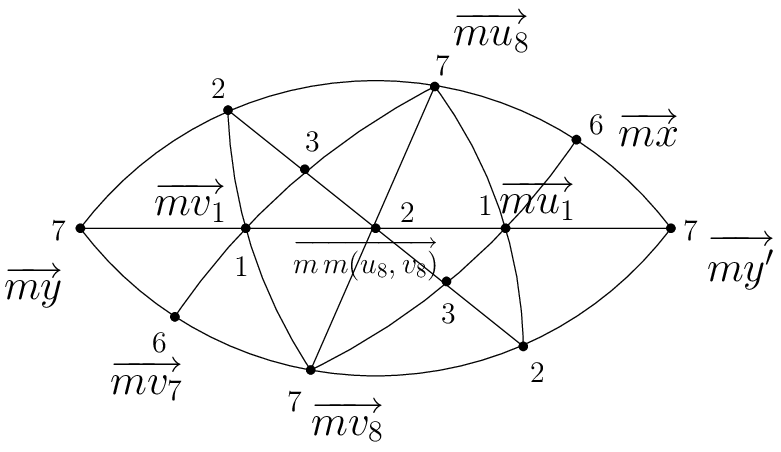}\end{center}

We describe in the following table
the 8-vertices $x$ in $\beta_8$. Let $\sigma$ be the face of
$\Si_{v_8} \triangle$ containing $\ora{v_8x}$ in its interior.
\begin{center}
\begin{tabular}{C{3cm} | C{6.5cm} | C{2.5cm} | C{2cm}}
&$x$		&	$d(x,v_8)$	&	Type of $\sigma$	\\\hline
8-vertices $x\neq v_8,\wht{v_8}$	&
	$(-1,-1,-1,-1,-1,-1,\phm1,-1)$\hfill		\linebreak
	$(-2,\phm0,\phm0,\phm0,\phm0,\phm0,\phm0,-2)$\hfill		\linebreak
	$(\phm0,\phm0,\phm0,\phm0,\phm0,\phm0,\phm2,-2)$	&
	$\pithird$\hfill					\linebreak
	$\pihalf$\hfill					\linebreak
	$\2pithird$						&
	7\hfill							\linebreak
	2\hfill							\linebreak
	7	\\
\end{tabular}
\end{center}

\label{page:7vertices}
We describe in the following table
the 7-vertices $x$ in $\beta_7(2,8)$ with
$d(x,v_7)>\pihalf$. Let $\sigma$ be the face of
$\Si_{v_7} \triangle(2,8)$ containing $\ora{v_7x}$ in its interior.
\begin{center}
\begin{tabular}{C{3cm} | C{6.5cm} | C{2.5cm} | C{2cm}}
&$x$		&	$d(x,v_7)$	&	Type of $\sigma$	\\\hline
7-vertices $x\neq v_7,\wht{v_7}$	&
	$(\phm0,\phm0,\phm1,\phm1,\phm1,\phm1,-1,-1)$\hfill		\linebreak
	$(\phm0,\phm0,\phm0,\phm0,\phm0,\phm2,-1,-1)$	&
	$\ac(-\frac{2}{3})$\hfill					\linebreak
	$\ac(-\third)$						&
	3\hfill							\linebreak
	6	\\
\end{tabular}
\end{center}

In order to make it easier to verify the table above, we present the complete table in
Appendix~\ref{app:E8}.

We describe in the following table
the 1-vertices $x$ in $\beta_1(2,7,8)$ with
$d(x,v_1)>\pihalf$. Let $\sigma$ be the face of
$\Si_{v_1} \triangle(2,7,8)$ containing $\ora{v_1x}$ in its interior.
\begin{center}
\begin{tabular}{C{3cm} | C{6.5cm} | C{2.5cm} | C{2cm}}
&$x$		&	$d(x,v_1)$	&	Type of $\sigma$	\\\hline
1-vertices $x\neq v_1,\wht{v_1}$	&
	$\half(-1,-1,-1,-1,\phm1,\phm3,\phm3,\phm3)$\hfill		\linebreak
	$(-1,-1,\phm1,\phm1,\phm1,\phm1,\phm1,\phm1)$\hfill		\linebreak
	$\half(-1,-1,\phm1,\phm1,\phm1,\phm3,\phm3,\phm3)$\hfill		\linebreak
	$\half(\phm1,\phm1,\phm1,\phm1,\phm1,\phm3,\phm3,\phm3)$	&
	$\ac(-\frac{3}{8})$\hfill					\linebreak
	$\2pithird$\hfill					\linebreak
	$\ac(-\frac{5}{8}) $\hfill					\linebreak
	$\ac(-\frac{7}{8})$						&
	56\hfill							\linebreak
	3\hfill							\linebreak
	36\hfill							\linebreak
	6	\\
\end{tabular}
\end{center}

In order to make it easier to verify the table above, we present the complete table in
Appendix~\ref{app:E8}.

We describe in the following table
the 6-vertices $x$ in $\beta_6(1,2,7,8)$ with
$d(x,v_6)>\pihalf$. Let $\sigma$ be the face of
$\Si_{v_6} \triangle(1,2,7,8)$ containing $\ora{v_6x}$ in its interior.
\begin{center}
\begin{tabular}{C{3cm} | C{6.5cm} | C{2.5cm} | C{2cm}}
&$x$		&	$d(x,v_6)$	&	Type of $\sigma$	\\\hline
6-vertices $x\neq v_6,\wht{v_6}$	&
	$(\phm0,\phm0,\phm0,\phm0,\phm6,-2,-2,-2)$\hfill		\linebreak
	$(\phm0,\phm0,\phm2,\phm4,\phm4,-2,-2,-2)$\hfill		\linebreak
	$(\phm1,\phm1,\phm3,\phm3,\phm5,-1,-1,-1)$	&
	$\ac(-\quart)$\hfill					\linebreak
	$\2pithird$\hfill					\linebreak
	$\ac(-\frac{3}{4}) $						&
	5\hfill							\linebreak
	34\hfill							\linebreak
	35	\\
\end{tabular}
\end{center}
 
Let us verify this last table. By considering the following 2-dimensional
bigons, we can see that if there are 6-vertices missing in the table above, they
must lie in the interior of $\beta_6(1,2,7,8)$.
\begin{center}
\hpic{\includegraphics[scale=0.7]{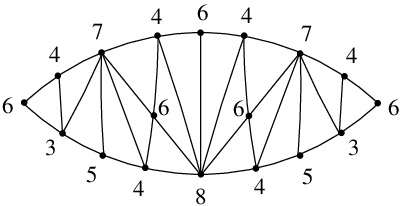}}\hspace{1cm}
\hpic{\includegraphics[scale=0.7]{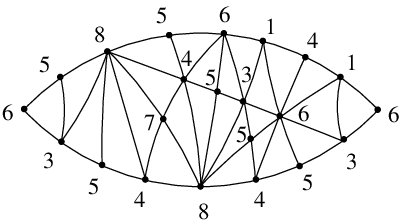}}\hspace{1cm}
\hpic{\includegraphics[scale=0.7]{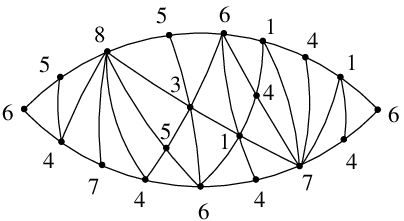}}
\end{center}
A 6-vertex $x$ in the interior of $\beta_6(1,2,7,8)$
should satisfy 
$$
 x_4+x_5+\dots+x_8 \overset{\text{(1)}}{=} x_1+x_2+x_3\;;\;\;\; \;\;
 x_1\overset{\text{(2)}}{=} x_2
 \overset{\text{(3)}}{<}x_3
 \overset{\text{(4)}}{<}x_4
 \overset{\text{(5)}}{<}x_5\;;\;\;\; \;\;
x_5 > x_6
 \overset{\text{(7)}}{=}x_7
 \overset{\text{(8)}}{=}x_8.
$$
In particular, we have four different values $x_2<x_3<x_4<x_5$. Hence,
$x$ cannot be a permutation of $(\pm4,\pm4,\pm4,0,0,0,0,0)$.

If $x$ is obtained from $(-3,-3,-3,-3,-3,-1,-1,-1)$ by
permutations of the coordinates and change of sign in an even number 
of places, then $x_1=x_2=-3$, $x_3=-1$, $x_4=1$ and $x_5=3$.
By the equalities (1), (2), (7) and (8), it follows
that $x_6=-\frac{11}{3}$, which is not possible.

If $x$ is obtained from $(\pm6,\pm2,\pm2,\pm2,0,0,0,0)$ by
permutations of the coordinates, then 
$(x_2,x_3,x_4,x_5)=(-6,-2,0,2)$, but $x_1=x_2=-6$ is not possible;
or $(x_2,x_3,x_4,x_5)=(-2,0,2,6)$.
In this case the equalities (1), (2), (7) and (8) imply 
$x_6=-4$, which is not possible.

If $x$ is obtained from $(-5,3,3,1,1,1,1,1)$ by
permutations of the coordinates and change of sign in an even number 
of places, then $x_2\in\{-5,-3,-1\}$. And $x_1=x_2=-5$ is not possible.
$x_2=-3$ implies $(x_1,x_2,x_3,x_4,x_5)=(-3,-3,-1,1,5)$ and 
equalities (1), (2), (7) and (8) imply $x_6=-\frac{13}{3}$. This is again
impossible.
$x_2=-1$ implies $(x_1,x_2,x_3,x_4,x_5)=(-1,-1,1,3,5)$ and 
equalities (1), (2), (7) and (8) imply $x_6=x_7=x_8=-3$, which cannot happen.

If $x$ is obtained from $(\pm4,\pm4,\pm2,\pm2,\pm2,\pm2,0,0)$ by
permutations of the coordinates, then
$(x_1,x_2,x_3,x_4,x_5)=(-4,-4,-2,0,2)$ or $(-2,-2,0,2,4)$. In both cases
the equalities $x_6=x_7=x_8$ cannot be satisfied.

So we have verified that $\beta_6(1,2,7,8)$ contains no 6-vertices
in its interior and therefore our table is complete.

\section{More information about $E_8$}\label{app:E8}

In this section, we complete some tables 
given in Appendix~\ref{app:E8coxeter}.
Although this information is not directly used in the proof of our main result,
we present it here in order to make it easier 
to verify the tables in Appendix~\ref{app:E8coxeter}.

The next table lists the 7-vertices $x$ in $\beta_7$ with $d(x,v_7)\geq\pihalf$.
The vertices marked with * are the
ones at distance $=\pihalf$ to $v_7$. 
Let $\sigma$ be the face of
$\Si_{v_7} \triangle$ containing $\ora{v_7x}$ in its interior.
Let $\sigma_x$ be the face of $\triangle$ spanned by
the initial part of the segment $v_7x$.

\begin{center}
\begin{tabular}{ R{6.5cm} | C{1.9cm} | C{7.7cm}}
$x$\hspace{3cm}\phm	&	Type of $\sigma$ & 
	$|Stab_{W_{E_8}}(v_7)\cdot x|=
		\frac{|Stab_{W_{E_8}}(v_7)|}{|Stab_{W_{E_8}}(\sigma_x)|}$	\\\hline
$(\phm1,\phm1,\phm1,\phm1,\phm1,\phm1,\phm0,\phm0) $ & &
	$|W_{E_6}||W_{A_1}|/(|W_{E_6}||W_{A_1}|)=1$\\
$(\phm0,\phm0,\phm1,\phm1,\phm1,\phm1,-1,-1) $ & 3 &
	$|W_{E_6}||W_{A_1}|/(|W_{A_1}||W_{A_4}||W_{A_1}|)=216$\\
* $\half(-1,-1,-1,\phm1,\phm1,\phm1,-3,-3) $ & 4&
	$|W_{E_6}||W_{A_1}|/(|W_{A_2}||W_{A_1}||W_{A_2}||W_{A_1}|)=\nolinebreak[4]720$\\	
$(\phm0,\phm0,\phm0,\phm0,\phm0,\phm2,-1,-1)$ & 6&
	$|W_{E_6}||W_{A_1}|/(|W_{D_5}||W_{A_1}|)=27$\\
$\half(\phm1,\phm1,\phm1,\phm1,\phm1,\phm1,-3,\phm3)$ & 8&
	$|W_{E_6}||W_{A_1}|/|W_{E_6}|=2$\\
$\half(-1,\phm1,\phm1,\phm1,\phm1,\phm1,-3,-3)$ & 12&
	$|W_{E_6}||W_{A_1}|/(|W_{A_4}||W_{A_1}|)=432$\\
$(\phm0,\phm1,\phm1,\phm1,\phm1,\phm1,-1,\phm0)$ & 28&
	$|W_{E_6}||W_{A_1}|/|W_{D_5}|=54$\\
$\half(-3,\phm1,\phm1,\phm1,\phm1,\phm1,-3,-1)  $ & 28&
	$|W_{E_6}||W_{A_1}|/|W_{D_5}|=54$\\
$(\phm0,\phm0,\phm0,\phm0,\phm1,\phm1,-2,\phm0)$ & 58&
	$|W_{E_6}||W_{A_1}|/(|W_{A_4}||W_{A_1}|)=432$\\
$(\phm0,\phm0,\phm0,\phm0,\phm0,\phm1,-2,\phm1)$ & 68&
	$|W_{E_6}||W_{A_1}|/|W_{D_5}|=54$\\
$\half(\phm1,\phm1,\phm1,\phm1,\phm1,\phm3,-3,\phm1) $ & 68&
	$|W_{E_6}||W_{A_1}|/|W_{D_5}|=54$\\
$(\phm0,\phm0,\phm0,\phm0,\phm0,\phm1,-2,-1)$ & 168&
	$|W_{E_6}||W_{A_1}|/|W_{A_4}|=864$\\
* $\phantom{\half}(-1,\phm0,\phm0,\phm0,\phm0,\phm1,-2,\phm0)$  & 268&
	$|W_{E_6}||W_{A_1}|/|W_{D_4}|=540$\\
$\half(-1,\phm1,\phm1,\phm1,\phm1,\phm3,-3,-1)$ & 268&
	$|W_{E_6}||W_{A_1}|/|W_{D_4}|=540$\\
\end{tabular}
\end{center}

Notice that since the antipode $\wht v_7$ of $v_7$ is also a 7-vertex, then 
the number of 7-vertices in $S$ at distance $\leq\pihalf$ to $v_7$
is the same as the
number of 7-vertices in $S$ at distance $\geq\pihalf$ to $v_7$. 
It follows that
the number of 7-vertices in $S$ is two times the number 
of 7-vertices in $S$ at distance $\leq\pihalf$ to $v_7$ 
minus the number of 7-vertices at distance $=\pihalf$ to $v_7$.
With this observation and the one at the end of the introductory section
of Appendix~\ref{app:coxeter},
we can verify the correctness of the list above:
$2(1+216+720+27+2+432+54+54+432+54+54+864+540+540)-720-540
=6720=\frac{|W_{E_8}|}{|Stab_{W_{E_8}}(v_7)|}=\#\{7\text{- vertices in }S\}$.

The next table lists the 8-vertices $x$ in $\beta_7$ with $d(x,v_7)\geq\pihalf$.
The vertices marked with * are the
ones at distance $=\pihalf$ to $v_7$. 
Let $\sigma$ be the face of
$\Si_{v_7} \triangle$ containing $\ora{v_7x}$ in its interior. 

\begin{center}
\begin{tabular}{ R{6.7cm} | C{2.5cm} }
$x$\hspace{3cm}\phm	&	Type of $\sigma$	\\\hline
* $(\phm0,\phm0,\phm0,\phm0,\phm0,\phm0,-2,-2)$ & 1 \\
$(-1,\phm1,\phm1,\phm1,\phm1,\phm1,-1,-1)$ & 2 \\
$(\phm1,\phm1,\phm1,\phm1,\phm1,\phm1,-1,\phm1)$ & 8 \\
* $(\phm0,\phm0,\phm0,\phm0,\phm0,\phm0,-2,\phm2)$ & 8 \\
$(\phm0,\phm0,\phm0,\phm0,\phm0,\phm2,-2,\phm0)$ & 68 \\
\end{tabular}
\end{center}

The next table lists the 1-vertices $x$ in $\beta_1$ with $d(x,v_1)\geq\pihalf$.
The vertices marked with * are the
ones at distance $=\pihalf$ to $v_1$. 
Let $\sigma$ be the face of
$\Si_{v_1} \triangle$ containing $\ora{v_1x}$ in its interior. 

\begin{center}
\begin{tabular}{ R{6.5cm} | C{2cm} | C{7.5cm}}
$x$\hspace{3cm}\phm	&	Type of $\sigma$ & 
	$|Stab_{W_{E_8}}(v_1)\cdot x|=
		\frac{|Stab_{W_{E_8}}(v_1)|}{|Stab_{W_{E_8}}(\sigma_x)|}$	\\\hline
$(\phm1,\phm1,\phm1,\phm1,\phm1,\phm1,\phm1,\phm1)$ & & 
	$|W_{A_7}|/|W_{A_7}|=1$\\
$\half(-5,\phm1,\phm1,\phm1,\phm1,\phm1,\phm1,\phm1)$ &   2&	
	$|W_{A_7}|/|W_{A_6}|=8$\\
$(-1,-1,\phm1,\phm1,\phm1,\phm1,\phm1,\phm1)$   & 3 &
	$|W_{A_7}|/(|W_{A_1}||W_{A_5}|)=28$\\
*\hspace{0.2cm} $(-1,-1,-1,-1,\phm1,\phm1,\phm1,\phm1)$ &   5& 
	$|W_{A_7}|/(|W_{A_3}||W_{A_3}|)=70$\\
$\half(\phm1,\phm1,\phm1,\phm1,\phm1,\phm3,\phm3,\phm3)$ &   6&
	$|W_{A_7}|/(|W_{A_2}||W_{A_4}|)=56$\\
$(-2,\phm0,\phm0,\phm0,\phm1,\phm1,\phm1,\phm1)$ &   25&
	$|W_{A_7}|/(|W_{A_2}||W_{A_3}|)=280$\\
$\half(-1,\phm1,\phm1,\phm1,\phm1,\phm1,\phm1,\phm5)$ &   28&
	$|W_{A_7}|/|W_{A_5}|=56$\\
$\half(-1,-1,\phm1,\phm1,\phm1,\phm3,\phm3,\phm3)$  &  36&
	$|W_{A_7}|/(|W_{A_1}||W_{A_2}||W_{A_2}|)=560$\\
$\half(-3,-3,\phm1,\phm1,\phm1,\phm1,\phm1,\phm3)$ &  38 & 
	$|W_{A_7}|/(|W_{A_1}||W_{A_4}|)=168$\\
$(\phm0,\phm0,\phm0,\phm1,\phm1,\phm1,\phm1,\phm2)$ &   48& 
	$|W_{A_7}|/(|W_{A_2}||W_{A_3}|)=280$\\
$\half(-1,-1,-1,\phm1,\phm1,\phm1,\phm1,\phm5)$  &  48&
	$|W_{A_7}|/(|W_{A_2}||W_{A_3}|)=280$\\
$\half(-1,-1,-1,-1,\phm1,\phm3,\phm3,\phm3)$  &  56&
	$|W_{A_7}|/(|W_{A_2}||W_{A_3}|)=280$\\
$\half(-1,-1,-1,-1,-1,\phm1,\phm1,\phm5)$  &  68&
	$|W_{A_7}|/(|W_{A_4}||W_{A_1}|)=168$\\
*\hspace{0.2cm} $(-2,-1,\phm0,\phm0,\phm0,\phm1,\phm1,\phm1)$  &  236&
	$|W_{A_7}|/(|W_{A_2}||W_{A_2}|)=1120$\\
$\half(-3,-1,\phm1,\phm1,\phm1,\phm1,\phm3,\phm3)$ &   237&
	$|W_{A_7}|/(|W_{A_3}||W_{A_1}|)=840$\\
$\half(-3,-1,-1,-1,\phm1,\phm1,\phm3,\phm3)$  &  257&
	$|W_{A_7}|/(|W_{A_2}||W_{A_1}||W_{A_1}|)=1680$\\
$(-1,\phm0,\phm0,\phm0,\phm1,\phm1,\phm1,\phm2)$ &   258 &
	$|W_{A_7}|/(|W_{A_2}||W_{A_2}|)=1120$\\
$(-1,-1,\phm0,\phm0,\phm0,\phm1,\phm1,\phm2)$  &  368  &
	$|W_{A_7}|/(|W_{A_2}||W_{A_1}||W_{A_1}|)=1680$\\
*\hspace{0.2cm} $(-1,-1,-1,\phm0,\phm0,\phm0,\phm1,\phm2)$  &  478&
	$|W_{A_7}|/(|W_{A_2}||W_{A_2}|)=1120$\\
\end{tabular}
\end{center}

We can verify this table as we did with the table above:
$2(1+8+28+70+56+280+56+560+168+280+280+280+168+1120+840+1680+1120+
1680+1120)-70-1120-1120=17280=\frac{|W_{E_8}|}{|W_{A_7}|}=\#\{1\text{- vertices in }S\}$.

\bigskip\bigskip
\textbf{{\em Acknowledgments.} }
The results in this paper are part of my doctoral thesis.
I am very thankful to my advisor Prof. Bernhard Leeb for his support
and encouragement, for introducing me to this problem 
in the joint work \cite{LeebRamos}, and for many useful
comments and corrections on previous versions of this paper.
I gratefully acknowledge the financial support 
from the joint program between the
Consejo Nacional de Ciencia y Tecnolog\'ia (Mexico) and the 
Deutscher Akademischer Austauschdienst, also from the
Deutsche Forschungsgemeinschaft and the Hausdorff Research Institute
for Mathematics.
The author is also grateful to the referee 
for useful suggestions concerning the presentation of this paper.


\begin{thebibliography}{}

\bibitem{AbramenkoBrown}
P.\ Abramenko, K.\ S.\ Brown,
{\em Buildings: Theory and Applications},
GTM 248, Springer 2008.

\bibitem{BalserLytchak}
A.\ Balser, A.\ Lytchak,
{\em Centers of convex subsets of buildings}, 
Ann.\ Glob.\ Anal.\ Geom.\ 28, No.\ 2, 201-209 (2005).

\bibitem{BateMartinRoehrle}
M. Bate, B. Martin, G. R\"ohrle,
{\em On Tits' Centre Conjecture for fixed point subcomplexes},
C.\ R.\ Math.\ Acad.\ Sci.\ Paris 347, no.\ 7-8, 353-356 (2009).

\bibitem{BateMartinRoehrle3}
M. Bate, B. Martin, G. R\"ohrle,
{\em Complete reducibility and separable field extensions},
C.\ R.\ Math.\ Acad.\ Sci.\ Paris 348, no.\ 9-10,  495-497 (2010).

\bibitem{BateMartinRoehrle2}
M. Bate, B. Martin, G. R\"ohrle,
{\em The strong Centre Conjecture: an invariant theory approach},
J.\ Algebra, to appear. Preprint. arXiv:1005.3212v3.

\bibitem{BateMartinRoehrleTange}
M.\ Bate, B.\ Martin, G.\ Roehrle, R.\ Tange,
{\em Closed Orbits and uniform S-instability in Geometric Invariant Theory},
Trans.\ Amer.\ Math.\ Soc., to appear. Preprint. arXiv:0904.4853v2.

\bibitem{BridsonHaefliger}
M.\ Bridson, A.\ Haefliger, 
{\em Metric spaces of non-positive curvature}, 
Springer 1999.

\bibitem{GroveBenson}
L.C.\ Grove, C.T.\ Benson,
{\em Finite reflection groups}, 
GTM 99, Springer 1971.

\bibitem{Kempf}
G.\ R.\ Kempf,
{\em Instability in invariant theory},
Ann. of Math. 108, 299-316 (1978).

\bibitem{KleinerLeeb}
B.\ Kleiner, B.\ Leeb, 
{\em Rigidity of quasi-isometries for symmetric spaces 
and Euclidean buildings}, 
Inst.\ Hautes \'Etudes Sci. Publ.\ Math.\ No.\ 86 (1997), 115--197 (1998).

\bibitem{KleinerLeebinvconv}
B.\ Kleiner, B.\ Leeb, 
{\em Rigidity of invariant convex sets in symmetric spaces}, 
Invent. Math. 163, No.\ 3, 657-676 (2006).

\bibitem{LeebRamos}
B.\ Leeb, C.\ Ramos-Cuevas,
{\em The center conjecture for spherical buildings of types $F_4$ and $E_6$}, 
Geom. Funct. Anal. 21, no.\ 3, 525-559 (2011). 

\bibitem{MuehlherrTits}
B.\ M\"uhlherr, J.\ Tits,
{\em The center conjecture for non-exceptional buildings}, 
J.\ Algebra 300, No.\ 2, 687-706 (2006).

\bibitem{MuehlherrWeiss}
B.\ M\"uhlherr, R.\ Weiss,
{\em Receding polar regions of a spherical building and the center conjecture}, 
Ann. Inst. Fourier, to appear. 

\bibitem{Mumford}
D.\ Mumford,
{\em Geometric Invariant Theory},
Springer, 1965.

\bibitem{ParkerTent_conv}
C.\ Parker, K.\ Tent,
{\em Convexity in buildings}, 
talk at the Oberwolfach conference 
{\em Buildings: Interactions with algebra and geometry},
Oberwolfach report no.\ 3/2008, \\
http://www.mfo.de/programme/schedule/2008/04/OWR\underline{\phantom{a}}2008\underline{\phantom{a}}03.pdf

\bibitem{Ramos}
C.\ Ramos-Cuevas,
{\em On convex subcomplexes of spherical buildings and Tits' Center Conjecture},
Ph.D. Thesis, University of Munich (2009).

\bibitem{Rousseau}
G.\ Rousseau,
{\em Immeubles sph\'eriques et th\'eorie des invariants},
C. R. Acad. Sci. Paris Ser. I 286, 347-250 (1978).

\bibitem{Scharlau}
R.\ Scharlau,
{\em A structure theorem for weak buildings of spherical type},
Geom. Dedicata 24, 77-84 (1987).

\bibitem{Serre}
J.-P.\ Serre,
{\em Compl\`ete r\'eductibilit\'e},
S\'em.\ Bourbaki, exp.\ 932, 
Ast\'erisque 299 (2005). 

\bibitem{Struyve}
K.\ Struyve,
{\em (Non)-completeness of $\R$-buildings and fixed point theorems},
Groups Geom.\ Dyn.\ 5, no.\ 1, 177-188 (2011). 

\bibitem{Titscenterconj}
J.\ Tits,
{\em Groupes semi-simples isotropes},
Coll. sur la th\'eorie des groupes alg\'ebriques, Bruxelles, 137-147 (1962).

\bibitem{Tits_bn}
J.\ Tits,
{\em Buildings of spherical type and finite BN-pairs},
LNM 386, Springer 1974.

\bibitem{Tits_w}
J.\ Tits,
{\em Endliche Spiegelungsgruppen, die als Weylgruppen auftreten}, 
Invent.\ Math.\ 43, 283-295 (1977).

\end{thebibliography}
\end{document}